\documentclass[a4paper,11pt]{article}

\usepackage[numbers]{natbib}

\usepackage{amssymb, adjustbox, amsfonts, amsmath, amsthm, bbm, caption, comment, empheq, mathrsfs, mathtools, tikz, tikz-3dplot, enumitem}
\usetikzlibrary{arrows.meta, decorations.markings}
\tikzset{midarrow/.style={
postaction={decorate, decoration={markings, mark=at position 0.5 with {\arrow[scale=1.2]{>}}}}}}

\newcommand{\uc}{\underline{c}}
\newcommand{\uw}{\underline{w}}
\newcommand{\ur}{\underline{\rho}}
\newcommand{\ov}{\overline{v}}
\newcommand{\oc}{\overline{c}}
\newcommand{\oV}{\overline{V}}
\newcommand{\Lip}{{\rm Lip}}
\newcommand{\vref}{v_{\rm ref}}
\newcommand{\Rref}{R_{\rm ref}}
\newcommand{\bU}{\mathbf{U}}
\newcommand{\bW}{\mathbf{W}}
\newcommand{\tbW}{\widetilde{\mathbf{W}}}

\newcommand{\ts}{\tilde{s}}
\newcommand{\tTV}{\widetilde{\rm{TV}}}

\newcommand{\calU}{\mathcal{U}}
\newcommand{\calI}{\mathcal{I}}
\newcommand{\R}{\mathbb R}
\newcommand{\N}{\mathbb N}
\newcommand{\Z}{\mathbb Z}
\newcommand{\D}{\mathbb{D}_0}
\newcommand{\Dr}{\mathbb{D}}
\newcommand{\PC}{\mathbf{PC}}
\newcommand{\Lloc}[1]{\mathbf{L^{#1}_{loc}}}
\newcommand{\BV}{\mathbf{BV}}
\newcommand{\C}[1]{\mathbf{C^{#1}}}
\newcommand{\Cc}[1]{\mathbf{C_c^{#1}}}
\newcommand{\TV}{{\rm{TV}}}
\newcommand{\RS}{\mathsf{RS}}

\newtheorem{theorem}{Theorem}
\newtheorem{lemma}[theorem]{Lemma}
\newtheorem{proposition}[theorem]{Proposition}
\newtheorem{example}[theorem]{Example}
\newtheorem{corollary}[theorem]{Corollary}
\newtheorem{definition}[theorem]{Definition}
\newtheorem{remark}[theorem]{Remark}

\usepackage[hidelinks]{hyperref}
\usepackage[total={190mm,257mm},left=10mm,top=20mm]{geometry}

\renewcommand{\geq}{\geqslant}
\renewcommand{\leq}{\leqslant}
\renewcommand{\d}{{\rm{d}}}
\renewcommand{\L}[1]{\mathbf{L^{#1}}}
\begin{document}
\let\WriteBookmarks\relax
\def\floatpagepagefraction{1}
\def\textpagefraction{.001}

\title{Existence and uniqueness result for a $2 \times 2$ system of conservation\\laws with discontinuous flux and application to vehicular traffic} 

\date{}

\author{
Felisia Angela Chiarello\\
Department of Engineering\\
University LUM, S.S.~100 km 18, Casamassima, 70010, Italy\\
chiarello@lum.it
\and
Simone Fagioli\\
Department of Engineering and Information Science and Mathematics,\\
University of L’Aquila, via Vetoio, Ed.~Coppito~1, Coppito, 67100, Italy\\
simone.fagioli@univaq.it
\and 
Massimiliano Daniele Rosini\\
Faculty of Mathematics and Computer Science,\\
Maria Curie-Skłodowska University,
Plac Marii Curie-Skłodowskiej 1, Lublin, 20031, Poland,
\\
Department of Management and Business Administration,\\ University \lq\lq G.~d'Annunzio\rq\rq\ of Chieti-Pescara, viale Pindaro, 42, Pescara, 65127, Italy\\
massimiliano.rosini@unich.it
}

\maketitle

\begin{abstract}
We study the Cauchy problem for a one-dimensional $2\times2$ non-symmetric Keyfitz-Kranzer system of conservation laws with a flux $f=f(x,\bU)$ that is discontinuous with respect to the spatial variable $x$.
No monotonicity assumption is imposed on the mapping $x\mapsto f(x,\bU)$.
We introduce a Kruzhkov-type entropy condition and establish global existence and uniqueness of entropy solutions for initial data with arbitrarily large total variation.
The proof relies on a wave-front tracking approximation.
The main technical novelty consists in the introduction of \emph{adapted} Riemann invariant coordinates, specifically designed to account for the discontinuities of the flux, which yield a uniform-in-time bound on the total variation of the approximate solutions $\bU^n(t)$.
As an application, we discuss the second-order ARZ model for vehicular traffic on inhomogeneous roads with abrupt \lq\lq collective\rq\rq\ changes in the speed law or road capacity.

\end{abstract}

\section{Introduction}

In this paper, we consider the Cauchy problem for a one dimensional $2\times2$ system of Temple type with discontinuous flux.
Temple systems are a particular class of strictly hyperbolic systems of conservation laws characterised by the fact that, for each genuinely non-linear characteristic field, shock and rarefaction curves coincide, see \cite{Temple02, Temple01}. 

Existence and well-posedness results for Temple systems have been established in several works, including \cite{BaitiBressan-Temple, Bianchini-Temple, BressanGoatin-Temple, ColomboCorli-Temple, Serre-Temple}.
These analyses have been extended to initial-boundary value problems, see for instance \cite{AnconaGoatin-Temple, ColomboRosini-Temple01, ColomboRosini-Temple02}. 
The vanishing viscosity approximation has also been studied, see for instance \cite{BianchiniBressan2002, HaspotJana2025, Serre1987}.
On the other hand, numerical approximation is still a delicate issue.
Indeed, a well-known problem associated with numerical approximation of Temple systems is that Godunov's and related methods produce spurious oscillations near contact discontinuities, since the numerical solution may leave the non-convex invariant region of the exact solution.
This leads to the introduction of \emph{ad hoc} numerical schemes.
For instance, in \cite{ChalonsGoatin2007} the authors propose a random sampling strategy combined with the Godunov method.
This approach was subsequently extended in \cite{Betancourt2017}.
A different strategy is proposed in \cite{DonadelloPolizziRazafisonRollandRosini-2026}, where the authors develop an adapted version of the Glimm scheme.

More recently, in \cite{BorscheGaravelloGunarso-Temple} the authors considered the Cauchy problem on networks for hyperbolic balance laws of Temple type under the assumption that the nodes are \emph{non-characteristic}, i.e., when the characteristic speeds are uniformly bounded away from zero.
Their framework also includes Temple systems with discontinuous fluxes.
Indeed, as noted in \cite[Remark 2.1]{BorscheGaravelloGunarso-Temple}, the flux function can be chosen differently on the various edges of the network; in particular, by considering a 1-to-1 junction one obtains a Temple system with discontinuous flux.
Another reference on Temple systems with discontinuous flux is \cite{RosiniAnnales}.
In that work, the author proposes six generalisations of the homogeneous Aw-Rascle-Zhang model \cite{AwRascle2000, zhang2002TRB}, a $2 \times 2$ system of Temple type, designed to account for road sections with different capacities and/or speed limits, which naturally lead to discontinuous fluxes.

The study of hyperbolic systems of conservation laws of the form
\[\partial_t\bU+\partial_xf(t,x,\bU)=0\]
with discontinuous flux function $f$ is a topic of intense current research, as such models arise in engineering applications and the applied sciences, see for instance the survey papers \cite{BurgerKarlsen2008, Mishra2017} and the references therein.
In fact, the flux $f$ may depend discontinuously on the unknown variable $\bU$ as in \cite{BurgerChalonsOrdonezVillada2021}, on the independent variable $(t,x)$ as in \cite{CrastaDeCiccoDePhilippisGhiraldin2016}, on the spatial position $x$ as in \cite{KarlsenMitrovic2025}, at some fixed space locations as in \cite{BandaHertyKlar2006, BurgerDiehlMarti}, or at moving locations that behave as free interfaces, see \cite{SunQuYuan2024}.
For numerical approximation schemes and their analysis in the presence of discontinuous fluxes, we refer the reader to
\cite{BaleLevequeMitranRossmanith2002, BurgerDiehlMartiVasquez2023, BurgerGarciaKarlsenTowers2008, BurgerKarlsenTowers2010, QiaoZhangLinWongChoi2017, ZhangWongShu2006}.

The main result of the present paper establishes the global existence and uniqueness of entropy solutions $\bU=(\rho,q)$ of the Cauchy problem for a one-dimensional $2\times2$ non-symmetric Keyﬁtz-Kranzer type system of conservation laws and a time-independent flux $f=f(x,\rho,q)$ of the form
\begin{equation}
\label{e:flux}
f(x,\rho,q) = c(x) \, V\left(\frac{q}{\rho} - p(\rho)\right) \bU,
\end{equation}
where the space-dependent coefficient function $c$ is piecewise constant and discontinuous at finitely many fixed points, and for initial data with arbitrarily large total variation.
We underline that our result does not rely on any monotonicity assumption on $f(x,\bU)$ with respect to $x$.
Moreover, unlike the non-characteristic setting, the second characteristic speed is allowed to vanish.
The proof requires some techniques typical of hyperbolic systems of conservation laws. 
To prove the existence of a solution, we use the wave-front tracking method to construct an approximate solution $\bU^n$, see \cite{Bressan-book, Holden2015} and the references therein. 
This method relies primarily on Riemann problems and their solutions as fundamental building blocks.
Away from the pointwise spatial discontinuities of the flux, these solutions are obtained by means of the standard Riemann solver.
At the discontinuities of $f$, we apply an adapted Riemann solver.
Such Riemann solvers are defined in Section~\ref{s:defRS}.
Unlike the standard Riemann solver, the adapted Riemann solver cannot be defined on the whole physical domain $\calU$. 
This leads to the introduction of a generalised notion of invariant domain and the restriction of $\bU$ to these domains.
We motivate our definition of the adapted Riemann solver by giving a detailed construction in Section~\ref{s:conRS}.
The global existence of the approximate solution $\bU^n$ is ensured by a time-independent bound on the number of wave-fronts.
We emphasise that, in our setting, the number of wave-fronts may increase after an interaction; however, we show that this can occur only at the discontinuities of the flux $f$, and that in such cases the number of waves increases by one.

A further ingredient of the proof is typically the availability of a uniform estimate on the total variation of the approximate solutions $\bU^n(t)$.
However, in our case, $\TV(\bU^n(t))$ may increase after an interaction.
For Temple systems, this issue is usually overcome by considering the total variation of $\bU^n$ expressed in the Riemann invariant coordinates $\bW$.
In contrast, we show that even $\TV(\bW(\bU^n(t)))$ may increase after an interaction, see Remark~\ref{r:Malmsteen}.
This forces us to introduce \lq\lq adapted\rq\rq\ Riemann invariant coordinates $\tbW$, that take into account the discontinuity points of the flux $f$, and that allow us to prove a uniform estimate for $\TV(\tbW(\bU^n(t),c))$.

These results allow us to obtain compactness using Helly's theorem.
We characterise the limit function $\bU$ by showing that it satisfies a Kruzhkov-type entropy condition \cite{Kruzhkov}, that is obtained by applying an argument similar to that in \cite[Section~2.3]{BCR-ARZ-M3AS}.
In particular, a compensating term is designed to account for the additional entropy dissipation at each discontinuity point of the flux.
Uniqueness is then obtained by exploiting the Kruzhkov-type entropy condition and applying \cite[Theorem~1.2]{Gusev}.

Throughout the proof, the choice of a logarithmic \lq\lq pressure\rq\rq\ function $p$ in \eqref{e:flux} plays a key role in our analysis.
This choice prevents the appearance of vacuum.
On the one hand, this may be regarded as a restrictive assumption, for instance, in the modelling of vehicular traffic, see Section~\ref{subsec:applications}.
On the other hand, the logarithmic pressure function has already been justified and adopted in several articles within this framework.
A further discussion of this aspect is provided in Remark~\ref{r:log}.

The paper is organised as follows. 
In the next section, we present the governing equations, list their main properties, and state the principal results in Theorems~\ref{t:uniqueness} and~\ref{t:LornaShore2}. 
The section concludes with a discussion of the application of our result to a generalised version of the ARZ model, presented in Section~\ref{subsec:applications}. 
In Section~\ref{s:RS}, we define the Riemann solver associated with the discontinuity points of the flux, providing a detailed construction and motivation. 
The proofs of Theorems~\ref{t:uniqueness} and~\ref{t:LornaShore2} are deferred to Sections~\ref{s:proofuni} and~\ref{s:proof}, respectively.

\section{The system and main notation}

The aim of this paper is the analysis of the one-dimensional $2\times2$ non-symmetric Keyﬁtz-Kranzer type system of conservation laws with discontinuous flux
\begin{equation}
\left\{\begin{array}{@{}l@{}}
\partial_t \rho + \partial_x \left( c \, V\left(\frac{q}{\rho} - p(\rho)\right) \rho \right)=0
\\[5pt]
\partial_t q + \partial_x \left( c \, V\left(\frac{q}{\rho} - p(\rho)\right) q \right)=0.
\end{array}\right.
\label{e:2x2}
\end{equation}
Here $p \colon (0,\Rref] \to (-\infty,0]$ is defined by
\begin{equation}
\label{e:pc}
p(\rho) \doteq \vref \, \ln(\rho/\Rref),
\end{equation}
where $\vref , \Rref > 0$ are given constants. Observe that for any $\rho\in(0,\Rref]$
\begin{gather}
\label{e:p}
\rho\,p'(\rho) = \vref,
\intertext{and for any $w\leq0$}
\label{e:pm1}
p^{-1}(w) = \Rref \, \exp\left(w/\vref\right) \in (0, \Rref].
\end{gather}
For notation simplicity, we extend $p^{-1}$ to $\R$ by means of the expression in \eqref{e:pm1}.

Given a fixed number of discontinuities $I_c\in\N$, the function $c \colon \R \to (0,+\infty)$ is piecewise constant, right-continuous and takes the form
\begin{equation}\label{e:def_c}
c(x) = \sum_{i = 0}^{I_{c}} c_{i+1/2} \, \mathbbm{1}_{\Xi_{i}}(x),
\end{equation}
for some constants $c_{i+1/2} \in (0,+\infty)$ such that $c_{i-1/2} \neq c_{i+1/2}$, and 
\begin{align*}
&\Xi_{0} \doteq (-\infty,\xi_{1}),&
&\Xi_{i} \doteq [\xi_i,\xi_{i+1}),\ i \in \{1,\ldots,I_{c}-1\},&
&\Xi_{I_{c}} \doteq [\xi_{I_{c}},+\infty),
\end{align*}
with $\xi_i < \xi_{i+1}$.
The \emph{velocity} function $V \colon [0,\ov] \to [0,+\infty)$, with $\ov>0$ be given, is $\C2$ and satisfies the following condition: for any $v \in [0,\ov]$ we have
\begin{align}
\label{e:V}
&V(0)=0,&
&V(v) < \vref \, V'(v),&
V'(v) > 0,&
&V'(v) > \vref \, V''(v).
\end{align}

\begin{example}
The function $V(v) = \sinh(v/\vref)$ satisfies \eqref{e:V} for any $\ov>0$.
On the other hand, the function $V(v) = v$ satisfies \eqref{e:V} only if $\ov \in(0, \vref)$.
\end{example}

For later use, for a fixed $w \le 0$ we introduce the function
\begin{equation}
\label{e:F}
\begin{array}{@{}r@{}r@{\,}c@{\,}l@{}}
F_{w} \colon & \left[ p^{-1}( w - \ov ) , p^{-1}( w ) \right]& \to &\left[ 0 , p^{-1}(w-\ov) \, V(\ov) \right]
\\
&\rho&\mapsto&\rho \, V\left( w - p(\rho) \right).
\end{array}
\end{equation}
Observe that the flux of $\rho$ is $c\,F_{q/\rho}(\rho)$.
We collect the main properties of $F_{w}$ in the following lemma.

\begin{lemma}
\label{l:F}
Fix $w\leq0$.
The function $F_{w}(\rho) \doteq \rho \, V( w - p(\rho))$ given in \eqref{e:F} is well defined, strictly decreasing, strictly concave down, $F_{w}\left(p^{-1}( w )\right) = 0$, and Lipschitz continuous with
\begin{equation}
\label{e:LipF}
\Lip(F_{w}) = \vref \, V'(0).
\end{equation}
Moreover, $F_{w}$ is invertible and its inverse function $F_{w}^{-1}$ is Lipschitz continuous with
\begin{equation}
\label{e:LipFm1}
\Lip(F_{w}^{-1}) = \frac{1}{\vref \, V'(\ov) - V(\ov)}.
\end{equation}
\end{lemma}
\begin{proof}
By \eqref{e:p}, \eqref{e:F} we have
\[F_{w}'(\rho) = V\left( w - p(\rho) \right) - \vref \, V'\left( w - p(\rho) \right).\]
By \eqref{e:V}\textsubscript{$2,4$}, the map $[0,\ov] \ni v \mapsto V(v) - \vref \, V'(v)$ takes values in $(-\infty,0)$ and is strictly increasing.
This readily implies \eqref{e:LipF} and \eqref{e:LipFm1}.
This concludes the proof.
\end{proof}
\noindent
It is worth noting that both $\Lip(F_{w})$ and $\Lip(F_{w}^{-1})$ do not depend on $w$.

\begin{remark}
Both choices $V(v) = v$ and $V(v) = \sinh(v/\vref)$ are well defined on $\R$.
Consequently, in both cases, for $w \le 0$, the function $F_{w}$ defined in \eqref{e:F} admits an extension to $(0,+\infty)$.
Enforcing the condition $F_{w}(\rho) \ge 0$ further confines the domain to $(0, p^{-1}(w)]$.
As observed in Lemma~\ref{l:F}, condition \eqref{e:V}\textsubscript{$2$} implies that $F_{w}'<0$ on the interval $[ p^{-1}( w - \ov ) , p^{-1}( w ) ]$.
It is worth noting that if $V(v) = v$, then $F_{w}'<0$ in $(p^{-1}(w-\vref),p^{-1}(w)]$, which contains $[ p^{-1}( w - \ov ) , p^{-1}( w ) ]$ whenever $\ov \in(0, \vref)$.
In this sense, for a fixed $w \le 0$, the quantity $p^{-1}(w-\vref)$ provides a strictly positive lower bound for the admissible values of $\rho$.
On the other hand, this property no longer holds when $V(v) = \sinh(v/\vref)$.
Indeed, in this case for any $w\leq0$ we have $F_{w}'<0$ in $(0,p^{-1}(w)]$.
For completeness, we point out that $p^{-1}( w - \vref )$ goes to zero as $w$ goes to $-\infty$.
Let us also observe that 
\begin{align*}
V(v) = v &\Longrightarrow \lim_{\rho\to0^+} F_{w}(\rho) = 0,
\\
V(v) = \sinh(v/\vref) &\Longrightarrow \lim_{\rho\to0^+} F_{w}(\rho) = \frac{1}{2} \, p^{-1}(w) > 0.
\end{align*}
\end{remark}

\subsection{Main properties of the system}
\label{s:mp2}

In this section, we collect the main properties of the system of conservation laws \eqref{e:2x2}.
The \lq\lq physical\rq\rq\ domain of the pair of conservative variables $\bU \doteq (\rho,q)$ is the compact set
\begin{equation}
\label{e:U}
\calU \doteq \left\{ (\rho,q) \in (0,\Rref] \times (-\infty,0] : 0 \leq \frac{q}{\rho}-p(\rho) \leq \ov ,\ \uw \leq \frac{q}{\rho} \leq 0 \right\},
\end{equation}
for a constant $\uw <0$.
A possible physical interpretation of the two constants $\ov$ and $\uw$ is the following: the quantity $\|c\|_{\L\infty(\R)} \, V(\ov)>0$ provides a uniform upper bound on the velocity, and 
\[\ur \doteq p^{-1}(\uw-\ov)>0\]
represents a uniform lower bound on the density as we show in the next lemma.
\begin{figure}[!htb]\centering
\begin{tikzpicture}[every node/.style={anchor=south west,inner sep=0pt, black},x=1mm, y=1mm]

\begin{scope}[scale=0.5]
\node at (0,0) {\includegraphics[width=50mm]{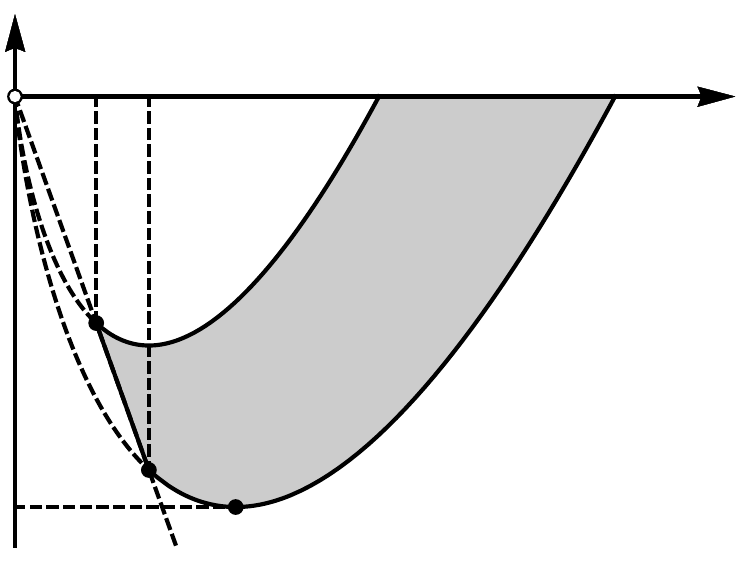}};
\node at (95,53){\strut$\rho$};
\node[left] at (-2,72) {\strut$q$};
\node[left] at (0,7) {\strut$-\Rref \, \vref / e$};
\node at (78,63) {\strut$\Rref$};
\node at (39,63) {\strut$p^{-1}(-\ov)$};
\node at (12,63) {\strut$\ur$};
\node at (23,-3) {\strut$q=\uw\,\rho$};
\node at (18,63) {\strut$p^{-1}(\uw)$};
\end{scope}

\begin{scope}[scale=0.5, shift={(150,0)}]
\node at (0,0) {\includegraphics[width=50mm]{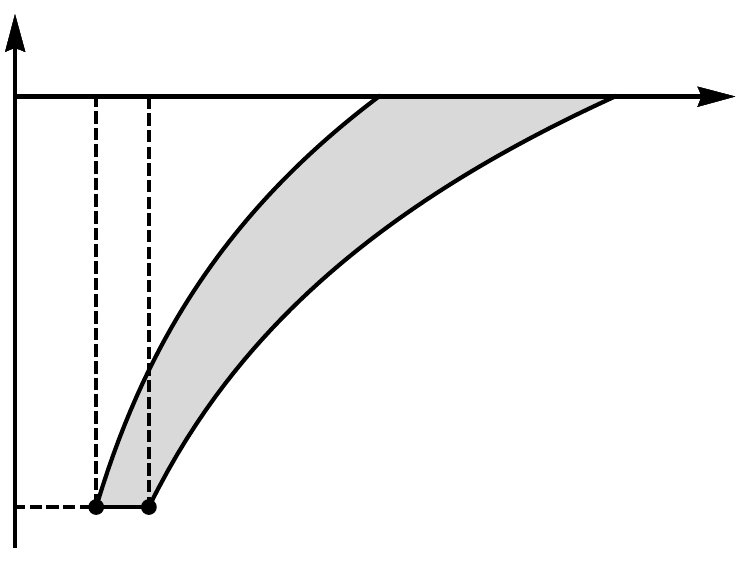}};
\node at (95,53){\strut$\rho$};
\node[left] at (-2,72) {\strut$q/\rho$};
\node[left] at (0,8) {\strut$\uw$};
\node at (78,63) {\strut$\Rref$};
\node at (39,63) {\strut$p^{-1}(-\ov)$};
\node at (12,63) {\strut$\ur$};
\node at (18,63) {\strut$p^{-1}(\uw)$};
\end{scope}
\end{tikzpicture}
\caption{The domain $\calU$ defined in \eqref{e:U} in the $(\rho,q)$ and $(\rho,q/\rho)$ coordinates.}
\label{f:U}
\end{figure}

\begin{lemma}
\label{l:U}
Given $\ov>0$ and $\uw<0$ fixed constants, the set $\calU$ defined in \eqref{e:U} has the qualitative shape illustrated in \figurename~\ref{f:U}, left, and is contained in $[\ur , \Rref]\times[-\Rref \, \vref / e,0]$.
\end{lemma}
\begin{proof}
Observe that if $\rho>0$ then
\[
0 \leq \frac{q}{\rho}-p(\rho) \leq \ov
\qquad\Longleftrightarrow\qquad
\rho\,p(\rho) \leq q \leq \rho\,\left(\ov+p(\rho)\right).
\]
Moreover by \eqref{e:pc} and \eqref{e:p} we have
\begin{align*}
\lim_{\rho\to0^+} \rho\,p(\rho) &= 0,&
\lim_{\rho\to0^+} \rho\,\left(\ov+p(\rho)\right) &= 0,
\\
\frac{\d}{\d\rho}\left(\rho\,p(\rho)\right) &= \vref + p(\rho),&
\frac{\d}{\d\rho}\left(\rho\,\left(\ov+p(\rho)\right)\right) &= \ov + \vref + p(\rho),
\\
\lim_{\rho\to0^+} \frac{\d}{\d\rho}\left(\rho\,p(\rho)\right) &= -\infty,&
\lim_{\rho\to0^+} \frac{\d}{\d\rho}\left(\rho\,\left(\ov+p(\rho)\right)\right) &= -\infty,
\\
\frac{\d^2}{\d\rho^2}\left(\rho\,p(\rho)\right) &= \frac{\vref}{\rho}>0,&
\frac{\d^2}{\d\rho^2}\left(\rho\,\left(\ov+p(\rho)\right)\right) &= \frac{\vref}{\rho}>0.
\end{align*}
In particular, by \eqref{e:pm1} and the above computations, the lowest value of $\rho \mapsto \rho\,p(\rho)$ is
\[q = [\rho\,p(\rho)]_{\rho=p^{-1}(-\vref)} = -p^{-1}(-\vref) \, \vref = -\Rref \, \vref / e,\]
where $e\approx 2.71828$ is the Euler number.
Hence $q\in[-\Rref \, \vref / e,0]$.
At last, we show that $\rho \in [\ur,\Rref]$.
If $\rho>0$ then
\[\uw \leq q/\rho \leq 0 
\qquad\Longleftrightarrow\qquad
\uw \, \rho \leq q \leq 0.\]
Moreover, we have
\[\left\{\begin{array}{@{}l@{}}
q = 0
\\
q = \rho \, p(\rho)
\end{array}\right.
\qquad\Longleftrightarrow\qquad
\left\{\begin{array}{@{}l@{}}
\rho = \Rref
\\
q = 0,
\end{array}\right.\]
and
\[\left\{\begin{array}{@{}l@{}}
q = \uw \, \rho
\\
q = \rho\,\left(\ov+p(\rho)\right)
\end{array}\right.
\qquad\Longleftrightarrow\qquad
\left\{\begin{array}{@{}l@{}}
\rho = p^{-1}(\uw-\ov) = \ur
\\
q = \uw \, p^{-1}(\uw-\ov).
\end{array}\right.\]
This concludes the proof.
\end{proof}

The Jacobian matrix of the flux in equation \eqref{e:2x2} possesses the following eigenvalues and corresponding eigenvectors
\begin{equation}
\label{e:eigen2}
\begin{aligned}
\lambda_{1}(\bU,c) &\doteq c \left(V\left(\frac{q}{\rho} - p(\rho)\right)-\vref \, V'\left(\frac{q}{\rho} - p(\rho)\right)\right),
\qquad
&
R_{1}(\bU) &\doteq \left(
\begin{array}{@{}c@{}}
\rho \\
q
\end{array}
\right),
\\
\lambda_{2}(\bU,c) &\doteq c \, V\left(\frac{q}{\rho} - p(\rho)\right),
&
R_{2}(\bU) &\doteq \left(
\begin{array}{@{}c@{}}
\rho \\
q + \vref \, \rho
\end{array}
\right).
\end{aligned}
\end{equation}

By \eqref{e:V}\textsubscript{$1,2,3$}, system \eqref{e:2x2} is strictly hyperbolic in $\calU$ since $\lambda_{1}(\bU,c) < 0 \leq \lambda_{2}(\bU,c)$ for any $(\bU,c) \in \calU \times (0,+\infty)$.
By \eqref{e:V}\textsubscript{$4$}, the $1$-characteristic field is genuinely non-linear since
\[
D_\bU\lambda_{1}(\bU,c) \cdot R_{1}(\bU) = -c \, \vref \left( V'\left(\frac{q}{\rho} - p(\rho)\right) - \vref \, V''\left(\frac{q}{\rho} - p(\rho)\right) \right) \neq 0
\]
for any $(\bU,c) \in \calU \times (0,+\infty)$.
On the other hand, the $2$-characteristic field is linearly degenerate, namely $D_\bU\lambda_{2}(\bU,c) \cdot R_{2}(\bU) = 0$.
By \cite[Theorem~8.2.1]{Dafermos-book}, any $1$-wave must be either a single centred $1$-rarefaction wave or a single compressive $1$-shock.
By \cite[Theorem~8.2.5]{Dafermos-book}, any $2$-wave is necessarily a $2$-contact discontinuity.
Notice that, in addition to $1$- and $2$-waves, the solution may contain stationary discontinuities located at the points $x=\xi_i$, $i\in\{1,\ldots,I_c\}$, where the function $c$ is discontinuous; we refer to such discontinuities as $c$-waves.

In view of \cite[Theorem~7.3.3]{Dafermos-book}, the functions
\begin{align}
\label{e:KillswitchEngage2}
&v(\bU) \doteq \frac{q}{\rho} - p(\rho),&
&w(\bU) \doteq \frac{q}{\rho},
\end{align}
form a coordinate system of Riemann invariants for \eqref{e:2x2} because $D_{\bU}v(\bU) \cdot R_2(\bU) = 0$, $D_{\bU}w(\bU) \cdot R_1(\bU) = 0$ and
\begin{align*}
&D_{\bU}v(\bU) \cdot R_{1}(\bU) = - \vref,&
&D_{\bU}w(\bU) \cdot R_{2}(\bU) = \vref,
\end{align*}
are different from zero for any $\bU \in \calU$.
Observe that $v$ is a $2$-Riemann invariant and $w$ is a $1$-Riemann invariant, see \cite[Definition~7.3.1]{Dafermos-book}.
Moreover, we highlight that the Riemann invariants and the eigenvectors do not depend on $c$; this simplifies the subsequent analysis.

The conservative variables $\bU$ can be expressed in terms of the coordinate system of Riemann invariants $\bW \doteq (v,w)$ given in \eqref{e:KillswitchEngage2} as follows:
\begin{align}
\label{e:TheMarsVolta2}
&\rho(\bW) \doteq p^{-1}(w-v),&
&q(\bW) \doteq w \, p^{-1}(w-v).
\end{align}
The domain $\calU$ defined in \eqref{e:U} can be reformulated using the coordinate system $\bW$ as
\begin{equation}
\label{e:W}
\mathcal{W} \doteq [0,\ov] \times [\uw,0].
\end{equation}
With a slight abuse of notation we introduce
\begin{align}
\label{e:WU-UW}
\bW(\bU) &\doteq \left( v(\bU) , w(\bU) \right),&
\bU(\bW) &\doteq \left( \rho(\bW) , q(\bW) \right).
\end{align}

By \cite[Theorem~7.6.6]{Dafermos-book}, the $1$-rarefaction curve through $\overline{\bU} = \left(\overline{\rho},\overline{q}\right) \in \calU$ is $\mathcal{R}_{\overline{\bU}}^1(s) \doteq \bigl( \rho^1_{\overline{\bU}}(s), q^1_{\overline{\bU}}(s) \bigr)$, where
\begin{align}
\label{e:C1}
\rho^1_{\overline{\bU}}(s) &\doteq \overline{\rho} \, e^s,&
q^1_{\overline{\bU}}(s) &\doteq \overline{q} \, e^s.
\end{align}
$\mathcal{R}_{\overline{\bU}}^1$ describes a straight line in the state space, hence it coincides with the $1$-shock curve through $\overline{\bU}$, $\mathcal{S}_{\overline{\bU}}^1$, by \cite[Theorem~8.2.7]{Dafermos-book}.
For notation simplicity, we denote them with $\mathcal{C}_{\overline{\bU}}^1$, see \figurename~\ref{f:constru}, right.

By \cite[Theorem~8.2.5]{Dafermos-book}, the $2$-contact discontinuity curve through $\overline{\bU} = \left(\overline{\rho},\overline{q},\right) \in \calU$ is $\mathcal{C}^2_{\overline{\bU}}(s) \doteq \bigl( \rho^2_{\overline{\bU}}(s), q^2_{\overline{\bU}}(s) \bigr)$, where
\begin{align}
\label{e:C2}
\rho^2_{\overline{\bU}}(s) &\doteq \overline{\rho} \, e^s,&
q^2_{\overline{\bU}}(s) &\doteq \left( \overline{q} + \overline{\rho} \, \vref \, s \right) e^s,
\end{align}
see \figurename~\ref{f:constru}, right.
Observe that $\lambda_{2}(\mathcal{C}^2_{\overline{\bU}}(s),c) = \lambda_{2}(\overline{\bU},c)$, see \cite[(8.2.31)]{Dafermos-book}.

System \eqref{e:2x2} falls in the class of Temple systems \cite{Temple01}, characterised by strict hyperbolicity and by coinciding shock and rarefaction curves for the genuinely non-linear characteristic fields.

\subsection{Main result}

In this section, we establish our main result, stated in Theorem~\ref{t:LornaShore2}, which provides an existence and uniqueness result for the Cauchy problem for \eqref{e:2x2} with initial condition
\begin{equation}
\label{e:ini2}
\bU(0,x) = \bU_o(x).
\end{equation}
For notation simplicity, let
\begin{align*}
\D &\doteq [0,+\infty) \times \R,&
\Dr &\doteq (0,+\infty) \times \R,
\end{align*}
be time-space domains.
Introduce $S \colon [0,\ov] \to [0,V(\ov) \, \exp(-\ov/\vref)]$ defined by
\begin{equation}
\label{e:S}
S(v) \doteq V(v) \, \exp(-v/\vref),
\end{equation}
and let
\begin{equation}
\label{e:ts}
\ts(\bU,c) \doteq c \, S\left(v(\bU)\right).
\end{equation}
Notice that $S$ is strictly increasing in $[0,\ov]$ by \eqref{e:V}\textsubscript{$2$}, and is therefore invertible.
Furthermore, we have
\begin{equation}
\label{e:SnarkyPuppy}
\ts(\bU,c) = c \, V\left(v(\bU)\right) \, \exp\left(-\frac{v(\bU)}{\vref}\right) = \frac{c \, V\left(v(\bU)\right) \rho}{p^{-1}\left(w(\bU)\right)} = \frac{c \, F_{w(\bU)}(\rho)}{p^{-1}\left(w(\bU)\right)},
\end{equation}
which follows from \eqref{e:pm1}, \eqref{e:TheMarsVolta2}\textsubscript{$1$} and \eqref{e:F}.

\begin{remark}
At first sight, the expressions for $\ts(\bU,c)$ in \eqref{e:SnarkyPuppy} may appear inconsistent: the first two depend solely on $v(\bU)$, whereas the last one does not explicitly involve it.
The apparent discrepancy arises from the fact that $F_{w(\bU)}(\rho) = \rho \, V(v(\bU))$ and $\rho/p^{-1}(w(\bU)) = \exp(-v(\bU)/\vref)$.
On the one hand, the first two formulations are simpler, as it involves only $v(\bU)$.
On the other hand, when dealing with the Rankine-Hugoniot condition \eqref{e:RH}, the last expression is more convenient.
\end{remark}

For completeness, we recall the definition of weak solution.

\begin{definition}
\label{d:ws2}
Let $\bU_o$ be in $\BV(\R;\calU)$.
A function $\bU = (\rho,q) \in \L\infty((0,+\infty);\BV(\R;\calU)) \cap \C0([0,+\infty);\Lloc1(\R;\calU))$ is a weak solution of the Cauchy problem \eqref{e:2x2}, \eqref{e:ini2} if for any test function $\varphi \in \Cc\infty(\D;\R)$ the following identity holds:
\begin{equation}
\label{e:ws2}
\iint_{\D} \left( \partial_t \varphi + c \, V\left(v(\bU)\right) \partial_x\varphi \right) \bU \,\d x\, \d t 
+ \int_{\R} \bU_{o}(x) \, \varphi(0,x) \, \d x = \mathbf{0}.
\end{equation}
\end{definition}

\begin{remark}
For almost every $t>0$ and any $i\in\{1,\ldots,I_{c}\}$, the traces of a weak solution $\bU(t)$ at  $\xi_{i}$ exist because $\bU$ is in $\L\infty((0,+\infty);\BV(\R;\calU))$.
Moreover, they satisfy
\begin{equation}
\label{e:Nkok}
\ts\left(\bU(t,\xi_{i}^{-}),c(\xi_{i}^{-})\right) = \ts\left(\bU(t,\xi_{i}^{+}),c(\xi_{i}^{+})\right).
\end{equation}
This follows from the Rankine-Hugoniot condition for a \emph{stationary} discontinuity at $x=\xi_{i}$
\[
\left\{\begin{array}{@{}l@{}}
c(\xi_{i}^{-}) \, V(v(\bU(t,\xi_{i}^{-}))) \, \rho(t,\xi_{i}^{-}) = c(\xi_{i}^{+}) \, V(v(\bU(t,\xi_{i}^{+}))) \, \rho(t,\xi_{i}^{+})
\\[5pt]
c(\xi_{i}^{-}) \, V(v(\bU(t,\xi_{i}^{-}))) \, q(t,\xi_{i}^{-}) = c(\xi_{i}^{+}) \, V(v(\bU(t,\xi_{i}^{+}))) \, q(t,\xi_{i}^{+}).
\end{array}\right.
\]
Indeed, this condition implies that either $V(v(\bU(t,\xi_{i}^{-}))) = V(v(\bU(t,\xi_{i}^{+}))) = 0$, or $w(\bU_{-}) = w(\bU_{+})$.
Then \eqref{e:Nkok} follows from \eqref{e:SnarkyPuppy}.
\end{remark}


It is well known that conservation laws may admit multiple weak solutions. 
This motivates the introduction of the entropy condition, an additional admissibility criterion based on well chosen entropy pairs and imposed to select physically relevant weak solutions.
We choose the family of entropy pairs $(\mathfrak{E}_\kappa, \mathfrak{Q}_\kappa)$, $\kappa\in[0,\ov]$, defined by
\begin{align}
\label{e:Deftones2}
\mathfrak{E}_\kappa(v) &\doteq 
\left|e^{-\kappa/\vref} - e^{-v/\vref}\right|,&
\mathfrak{Q}_\kappa(v,c) &\doteq 
-c \left| S(\kappa) - S(v) \right|.
\end{align}
Observe that $(\mathfrak{E}_\kappa, \mathfrak{Q}_\kappa)$ gives indeed an entropy pair for any $\kappa\in[0,\ov]$, as it satisfies condition \cite[(7.4.12)]{Dafermos-book}:
\begin{align*}
&\partial_{v}\mathfrak{Q}_\kappa(v,c) = \lambda_1\left(\bU(\bW),c\right) \partial_{v}\mathfrak{E}_\kappa(v),&
&\partial_{w}\mathfrak{Q}_\kappa(v,c) = \lambda_2\left(\bU(\bW),c\right) \partial_{w}\mathfrak{E}_\kappa(v),
\end{align*}
where $\lambda_{1}(\bU,c)$, $\lambda_{2}(\bU,c)$ are the eigenvalues given in \eqref{e:eigen2}, and $\bU(\bW)$ is the change of variables defined in \eqref{e:WU-UW}\textsubscript{$2$}.
Notice that $\mathfrak{E}_{\kappa}(v) \geq 0 \geq \mathfrak{Q}_{\kappa}(v,c)$ for every $v \in [0,\ov]$.
We give now the definition of entropy solution.
\begin{definition}
\label{d:entro2}
Let $\bU_o$ be in $\BV(\R;\calU)$.
A weak solution $\bU = (\rho,q)$ 
of the Cauchy problem \eqref{e:2x2}, \eqref{e:ini2} in the sense of Definition~\ref{d:ws2} is an entropy solution if for any test function $\varphi \in \Cc\infty(\D;[0,+\infty))$ and any $\kappa \in [0,\ov]$
\begin{equation}
\label{e:entro2}
\begin{aligned}
\iint_{\D} \left( \mathfrak{E}_\kappa\left(v(\bU)\right) \partial_t\varphi + \mathfrak{Q}_\kappa\left(v(\bU),c\right) \partial_x\varphi \right) \d x\, \d t 
+ \int_{\R} \mathfrak{E}_\kappa\left(v(\bU_{o})\right) \, \varphi(0,x) \, \d x
\\
+ S(\kappa) \sum_{i=1}^{I_{c}} \left|c(\xi_{i}^{+})-c(\xi_{i}^{-})\right| \int_{0}^{+\infty} \varphi(t,\xi_{i}) \, \d t &\geq 0,
\end{aligned}
\end{equation}
where $\mathfrak{E}_\kappa$ and $\mathfrak{Q}_\kappa$ are defined in \eqref{e:Deftones2}, and $S$ in \eqref{e:S}.
\end{definition}

The entropy condition \eqref{e:entro2} guarantees the uniqueness of the $v$-component of an entropy solution. 
This is a consequence of the following theorem.

\begin{theorem}
\label{t:uniqueness}
The entropy solutions $\bU_{a}$ and $\bU_{b}$ corresponding to two initial data $\bU_{o}^{a}$ and $\bU_{o}^{b}$ satisfy the following estimate: for any $T>0$ and $R>0$ we have
\begin{equation}
\label{e:Lipesti}
\int_{-R}^{R} \left|v^{b}(T,x) - v^{a}(T,x)\right| \d x \leq e^{\ov/\vref} \int_{-R-T \, \oc \, L}^{R+T \, \oc \, L} \left|v^{b}(0,x) - v^{a}(0,x)\right| \d x,
\end{equation}
where $L \doteq \vref \, V'(0) > 0$.
\end{theorem}
The proof is based on the doubling-of-variables technique and is postponed to Section~\ref{s:proofuni}.
Observe that $L = \max_{v \in [0,\ov]} \left| \vref \, V'(v) - V(v) \right|$ by \eqref{e:V}\textsubscript{$2,4$}.

We introduce the coordinate system
\begin{equation}
\label{e:tbW}
\tbW(\bU,c) \doteq \left(\ts(\bU,c),w(\bU)\right),
\end{equation}
where $\ts$ and $w$ are given in \eqref{e:ts} and \eqref{e:KillswitchEngage2}\textsubscript{$2$}, respectively.
Observe that
\begin{align}
\label{e:TheMarsVolta}
\rho &= 
F_{w(\bU)}^{-1}\left( \frac{\ts(\bU,c)}{c} \, p^{-1}(w(\bU)) \right),&
q &= 
w(\bU) \, F_{w(\bU)}^{-1}\left( \frac{\ts(\bU,c)}{c} \, p^{-1}(w(\bU)) \right).
\end{align}
Hence, for a fixed $c>0$, the domain $\calU$ expressed in the coordinate system $\tbW(\,\cdot\,,c)$ corresponds to
\begin{equation}
\label{e:wW}
\widetilde{\mathcal{W}} \doteq \left[0,c \,S(\ov)\right] \times [\uw,0].
\end{equation}
To prove this, it is sufficient to recall that if $\bU \in \calU$, then $w(\bU) \in [\uw,0]$, $v(\bU) \in [0,\ov]$, and that $S$ defined in \eqref{e:S} is strictly increasing.

In what follows, for any $(a,b) \in \R^{2}$ we let
\[\|(a,b)\| \doteq |a|+|b|.\]
We are now in a position to state the main result of the paper.

\begin{theorem}
\label{t:LornaShore2}
Assume that $V$ satisfies \eqref{e:V} and $c$ is the piecewise constant function in \eqref{e:def_c}. Then, the following hold.
\begin{enumerate}[leftmargin=*]
\item 
The system of conservation laws \eqref{e:2x2} is strictly hyperbolic and falls in the class of Temple systems in the domain $\calU$ given in \eqref{e:U}.
\item 
Fix two constants $\uc, \oV>0$ such that
\begin{equation}
\label{e:cV}
\oV \leq \uc \, S(\ov).
\end{equation}
If $c\geq\uc$, then for every initial datum $\bU_o = (\rho_o,q_o) \in \Lloc1(\R;\calU)$ such that
\begin{equation}
\label{e:inibounds}
\tbW(\bU_o,c) \in \BV(\R;[0,\oV] \times [\uw,0]),
\end{equation}
the Cauchy problem \eqref{e:2x2}, \eqref{e:ini2} admits a unique entropy solution $\bU \in \L\infty(\D;\calU)$ in the sense of Definition~\ref{d:entro2}.
Furthermore, there exists a constant $L_{\bU}>0$ such that for all $0\leq s<t$ and $a<b$ we have
\begin{align}
\label{e:tTVbound}
\TV\left( \tbW\left(\bU(t),c\right) \right) &\leq \TV\left( \tbW(\bU_{o},c)\right),
\\\nonumber
\int_a^b \left\| \bU(t,x)-\bU(s,x) \right\| \, \d x &\leq L_{\bU} \, (t-s),
\\\nonumber
\left\| \tbW\left(\bU(t),c\right) \right\|_{\L\infty(\R)}&\leq \left\| \tbW\left(\bU_o,c\right)\right\|_{\L\infty(\R)}.
\end{align}
\end{enumerate}
\end{theorem}

\begin{proof}
The first statement is already proved in Section~\ref{s:mp2}.
The proof of the second statement is deferred to Section~\ref{s:proof}.
The explicit expression of $L_{\bU}$ is given in the proof of Lemma~\ref{l:LipUt}.
\end{proof}

We conclude this section by observing that, if the assumption \eqref{e:inibounds} in Theorem~\ref{t:LornaShore2} is replaced by
\begin{align*}
&\bU_o \in \BV(\R;\calU),&
&\ts(\bU_{o},c) \in [0,\oV],
\end{align*}
then the theorem remains valid and \eqref{e:tTVbound} can be replaced by
\[
\TV\left( \tbW\left(\bU(t),c\right) \right)
\leq C \left( \TV(\bU_o) + \TV(c)\right),
\]
where the constant $C>0$ is given in the proof of the following Proposition.

\begin{proposition}
If $\bU \in \BV(\R;\calU)$ then $\tbW(\bU,c) \in \BV(\R;\widetilde{\mathcal{W}})$, where $\widetilde{\mathcal{W}}$ is given in \eqref{e:wW}.
Furthermore, there exists a constant $C>0$ such that for any $\bU \in \BV(\R;\calU)$ we have
\[\TV\left( \tbW(\bU,c) \right) \leq C \left( \TV(\bU) + \TV(c) \right).\]
\end{proposition}

\begin{proof}
By \eqref{e:ts} and \eqref{e:KillswitchEngage2}\textsubscript{$1$}, for any $\bU^a,\bU^b \in \calU$ we have
\begin{align*}
\left|\ts(\bU^a,c^a) - \ts(\bU^b,c^b) \right| 
&\leq
\left|\ts(\bU^a,c^a) - \ts(\bU^b,c^a) \right| + 
\left|\ts(\bU^b,c^a) - \ts(\bU^b,c^b) \right| 
\\&\leq 
\oc \, V'(0) \left| v(\bU^a) - v(\bU^b) \right| +
\left|c^a-c^b\right| S(\ov)
\end{align*}
because $\Lip(S) = V'(0)$.
Furthermore, we have
\[
\left| v(\bU^a) - v(\bU^b) \right|
\leq
\frac{\left| q^a-q^b\right|}{\rho^a} - q^b \, \frac{\left|\rho^a-\rho^b\right|}{\rho^a\,\rho^b} + \left|p(\rho^a)-p(\rho^b)\right| 
\leq
\frac{1}{\ur} \left( \left| q^a-q^b\right| + \vref \left(\frac{\Rref}{e\,\ur} + 1 \right) \left|\rho^a-\rho^b\right|\right)
\]
because $\Lip(p) = \vref/\ur$ and by Lemma~\ref{l:U}.
At last, observe that for any $\bU^a,\bU^b \in \calU$ we have
\[
\left|w(\bU^a)-w(\bU^b)\right| \leq 
\frac{1}{\ur} \left( \left| q^a-q^b\right| + \vref \, \frac{\Rref}{e\,\ur} \left|\rho^a-\rho^b\right|\right).
\]
Therefore, it is sufficient to take
\[C \doteq \max\left\{ \frac{\vref}{\ur} \left(
\oc \, V'(0) \left(\frac{\Rref}{e\,\ur} + 1 \right)
+ \frac{\Rref}{e\,\ur}\right)
,
\frac{\oc \, V'(0) + 1}{\ur}
,
S(\ov)
\right\}.\]
This concludes the proof.
\end{proof}

\subsection{Application to the ARZ model with discontinuous flux}
\label{subsec:applications}

The homogeneous version of the second-order Aw-Rascle-Zhang (ARZ) model~\cite{AwRascle2000, zhang2002TRB} for vehicular traffic in the \lq\lq isothermal\rq\rq\ case, which corresponds to choosing $\gamma=0$ in \cite[(2.7)]{Aw2002}, is given by
\begin{equation}
\left\{\begin{array}{@{}l@{}}
\partial_t \rho + \partial_x \left( v \, \rho \right)=0,
\\[5pt]
\partial_t \left(\rho\,w\right) + \partial_x \left( v \, \rho \, w \right)=0.
\end{array}\right.
\label{e:ARZ}
\end{equation}
In this context, $\rho$ and $v=w-p(\rho)$ represent the density and the velocity of the vehicles, $p$ is velocity offset (also referred to as the hesitation function) and $w$ is the Lagrangian marker.

System~\eqref{e:2x2} can be interpreted as an extension of the ARZ model~\eqref{e:ARZ} to the case of a space-dependent speed law $v=c(x)\,V(w-p(\rho))$, which admits the following physical interpretation. 
The classical ARZ model describes traffic flow on a homogeneous road, where the speed law is uniform.
In contrast, the introduction of the coefficient $c(x)$ in \eqref{e:2x2} allows the model to account for road inhomogeneities through abrupt spatial variations of the speed law.
For instance, consider a road with variable size $N(x) \in \R$.
In this case, the continuity equations that replace \eqref{e:ARZ} read as follows:
\begin{equation}
\left\{\begin{array}{@{}l@{}}
N(x) \, \partial_t \rho + \partial_x \left( N(x) \, v \, \rho \right)=0,
\\[5pt]
N(x) \, \partial_t \left(\rho\,w\right) + \partial_x \left( N(x) \, v \, \rho \, w \right)=0.
\end{array}\right.
\label{e:ARZdisco}
\end{equation}
Assume that $N(x)\geq \underline{N}$ with $\underline{N} > 0$.
Introduce a new space variable
\[\tilde{x} = \tilde{x}(x) = \frac{1}{\underline{N}} \int_0^x N(\zeta) \, \d \zeta.
\]
System \eqref{e:ARZdisco} then becomes
\[\left\{\begin{array}{@{}l@{}}
\partial_t \rho + \partial_{\tilde{x}} \left( \frac{N(x)}{\underline{N}} \, v \, \rho \right)=0,
\\[5pt]
\partial_t \left(\rho\,w\right) + \partial_{\tilde{x}} \left( \frac{N(x)}{\underline{N}} \, v \, \rho \, w \right)=0.
\end{array}\right.\]
Since the width of a road typically varies with the number of lanes, it is physically reasonable to assume that $N$ is piecewise constant and takes values in $\N$.
It is then easy to recast the above system into one of the form \eqref{e:2x2}.

\begin{remark}
\label{r:log}
The prototype of the velocity offset originally indicated in \cite[(2.2)]{AwRascle2000} takes the form
\begin{equation}
\label{e:pexpo}
p(\rho) = \rho^\gamma,
\end{equation}
where $\gamma>0$ is a parameter of the model.
This probably explains why such expression has been widely taken into account in the literature, see for instance \cite{BCR-ARZ-M3AS, AndreianovDonadelloRazafisonRollandRosini2016, BenyahiaDonadelloDymskiRosini2018, BrianiCristianiFranzinaIgnoto-2025, DiFraFagioliRosini17, DonadelloPolizziRazafisonRollandRosini-2026, DymskiGoatinRosini2018, GottlichHertyMoutariWeissen2021, HamoriTan-2025, RosiniAnnales, Shen2018, ShengZhang2022, ZhangSheng2021} (list far from being exhaustive!).
On the other hand, Aw and Rascle together with collaborators highlighted in \cite{Aw2002} the importance of the logarithmic velocity offset \eqref{e:pc}.
However, the logarithmic velocity offset has not been given enough attention, see for instance \cite{AndreianovDonadelloRosini2021, BagneriniColomboCorli2006, BenyahiaRosini2016, BenyahiaRosini2017, BenyahiaRosini2020, ChalonsGoatin2007, DalSantoRosiniDymskiBenyahia2017, Goatin2006, Sun2024} for some exceptions.

One of the advantages of working with the power-law velocity offset \eqref{e:pexpo} is that it yields the uniform bound
\[
\|v(t)\|_{\L\infty(\R)} \le \|w_o\|_{\L\infty(\R)}.
\]
However, while the function $p$ is well defined at vacuum, i.e.\ at $\rho=0$, the system~\eqref{e:ARZ} is ill posed in the presence of vacuum, since uniqueness of entropy solutions is not guaranteed, see \cite[Remark~2.1]{BCR-ARZ-M3AS}.

An advantage of the logarithmic velocity offset \eqref{e:pc} is that vacuum does not occur and $\rho > 0$.
On the other hand, there is no \emph{a priori} strictly positive uniform lower bound for the density; moreover, there is no an \emph{a priori} uniform upper bound for the velocity for $(\rho,q) \in (0,+\infty) \times (-\infty,0]$.
This motivates the definition \eqref{e:U} of $\calU$, which guarantees a strictly positive uniform lower bound for the density represented by $\ur>0$, see Lemma~\ref{l:U}, and a finite propagation speed of the waves.
Indeed by \eqref{e:eigen2} and \eqref{e:V}\textsubscript{$1,2,3$} we have
\begin{equation}
\label{e:BlacklitCanopy}
- \|c\|_{\L\infty(\R)} \, \vref \, V'(0) \leq \lambda_{1}(\bU,c) < 0 \leq \lambda_{2}(\bU,c) \leq \|c\|_{\L\infty(\R)} \, V(\ov).
\end{equation}
\end{remark}

We conclude this section by mentioning that system \eqref{e:2x2} also extends the generalised ARZ model \cite[(8)]{Fan2014}, the generic second-order model \cite{lebacque2007generic}, and the model introduced in \cite[(3.16)]{ChiarelloCMS23} to the case of discontinuous fluxes.
We omit the details for the sake of brevity.

\section{Our Riemann Solver}
\label{s:RS}

In this section we restrict our attention to Riemann problems for \eqref{e:2x2}, that is, Cauchy problems for \eqref{e:2x2} with initial condition \eqref{e:ini2} involving piecewise constant initial data $\bU_o$ with at most one discontinuity at $x=0$, namely Cauchy problems for \eqref{e:2x2} with initial condition
\begin{equation}
\label{e:iRie2}
\bU(0,x)=
\left\{\begin{array}{@{}l@{\quad\hbox{if }}l@{}}
\bU_L&x<0,\\
\bU_R&x\geq0,
\end{array}\right.
\end{equation}
where $\bU_L, \bU_R \in \calU$ are given constants.
In addition, we consider the case when $c$ admits at most one discontinuity at $x=0$, that is
\begin{align}
\label{e:cRie2}
c(x) \doteq
\left\{\begin{array}{@{}l@{\quad\hbox{if }}l@{}}
c^-&x<0,\\
c^+&x\geq0,
\end{array}\right.
\end{align}
for some $c^-,c^+>0$ given constants.
These assumptions on $\bU_o$ and $c$ allow to explicitly construct a self-similar weak solution of the corresponding Riemann problems and to characterise it through a Kruzhkov-type entropy condition, see Proposition~\ref{p:InFlames2} below.

A further motivation for this section is to show that the set $\calU$ defined in \eqref{e:U} fails to be an invariant domain, see Example~\ref{ex:Ozzy} below.
This observation forces us to relax the notion of invariant domain, see \eqref{e:Marillion}, and leads us to consider a suitable restriction of the Riemann solver to a domain that does possess the (relaxed) invariance property.

\subsection{Construction of the Riemann solver}
\label{s:conRS}

We now describe in detail the construction of our Riemann solver for \eqref{e:2x2}, in order to motivate our choice.
A key role in this construction is played by the Rankine-Hugoniot condition: any discontinuity $(\bU_-,\bU_+)$ originating from $(t,x) = (0,0)$ and propagating with speed $\sigma \in \R$ satisfies
\begin{equation}
\label{e:RH}
\left\{\begin{array}{@{}l@{}}
c(\sigma^-) \, V\left(v(\bU_-)\right) \rho_- - c(\sigma^+) \, V\left(v(\bU_+)\right) \rho_+ = \sigma \, (\rho_- - \rho_+)
\\[5pt]
c(\sigma^-) \, V\left(v(\bU_-)\right) q_- - c(\sigma^+) \, V\left(v(\bU_+)\right) q_+ = \sigma \, (q_- - q_+).
\end{array}\right.
\end{equation}

As a first step in the construction, we look for the admissible $c$-waves $(\hat{\bU},\check{\bU})$ at the interface $x=0$.
Since $\lambda_{1}(\bU,c) < 0 \leq \lambda_{2}(\bU,c)$ for any $(\bU,c) \in \calU \times (0,+\infty)$, every $1$-wave has strictly negative propagation speed, whereas every $2$-contact discontinuity has non-negative propagation speed.
This forces $\hat{\bU}$ to lie on the $1$-curve $\mathcal{C}_{\bU_L}^1$ given in \eqref{e:C1}, $\check{\bU}$ to lie on the $2$-curve $\mathcal{C}_{\bU_R}^2$ expressed in \eqref{e:C2}, and $(\hat{\bU}, \check{\bU})$ to satisfy \eqref{e:RH} with $\sigma=0$. 
This leads to seek for $\hat{s},\check{s}\in\R$ such that
\begin{align}
\label{e:condUhUc}
&\hat{\bU} = \mathcal{C}_{\bU_L}^1(\hat{s}),&
&c^- \, V\bigl(v(\hat{\bU})\bigr) \, \hat{\bU} = c^+ \, V\bigl(v(\check{\bU})\bigr) \, \check{\bU},&
&\check{\bU} = \mathcal{C}_{\bU_R}^2(\check{s}).
\end{align}
By solving the above system we get
\begin{align}
\label{e:UhUc}
\hat{\bU} &= F_{w(\bU_L)}^{-1}\left( \frac{c^+}{c^-} \, F_{w(\bU_L)}\left(\rho_R \, \frac{p^{-1}(w(\bU_L))}{p^{-1}(w(\bU_R))}\right) \right) \left(1,w(\bU_L)\right),&
\check{\bU} &= \rho_R \, \frac{p^{-1}(w(\bU_L))}{p^{-1}(w(\bU_R))} \left(1,w(\bU_L)\right),
\end{align}
where $F_{w}$ is defined in \eqref{e:F}.
In particular, we emphasise that the discontinuous coefficient $c$ affects the determination of $\hat{\bU}$ only, and does not influence $\check{\bU}$.

For completeness, we express $\hat{\bU}$ and $\check{\bU}$ in the Riemann invariant coordinates \eqref{e:KillswitchEngage2}.
First, observe that
\begin{align}
\label{e:Frost0}
&v(\check{\bU}) = v(\bU_R),&
&w(\hat{\bU}) = w(\check{\bU}) = w(\bU_L),
\end{align}
see again \cite[Theorems~7.3.3 and~8.2.5]{Dafermos-book}.
Furthermore, we have
\begin{align*}
v(\hat{\bU}) 
&= w(\bU_L) - p\left( F_{w(\bU_L)}^{-1}\left( \frac{c^+}{c^-} \, F_{w(\bU_L)} \left( p^{-1}\left(w(\bU_R)-v(\bU_R)\right) \frac{p^{-1}(w(\bU_L))}{p^{-1}(w(\bU_R))} \right) \right) \right)
\\
&= w(\bU_L) - p\left( F_{w(\bU_L)}^{-1}\left( \frac{c^+}{c^-} \, F_{w(\bU_L)} \left( p^{-1}\left(w(\bU_L)-v(\bU_R)\right) \right) \right) \right).
\end{align*}

\begin{figure}[!htb]\centering
\maxsizebox{.8\linewidth}{!}{
\begin{tikzpicture}[every node/.style={anchor=south west,inner sep=0pt, black},x=1mm, y=1mm]

\begin{scope}[scale=.8]
\node at (0,0) {\includegraphics[width=64mm]{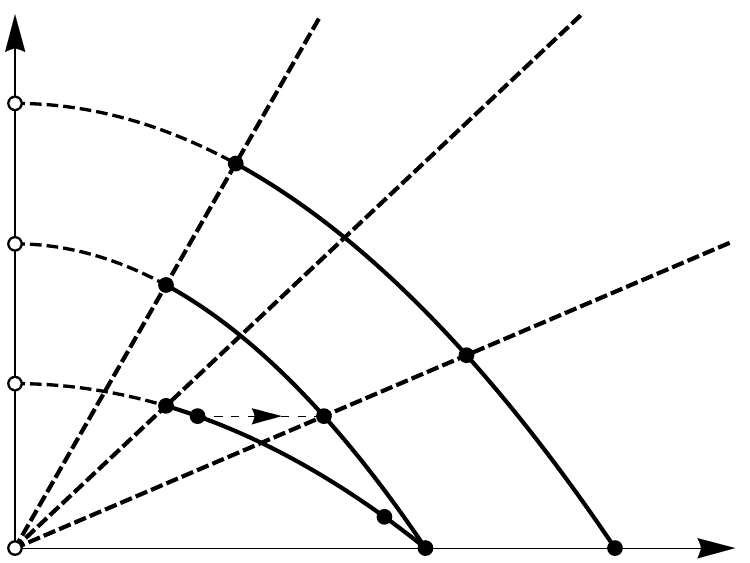}};
\node[below] at (78.4,0){\strut$\rho$};
\node[left] at (0,58) {\strut$c\,F_{w}$};
\node[left] at (0,20) {\strut$c^{-}F_{w^L}$};
\node[left] at (0,35) {\strut$c^{+}F_{w^L}$};
\node[left] at (0,50) {\strut$c^{+}F_{w^R}$};
\node[below] at (65,0) {\strut$p^{-1}(w^R)$};
\node[below] at (45,0) {\strut$p^{-1}(w^L)$};
\node[right, inner sep=2pt, rotate=60] at (30,43) {\strut$c^{+} \, \rho \, V(\ov)$};
\node[right, inner sep=2pt, rotate=43] at (48,41) {\strut$c^{-} \, \rho \, V(\ov)$};
\node[right, inner sep=2pt, rotate=23] at (56,28) {\strut$c^{+} \, \rho \, V(v(\bU_{R}))$};
\node at (18,8.5) {\strut$\hat{\bU}$};
\node at (37.5,10.5) {\strut$\check{\bU}$};
\node at (52,17.5) {\strut$\bU_R$};
\node at (33,1) {\strut$\bU_L$};
\end{scope}

\begin{scope}[shift={(80,0)}]

\begin{scope}[scale=.8]
\node at (0,0) {\includegraphics[width=64mm]{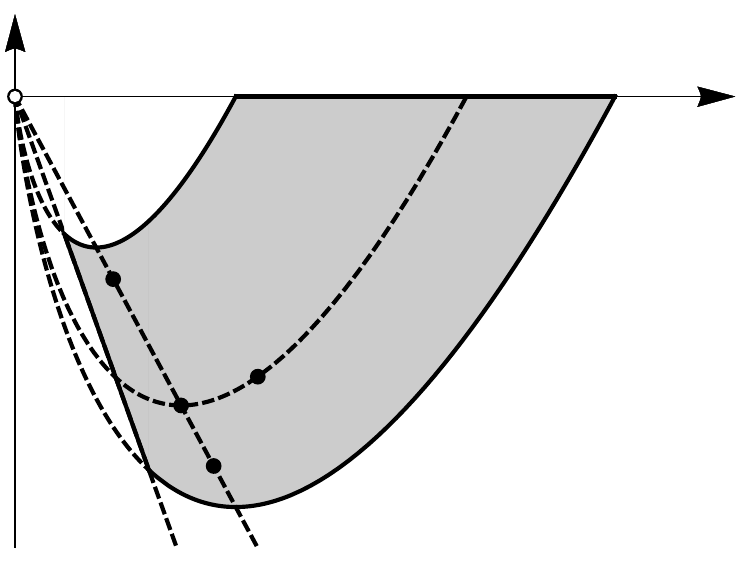}};
\node[below] at (78.4,48){\strut$\rho$};
\node[left] at (0,58) {\strut$q$};
\node[below] at (30,56) {\strut$p^{-1}(-\ov)$};
\node[below] at (65,56) {\strut$\Rref$};
\node[right, inner sep=2pt] at (47,41) {\strut$\mathcal{C}^2_{\bU_{R}}$};
\node[right, inner sep=2pt] at (28,0) {\strut$\mathcal{C}^1_{\bU_{L}}$};
\node at (19,17) {\strut$\check{\bU}$};
\node at (14,28) {\strut$\hat{\bU}$};
\node at (29,15) {\strut$\bU_R$};
\node at (24.5,8) {\strut$\bU_L$};
\end{scope}

\end{scope}

\end{tikzpicture}}
\caption{Left: Graphical illustration in the $(\rho,c\,F_{w})$-plane of the construction of the intermediate states $\hat{\bU}$ and $\check{\bU}$, as described in Remark~\ref{r:constru}.
Right: Representation of the states $\hat{\bU}$ and $\check{\bU}$ in the $(\rho,q)$-plane.
Here $V(v) = \sinh(v/\vref)$, $c^+ > c^-$, $w^L = w(\bU_L) < w^R = w(\bU_R)$.
}
\label{f:constru}
\end{figure}

\begin{remark}
\label{r:constru}
The intermediate states $\hat{\bU}$ and $\check{\bU}$ can be obtained as follows, see \figurename~\ref{f:constru}, left.
By \eqref{e:Frost0} we have that $\bW(\check{\bU}) = (v(\bU_R),w(\bU_L))$.
Once $\check{\bU}$ is determined, $\hat{\bU}$ is obtained as follows.
By \eqref{e:Frost0}\textsubscript{$2$} we have $w(\hat{\bU}) = w(\check{\bU})$.
Denote this common value by $w$; then by \eqref{e:UhUc} we have
\begin{equation}
\label{e:Testament}
c^- F_{w}(\hat{\rho}) = c^+ F_{w}(\check{\rho}),
\end{equation}
and therefore
\begin{equation}
\label{e:TheCult}
\hat{\rho} = F_{w}^{-1}\left( \frac{c^+}{c^-} \, F_{w}(\check{\rho}) \right).
\end{equation}
\end{remark}

From equation \eqref{e:TheCult} readily follows that:
\begin{align}
&\hat{\rho} < \check{\rho} \Longleftrightarrow (c^--c^+) \, v(\check{\bU}) < 0,&
&\hat{\rho} = \check{\rho} \Longleftrightarrow (c^--c^+) \, v(\check{\bU}) = 0.
\label{e:UI=UJ}
\end{align}

In addition to verifying system \eqref{e:condUhUc}, we must also ensure that $\hat{\bU},\check{\bU} \in \calU$.
From \eqref{e:Frost0} and \eqref{e:W}, we deduce that $\check{\bU} \in \calU$, since $v(\check{\bU}) \in [0,\ov]$ and $w(\check{\bU}) \in [\uw,0]$, as $\bU_L,\bU_R \in \calU$.
The case $\hat{\bU} = \check{\bU}$ is trivial. Hence, it is sufficient to consider the case $\hat{\bU} \neq \check{\bU}$.
By \eqref{e:UI=UJ} we have that in this case $(c^--c^+) \, v(\check{\bU}) \neq 0$.
Consequently, $c^-\neq c^+$ and \eqref{e:Frost0} implies that any $c$-wave $(\hat{\bU},\check{\bU})$ satisfies
\begin{align}
\label{e:Frost}
&v(\check{\bU}) = v(\bU_R) > 0,&
&w(\hat{\bU}) = w(\check{\bU}) = w(\bU_L).
\end{align}
The following example shows that the condition $\hat{\bU} \in \calU$ is not automatically satisfied by solutions of \eqref{e:condUhUc}.

\begin{example}
\label{ex:Ozzy}
Assume that $c^- < c^+$.
Take $\bU_o \in \calU$ so that $v(\bU_o) \equiv \ov$ and consider the Riemann problem \eqref{e:2x2}, \eqref{e:iRie2}, \eqref{e:cRie2} with $\bU_L = \bU_R = \bU_o$.
In this case, either $\hat{\bU}$ given by \eqref{e:UhUc} is not well defined, or $v(\hat{\bU}) > \ov$ and therefore $\hat{\bU}$ doesn't belong to $\calU$.
\end{example}

On the one hand, $\calU$ is the natural physical domain for the conserved variables $\bU$. 
On the other hand, it turns out that, for data in $\calU$, the above construction may fail, as illustrated by Example~\ref{ex:Ozzy}.
A possible explanation of this drawback is that $\calU$ does not take into account the discontinuous coefficient $c$.
This motivates restricting the model to a suitable invariant domain contained in $\calU$.

In the framework of Temple systems, it is often convenient to construct invariant domains in the coordinates defined by the Riemann invariants \eqref{e:KillswitchEngage2}.
In the present case, however, it turns out to be more appropriate to consider a coordinate system obtained by modifying the first Riemann invariant $v$ so as to account for the discontinuous coefficient~$c$.
More precisely, we consider the coordinate system $\tbW(\bU,c)$ given in \eqref{e:tbW}.
Notice that by \eqref{e:Frost}\textsubscript{$2$}, \eqref{e:Testament}, \eqref{e:SnarkyPuppy} we have
\begin{align}
\label{e:Leprous}
&\ts(\hat{\bU},c^-) = \ts(\check{\bU},c^+) = \ts(\bU_{R},c^+),&
&w(\hat{\bU}) = w(\check{\bU}) = w(\bU_L).
\end{align}

\begin{figure}[!htb]\centering

\begin{tikzpicture}[every node/.style={anchor=south west,inner sep=0pt, black},x=1mm, y=1mm]
\node at (0,-5) {\includegraphics[width=58mm]{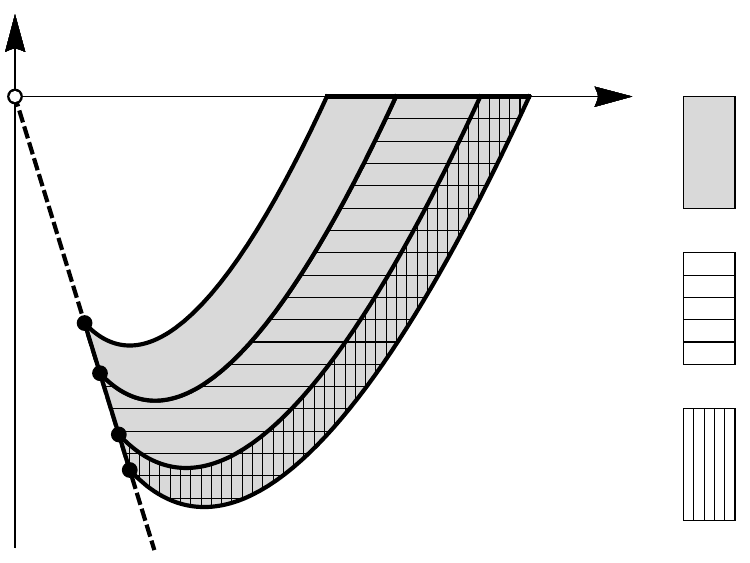}};
\node at (49,26){\strut$\rho$};
\node at (-3,35) {\strut$q$};
\node at (38,32) {\strut$\Rref$};
\node at (17,32) {\strut$p^{-1}(-\ov)$};
\node at (59,25) {\strut$\calU$};
\node at (59,13) {\strut$\calU_{c^{-}}$};
\node at (59,0) {\strut$\calU_{c^{+}}$};
\end{tikzpicture}
\caption{The domains $\calU$, $\calU_{c^{-}}$ and $\calU_{c^{+}}$ defined in \eqref{e:U} and \eqref{e:dominva}. Here $\uc \leq c^{-} < c^{+}$ and therefore $\calU_{c^{+}} \subset \calU_{c^{-}} \subset \calU$.}
\label{f:Owane}
\end{figure}

In the following lemma we show that if
\begin{align}
\label{e:StevenWilson}
&\uc, \oV > 0,&
&\oV \leq \uc \, S(\ov),&
&c^-,c^+ \geq \uc,
\end{align}
then we can consider
\begin{equation}
\label{e:dominva}
\calU_{c^{\pm}} \doteq 
\left\{ \bU \in \calU : \ts(\bU,c^\pm) \in [0 , \oV] \right\},
\end{equation}
see \figurename~\ref{f:Owane}, and get the following result.

\begin{lemma}
\label{l:traces}
Assume \eqref{e:StevenWilson}.
If $(\bU_L,\bU_R) \in \calU_{c^-} \times \calU_{c^+}$, with $\calU_{c^{\pm}}$ given in \eqref{e:dominva}, then the pair $(\hat{\bU},\check{\bU})$ determined by \eqref{e:UhUc} lies in $\calU_{c^-} \times \calU_{c^+}$.
\end{lemma}
\begin{proof}
It readily follows from the construction that $\check{\bU} \in \calU$, see \figurename~\ref{f:constru}, right.
By \eqref{e:Leprous}\textsubscript{$1$} and \eqref{e:dominva} we have $\ts(\hat{\bU},c^-) = \ts(\check{\bU},c^+) = \ts(\bU_{R},c^+) \in [0,\oV]$.
It remains to prove that $\hat{\bU} \in \calU$ by showing that $v(\hat{\bU}) \in [0,\ov]$.
By \eqref{e:Frost0}\textsubscript{$2$} we have $w(\hat{\bU}) = w(\check{\bU}) = w(\bU_L) \in [\uw,0]$; denote by $w$ this common value.
Recall that $\ts(\check{\bU},c^+) \leq \oV$, hence by \eqref{e:SnarkyPuppy}, \eqref{e:cV}, \eqref{e:S} and \eqref{e:pm1}
\begin{align*}
&0 \leq \ts(\check{\bU},c^+) = c^+ \, \frac{F_{w}(\check{\rho})}{p^{-1}(w)} \leq \oV \leq \uc \, S(\ov) \leq c^- \, S(\ov)
\\ \Longrightarrow{}&
0 \leq \frac{c^+}{c^-} \, F_{w}(\check{\rho}) \leq 
p^{-1}(w) \, V(\ov) \, \frac{\Rref}{p^{-1}(\ov)} = 
p^{-1}(w-\ov) \, V(\ov)
\\ \Longrightarrow{}&
\frac{c^+}{c^-} \, F_{w}(\check{\rho}) \in \left[0,p^{-1}(w-\ov) \, V(\ov)\right].
\end{align*}
By \eqref{e:TheCult} and Lemma~\ref{l:F} we have then that
\[\hat{\rho} \in \left[ p^{-1}( w - \ov ) , p^{-1}( w ) \right] \subset (0,\Rref].\]
The above condition also implies that $0\leq v(\hat{\bU}) = w - p(\hat{\rho}) \leq \ov$.
At last, by \eqref{e:Frost}\textsubscript{$2$} we also have $\hat{q} = w \, \hat{\rho} \leq 0$.
This concludes the proof.
\end{proof}

Since $S$ given in \eqref{e:S} is strictly increasing, we have that $\ts(\bU,c^{\pm}) \in [0,\oV]$ if and only if $v(\bU) \in [0,S^{-1}(\oV/c^{\pm})]$.
Observe that \eqref{e:cV} ensures that $S^{-1}(\oV/c) \in (0,\ov]$.
Furthermore, if $c^{+} > c^{-}$ then $S^{-1}(\oV/c^{+}) < S^{-1}(\oV/c^{-})$ and therefore $\calU_{c^{+}} \subset \calU_{c^{-}}$, see \figurename~\ref{f:Owane}.

We can now construct our Riemann solver
\[\RS_{c^-,c^+} \colon \calU_{c^{-}}\times \calU_{c^{+}} \to \BV(\R;\calU),\] 
associated with given constants $c^-,c^+ \geq \uc$, as follows.
In a rough formulation, for any $(\bU_L,\bU_R) \in \calU_{c^-} \times \calU_{c^+}$, \emph{typically} we define $\RS_{c^-,c^+}[\bU_L,\bU_R]$ as the juxtaposition of the $1$-wave $(\bU_L,\hat{\bU})$, the $c$-wave $(\hat{\bU},\check{\bU})$, and the $2$-contact discontinuity $(\check{\bU},\bU_R)$, where $\hat{\bU}$ and $\check{\bU}$ are defined in~\eqref{e:UhUc}. 
Moreover, if one of the following equivalent conditions is satisfied
\begin{align*}
&\rho_L<\hat{\rho},
&v(\bU_L) > v(\hat{\bU}),&
&\ts(\bU_L,c^-) > \ts(\hat{\bU},c^-),
\end{align*}
then the $1$-wave $(\bU_L,\hat{\bU})$ is a $1$-shock, otherwise it is a $1$-rarefaction.

By Lemma~\ref{l:traces} we have that $(\hat{\bU},\check{\bU}) \in \calU_{c^-} \times \calU_{c^+}$.
Furthermore, since $v$ is strictly increasing along any $1$-rarefaction, it follows that also $\ts$ satisfies this property, as the map $S$ given in \eqref{e:S} is strictly increasing.
It is then clear that $\calU_{c^-} \times \calU_{c^+}$ is an invariant domain for $\RS_{c^-,c^+}$, in the sense that
\begin{equation}
\label{e:Marillion}
(\bU_L,\bU_R) \in \calU_{c^-} \times \calU_{c^+}
\Longrightarrow
\left\{\begin{array}{@{}l@{}}
\nu < 0 \Longrightarrow \RS_{c^-,c^+}[\bU_L,\bU_R](\nu) \in \calU_{c^{-}}
\\
\nu \geq 0 \Longrightarrow \RS_{c^-,c^+}[\bU_L,\bU_R](\nu) \in \calU_{c^{+}}.
\end{array}\right.
\end{equation}

\begin{remark}
For completeness, we remark that for any $\overline{w} \in [\uw,0]$ the set
\[\left\{ \bU \in \calU_{c^-} : w(\bU) \in [\uw , \overline{w}] \right\} \times \left\{ \bU \in \calU_{c^+} : w(\bU) \in [\uw , \overline{w}] \right\}\]
is an invariant domain for $\RS_{c^-,c^+}$.
Moreover, if $c^+ \leq c^-$ then $\calU\times\calU$ is an invariant domain for $\RS_{c^-,c^+}$.
\end{remark}

We now determine the propagation speed of the waves.
\begin{itemize}
\item
If $(\bU_L,\hat{\bU})$ is a $1$-rarefaction, then its cone of propagation is
\[\lambda_{1}(\bU_L,c^-) \leq x/t \leq \lambda_{1}(\hat{\bU},c^-).\] 
Let us underline that by \eqref{e:Frost0}\textsubscript{$2$}
\begin{align*}
&\lambda_{1}(\bU_L,c^-) < \lambda_{1}(\hat{\bU},c^-)
\Longleftrightarrow\\&
V\left( w(\bU_L) - p(\rho_L) \right) - \vref \, V'\left( w(\bU_L) - p(\rho_L) \right) < 
V\left( w(\bU_L) - p(\hat{\rho}) \right) - \vref \, V'\left( w(\bU_L) - p(\hat{\rho}) \right),
\end{align*}
that holds true for $\rho_L > \hat{\rho}$ because the function $\rho \mapsto V\left( w - p(\rho) \right) - \vref \, V'\left( w - p(\rho) \right)$ is strictly decreasing. 
This follows from the fact that $v \mapsto V(v)-\vref \,V'(v)$ is increasing by \eqref{e:V}\textsubscript{$4$} and that $p$ is increasing by \eqref{e:pc}.
\item
If $(\bU_L,\hat{\bU})$ is a $1$-shock, then its propagation speed $\sigma = \sigma(\bU_L,\hat{\bU},c^{-})$ satisfies the Rankine-Hugoniot condition \eqref{e:RH}, hence by \eqref{e:KillswitchEngage2}\textsubscript{$2$} and \eqref{e:Frost}\textsubscript{$2$} we have
\begin{align}\nonumber
&\left\{\begin{array}{@{}l@{}}
c^- \, V\bigl(v(\bU_L)\bigr) \, \rho_L - c^- \, V\bigl(v(\hat{\bU})\bigr) \, \hat{\rho} = \sigma \, ( \rho_L-\hat{\rho})
\\[5pt]\nonumber
c^- \, V\bigl(v(\bU_L)\bigr) \, q_L - c^- \, V\bigl(v(\hat{\bU})\bigr) \, \hat{q} = \sigma \, (q_L-\hat{q})
\end{array}\right.
\\\nonumber
\Longleftrightarrow\ &\left\{\begin{array}{@{}l@{}}
c^- \left( F_{w(\bU_L)}(\rho_L) - F_{w(\bU_L)}(\hat{\rho}) \right) = \sigma \, ( \rho_L-\hat{\rho})
\\[5pt]
c^- \left( F_{w(\bU_L)}(\rho_L) - F_{w(\bU_L)}(\hat{\rho}) \right) w(\bU_L) = \sigma \, (q_L-\hat{q})
\end{array}\right.
\\ \Longleftrightarrow\ &
\sigma(\bU_L,\hat{\bU},c^{-}) = c^- \, \frac{F_{w(\bU_L)}(\rho_L) - F_{w(\bU_L)}(\hat{\rho})}{\rho_L-\hat{\rho}},
\label{e:BringMeTheHorizon}
\end{align}
where $F_{w}$ is defined in \eqref{e:F}.
\item
The $2$-contact discontinuity $(\check{\bU},\bU_R)$ has propagation speed
\[\lambda_{2}\bigl(\check{\bU},c^{+}\bigr) = \lambda_{2}\bigl(\bU_R,c^{+}\bigr) = c^+ \, V\bigl(v(\bU_R)\bigr).\]
\end{itemize}

Observe that $\rho_L,\hat{\rho} \in \left[ p^{-1}( w(\bU_L) - \ov ) , p^{-1}( w(\bU_L) ) \right]$ because $(\rho_L,q_L),(\hat{\rho},\hat{q}) \in \calU$ and $w(\hat{\bU}) = w(\bU_L)$.
This implies that $\sigma(\bU_L,\hat{\bU},c^{-}) < 0$ by Lemma~\ref{l:F}.
Recall also that $\lambda_{1}(\bU,c^{-}) < 0$ for any $\bU \in \calU$.
Therefore, any $1$-wave has negative propagation speed.
On the contrary, as $\lambda_{2}(\bU,c^{+}) \geq 0$ for any $\bU \in \calU$, any $2$-contact discontinuity propagates with non-negative speed.

For completeness, observe that the $2$-contact discontinuity is stationary if and only if $v(\bU_R) = v(\check{\bU}) = 0$ and $w(\bU_R) \neq w(\check{\bU})$, but in this case $\hat{\bU}=\check{\bU}$ by \eqref{e:UI=UJ}\textsubscript{$2$} and no $c$-wave is involved.

\subsection{Definition and main properties of the Riemann solver}
\label{s:defRS}

The arguments developed in the previous section lead to the following definition.

\begin{definition}
\label{d:RS}
Assume \eqref{e:StevenWilson}.
For any $(\bU_L,\bU_R) \in \calU_{c^{-}}\times\calU_{c^{+}}$, introduce $(\hat{\bU},\check{\bU}) \in \calU_{c^{-}}\times\calU_{c^{+}}$ implicitly defined by the following conditions:
\begin{align}
\label{e:Avatar}
&w(\hat{\bU}) = w(\check{\bU}) = w(\bU_L),&
&v(\check{\bU}) = v(\bU_R),&
&c^- \, F_{w(\bU_L)}(\hat{\rho}) = c^+ \, F_{w(\bU_L)}(\check{\rho}).
\end{align}
The Riemann solver $\RS_{c^-,c^+}$ computed at $(\bU_L,\bU_R)$ is defined by
\[\RS_{c^-,c^+}[\bU_L,\bU_R](\nu) \doteq
\left\{\begin{array}{@{}l@{\quad\hbox{if }}l@{}}
\RS_{c^-}[\bU_L,\hat{\bU}] & \nu < 0,
\\
\RS_{c^+}[\check{\bU},\bU_R] & \nu \geq 0,
\end{array}\right.\]
where 
\[\RS_{c^+}[\check{\bU},\bU_R](\nu) \doteq
\left\{\begin{array}{@{}l@{\quad\hbox{if }}l@{}}
\check{\bU} & \nu < \lambda_{2}(\bU_R,c^+),
\\
\bU_R & \nu \geq \lambda_{2}(\bU_R,c^+),
\end{array}\right.\]
and $\RS_{c^-} \colon \{(\bU_{L},\hat{\bU}) \in \calU_{c^{-}}\times\calU_{c^{-}} : w(\bU_{L}) = w(\hat{\bU})\} \to \BV(\R;\calU)$ is defined according to the following cases:
\begin{enumerate}
\item 
If $\rho_L < \hat{\rho}$, then we let
\[\RS_{c^-}[\bU_L,\hat{\bU}](\nu) \doteq
\left\{\begin{array}{@{}l@{\quad\hbox{if }}l@{}}
\bU_L & \nu < c^- \, \frac{F_{w(\bU_L)}(\rho_L) - F_{w(\bU_L)}(\hat{\rho})}{\rho_L-\hat{\rho}},
\\[7pt]
\hat{\bU} & \nu \geq c^- \, \frac{F_{w(\bU_L)}(\rho_L) - F_{w(\bU_L)}(\hat{\rho})}{\rho_L-\hat{\rho}}.
\end{array}\right.\]

\item 
If $\rho_L \geq \hat{\rho}$, then we let
\[\RS_{c^-}[\bU_L,\hat{\bU}](\nu) \doteq
\left\{\begin{array}{@{}l@{\quad\hbox{if }}l@{}}
\bU_L & \nu \leq \lambda_{1}(\bU_L,c^-),
\\
\bU & \lambda_{1}(\bU_L,c^-) \leq \nu \leq \lambda_{1}(\hat{\bU},c^-),\ \nu = \lambda_{1}(\bU,c^-),
\\
\hat{\bU} & \nu \geq \lambda_{1}(\hat{\bU},c^-).
\end{array}\right.\]
\end{enumerate}
\end{definition}

\begin{remark}
\label{r:traces}
Some comments on the above definition are in order.
For notation simplicity, let
\[\bU_\pm = \RS_{c^-,c^+}[\bU_L,\bU_R](0^\pm).\]
\begin{enumerate}[label={\bf{(\Alph*)}}]
\item
Definition~\ref{d:RS} simplifies when $c^{-} = c^{+}$, since in this case we have $\hat{\bU} = \check{\bU}$, see \eqref{e:UI=UJ}\textsubscript{$2$} and \eqref{e:Frost}\textsubscript{$2$}.
\item\label{i:traces}
For any $(\bU_L,\bU_R) \in \calU_{c^{-}}\times\calU_{c^{+}}$ we have
\begin{equation}
\label{e:Slipknot}
\bU_- = \hat{\bU},
\end{equation}
also in the case $\bU_{L} = \hat{\bU}$.
On the other hand, we have 
\[\bU_+ \in \{\check{\bU} , \bU_{R} \}.\]
Moreover, $\bU_+ = \bU_{R} \neq \check{\bU}$ if and only if $v(\bU_{R}) = 0$ and $w(\bU_{R}) \neq w(\bU_{L})$. 
Furthermore, in this case also $v(\check{\bU}) = 0$ by \eqref{e:Avatar}\textsubscript{$2$}, hence by \eqref{e:Avatar}\textsubscript{$3$}
\[
c^- \, F_{w(\bU_L)}(\hat{\rho}) = c^+ \, F_{w(\bU_L)}(\check{\rho}) = 0,
\]
and therefore also $v(\hat{\bU}) = 0$.
In particular, in this case we have
\[
F_{w(\bU_-)}(\rho_-) = 0 = F_{w(\bU_+)}(\rho_+).
\]
\item
\label{i:Rammstein}
For any $(\bU_L,\bU_R) \in \calU_{c^{-}}\times\calU_{c^{+}}$ we have
\begin{equation}
\label{e:pezza}
\left( w(\hat{\bU}) - w(\bU_+) \right) v(\hat{\bU}) = 0.
\end{equation}
Indeed, if $\bU_+ = \check{\bU}$, then $w(\hat{\bU}) = w(\bU_+)$ by \eqref{e:Avatar}\textsubscript{$1$}.
On the other hand, if $\bU_+ = \bU_{R} \neq \check{\bU}$, then as already observed in \ref{i:traces}, we have $v(\bU_{R}) = 0$ which implies that $v(\check{\bU}) = 0$ by \eqref{e:Avatar}\textsubscript{$2$}, and therefore $v(\hat{\bU}) = 0$ by \eqref{e:Avatar}\textsubscript{$3$}.
\item 
For any $(\bU_L,\bU_R) \in \calU_{c^{-}}\times\calU_{c^{+}}$ we have
\begin{equation}
\label{e:Katatonia}
c^- \, F_{w(\bU_-)}(\rho_-) = c^+ \, F_{w(\bU_+)}(\rho_+).
\end{equation}
Indeed, if $\bU_+ = \check{\bU}$ then it readily follows from \eqref{e:Avatar}\textsubscript{$1,3$} and \eqref{e:Slipknot}.
On the other hand, if $\bU_+ = \bU_{R} \neq \check{\bU}$, then as already observed in \ref{i:Rammstein}, we have $v(\bU_{R}) = 0 = v(\hat{\bU})$ and therefore $F_{w(\bU_-)}(\rho_-) = 0 = F_{w(\bU_+)}(\rho_+)$ by \eqref{e:Slipknot}.

\begin{figure}[!htb]\centering
\resizebox{.9\linewidth}{!}{\begin{tikzpicture}[every node/.style={anchor=south west,inner sep=3pt, black},x=6mm, y=6mm, black]
\draw[-{Latex[scale=1.6]}] (0,0) -- (12.8,0) node[below] {$\ts$};
\draw[-{Latex[scale=1.6]}] (0,-4) -- (0,1) node[left] {\strut$w$};
\coordinate (A) at (2,-3);
\coordinate (B) at (7,-3);
\coordinate (C) at (7,-1);
\draw[fill=black] (A) circle (2pt) node[below] {\strut$(\bU_L,c^-)$};
\draw[fill=black] (B) circle (2pt) node[below] {\strut$(\hat{\bU},c^-), (\check{\bU},c^+)$};
\draw[fill=black] (C) circle (2pt) node[above] {\strut$(\bU_R,c^+)$};
\draw[thick, midarrow] (A) -- (B);
\draw[dashed, thick, midarrow] (B) -- (C);

\begin{scope}[x=1mm, y=1mm, shift={(-1.7,-94)}]
\node[inner sep=0pt] at (0,0) {\includegraphics[width=80mm]{figure/solution_LIJR_sinh}};
\node[below] at (78.4,0){\strut$\rho$};
\node[left] at (1,58) {\strut$c\,F_{w}$};
\node[left] at (1,20) {\strut$c^{-}F_{w_L}$};
\node[left] at (1,35) {\strut$c^{+}F_{w_L}$};
\node[left] at (1,50) {\strut$c^{+}F_{w_R}$};
\node[below] at (65,0) {\strut$p^{-1}(w_R)$};
\node[below] at (45,0) {\strut$p^{-1}(w_L)$};
\node[right, inner sep=2pt, rotate=60] at (29,43) {\strut$c^{+} \, \rho \, V(\ov)$};
\node[right, inner sep=2pt, rotate=43] at (48,41) {\strut$c^{-} \, \rho \, V(\ov)$};
\node[right, inner sep=2pt, rotate=23] at (56,28) {\strut$c^{+} \, \rho \, V(v(\bU_{R}))$};
\node at (17,8.5) {\strut$\hat{\bU}$};
\node at (32.5,15) {\strut$\check{\bU}$};
\node at (51.5,17.5) {\strut$\bU_R$};
\node at (33.5,1) {\strut$\bU_L$};
\end{scope}

\begin{scope}[shift={(16,0)}]
\draw[-{Latex[scale=1.6]}] (0,0) -- (12.8,0) node[below] {\strut$\ts$};
\draw[-{Latex[scale=1.6]}] (0,-4) -- (0,1) node[left] {\strut$w$};
\coordinate (A) at (8,-3);
\coordinate (B) at (3,-3);
\coordinate (C) at (3,-1);
\draw[fill=black] (A) circle (2pt) node[below] {\strut$(\bU_L,c^-)$};
\draw[fill=black] (B) circle (2pt) node[below] {\strut$(\check{\bU},c^+), (\hat{\bU},c^-)$};
\draw[fill=black] (C) circle (2pt) node[above] {\strut$(\bU_R,c^+)$};
\draw[dashed, thick, midarrow] (A) -- (B);
\draw[dashed, thick, midarrow] (B) -- (C);
\end{scope}

\begin{scope}[x=1mm, y=1mm, shift={(94.3,-94)}]
\node[inner sep=0pt] at (0,0) {\includegraphics[width=80mm]{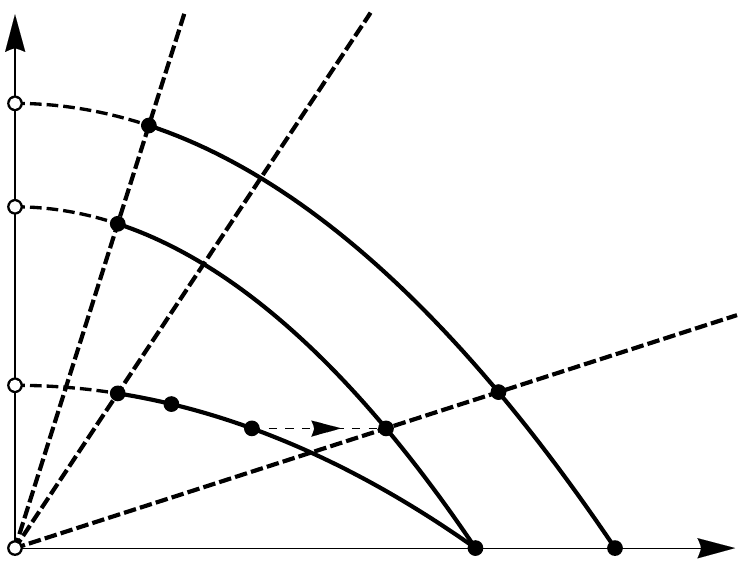}};
\node[inner sep=0pt, below] at (78.4,0){\strut$\rho$};
\node[left] at (1,58) {\strut$c\,F_{w}$};
\node[left] at (1,20) {\strut$c^{-}F_{w_L}$};
\node[left] at (1,38) {\strut$c^{+}F_{w_L}$};
\node[left] at (1,50) {\strut$c^{+}F_{w_R}$};
\node[below] at (65,0) {\strut$p^{-1}(w_R)$};
\node[below] at (45,0) {\strut$p^{-1}(w_L)$};
\node[right, inner sep=2pt] at (15,62) {\strut$c^{+} \, \rho \, V(\ov)$};
\node[right, inner sep=2pt, rotate=58] at (32,42) {\strut$c^{-} \, \rho \, V(\ov)$};
\node[right, inner sep=2pt, rotate=18] at (56,16) {\strut$c^{+} \, \rho \, V(v(\bU_{R}))$};
\node at (25,14) {\strut$\hat{\bU}$};
\node at (39,14) {\strut$\check{\bU}$};
\node at (50.5,18.5) {\strut$\bU_R$};
\node at (17,16) {\strut$\bU_L$};
\end{scope}

\end{tikzpicture}}
\caption{Representation of $\RS_{c^-,c^+}[\bU_L,\bU_R]$ in the $\tbW$-coordinates introduced in~\eqref{e:tbW} (first row), and in the $(\rho,c\,F_{w})$-coordinates (second row), in the case of a $1$-rarefaction (left) and of a $1$-shock (right), followed by a $2$-contact discontinuity.
Here $V(v) = \sinh(v/\vref)$, $c^+ > c^-$, $w_L = w(\bU_L)$ and $w_R = w(\bU_R)$.}
\label{f:Puscifer}
\end{figure}

\item
\label{i:Puscifer}
As already pointed out in the construction of the Riemann solver, see \figurename~\ref{f:Puscifer}, we have that:
\begin{itemize}
\item
Along a $1$-rarefaction, both $v$ and $\ts$ are strictly increasing, while $w$ remains constant.
\item
Across a $1$-shock, both $v$ and $\ts$ are strictly decreasing, while $w$ remains constant.
\item
Across the stationary $c$-wave $(\hat{\bU}, \check{\bU})$, if $c^{+} > c^{-}$ then $v$ is strictly decreasing, if $c^{+} < c^{-}$ then $v$ is strictly increasing. In both cases, $\ts$ and $w$ remain constant.
\item
Across a $2$-contact discontinuity, $v$ and $\ts$ remain constant, whereas $w$ varies.
\end{itemize}
\item
The Riemann solver $\RS_{c^-,c^+}$ is $\Lloc1$-continuous and coherent in the sense of \cite[Definition~3.7]{CorliRazafisonUlrich2025}.
The proof follows standard arguments and is omitted for brevity.
\end{enumerate}
\end{remark}

\subsection{Consistency of the Riemann solver}

In the next proposition we show that if $(\bU_L,\bU_R) \in \calU_{c^{-}}\times\calU_{c^{+}}$, then $\bU(t,x) = \RS_{c^-,c^+}[\bU_L,\bU_R](x/t)$ is a weak solution of the Riemann problem \eqref{e:2x2}, \eqref{e:iRie2}, \eqref{e:cRie2} in the sense of Definition~\ref{d:ws2}.

\begin{proposition}
\label{p:InFlames2a}
Assume \eqref{e:StevenWilson}.
For any $(\bU_L,\bU_R) \in \calU_{c^{-}}\times\calU_{c^{+}}$, the function $\bU(t,x) = \RS_{c^-,c^+}[\bU_L,\bU_R](x/t)$ is a weak solution of the Riemann problem \eqref{e:2x2}, \eqref{e:iRie2}, \eqref{e:cRie2} in the sense of Definition~\ref{d:ws2}.
\end{proposition}

\begin{proof}
By construction the initial condition \eqref{e:ini2} holds almost everywhere.
For this reason we can consider the integral condition \eqref{e:ws2} for test functions $\varphi$ with compact support in $\Dr$.
The claim then follows by observing that any discontinuity $(\bU_-,\bU_+)$ of $\bU(t,x)$ with propagation speed $\sigma \in \R$ satisfies the Rankine-Hugoniot condition \eqref{e:RH}.
The case of a $1$-shock is straightforward, since \eqref{e:RH} holds by construction, see \eqref{e:BringMeTheHorizon}.
For a $2$-contact discontinuity $(\check{\bU},\bU_R)$ with $\sigma = \lambda_{2}(\check{\bU}) = \lambda_{2}(\bU_R) = c^+ \, V(v) > 0$, where $v=v(\check{\bU})=v(\bU_R)$, we have $c(\sigma^\pm) = c^+$, and \eqref{e:RH} easily follows.
Finally, in the case of a stationary discontinuity it is sufficient to recall \eqref{e:Katatonia}.
This concludes the proof.
\end{proof}

In the next proposition we show that if $(\bU_L,\bU_R) \in \calU_{c^{-}}\times\calU_{c^{+}}$, then $\bU(t,x) = \RS_{c^-,c^+}[\bU_L,\bU_R](x/t)$ is an entropy solution in the sense of Definition~\ref{d:entro2} of the Riemann problem \eqref{e:2x2}, \eqref{e:iRie2}, \eqref{e:cRie2}, and that \eqref{e:entro2} is not satisfied by any self-similar weak solution to \eqref{e:2x2}, \eqref{e:iRie2}, \eqref{e:cRie2} that contains $1$-shocks not present in $\RS$-solutions.

\begin{proposition}
\label{p:InFlames2}
Assume \eqref{e:StevenWilson}.
Fix $(\bU_L,\bU_R) \in \calU_{c^{-}}\times\calU_{c^{+}}$ and take $\bU(t,x) = \RS_{c^-,c^+}[\bU_L,\bU_R](x/t)$.
\begin{enumerate}
\item 
$\bU(t,x)$ is an entropy solution of the Riemann problem \eqref{e:2x2}, \eqref{e:iRie2}, \eqref{e:cRie2} in the sense of Definition~\ref{d:entro2}.
\item 
Any $1$-shock $(\bU_{-},\bU_{+})$ with $\rho_->\rho_+$ does not satisfy the entropy condition \eqref{e:entro2}.
\end{enumerate}
\end{proposition}

\begin{proof}
For notational simplicity, we denote $v(\bU_{\pm})$ by $v_{\pm}$, $v(\hat{\bU})$ by $\hat{v}$, and adopt analogous notation for other quantities.
The former claim follows by proving that any discontinuity $(\bU_-,\bU_+)$ of $\bU(t,x)$ with propagation speed $\sigma \in \R$ satisfies $\Upsilon_\kappa(\bU_-,\bU_+,\sigma) \geq 0$ for any $\kappa\in[0,\ov]$, where
\begin{equation}
\label{e:Uk}
\Upsilon_\kappa(\bU_-,\bU_+,\sigma) \doteq 
\begin{cases}
\sigma \left( \mathfrak{E}_\kappa(v_+) - \mathfrak{E}_\kappa(v_-) \right) - \left( \mathfrak{Q}_\kappa(v_+,c(\sigma)) - \mathfrak{Q}_\kappa(v_-,c(\sigma)) \right)
&\hbox{if } \sigma\neq0,
\\
S(\kappa) \left|c^{+}-c^{-}\right| - \left( \mathfrak{Q}_\kappa(v_+,c^+) - \mathfrak{Q}_\kappa(v_-,c^-) \right)
&\hbox{if } \sigma=0.
\end{cases}
\end{equation}
Assume that $\sigma\neq0$.
The entropy production across any $2$-contact discontinuity $(\check{\bU},\bU_R)$ with propagation speed $\sigma = c^{+} \, V(\check{v}) = c^{+} \, V(v_R) > 0$ is zero, namely $\Upsilon_\kappa(\check{\bU},\bU_R,\sigma) = 0$, because $\check{v} = v_{R}$.
Consider the $1$-shock $(\bU_L,\hat{\bU})$.
In this case, $v_L > \hat{v}$ and $\sigma = \sigma(\bU_L,\hat{\bU},c^{-})<0$.
Take $w = w_L = \hat{w}$ and let $\rho_{\kappa} \doteq p^{-1}(w-\kappa)$.
If $\kappa \in (\hat{v} , v_L)$ then $\rho_{\kappa} \in (\rho_L,\hat{\rho})$, hence by Lemma~\ref{l:F} we have
\[F_{w}(\rho_{\kappa}) > F_{w}(\rho_L) + \frac{F_{w}(\hat{\rho})-F_{w}(\rho_L)}{\hat{\rho}-\rho_L} \, (\rho_{\kappa}-\rho_L).\]
By \eqref{e:F}, \eqref{e:Deftones2}, \eqref{e:BringMeTheHorizon} and \eqref{e:TheMarsVolta2}\textsubscript{$1$} we have then
\begin{align*}
\frac{p^{-1}(w)}{c^{-}} \, \Upsilon_\kappa\left(\bU_L,\hat{\bU},\sigma\right)
&=
\frac{F_{w}(\hat{\rho})-F_{w}(\rho_L)}{\hat{\rho}-\rho_L} \left( \hat{\rho} - 2\rho_{\kappa} + \rho_{L} \right) 
+2F_{w}(\rho_{\kappa})-F_{w}(\hat{\rho})-F_{w}(\rho_{L})
\\&=
2 \left(
F_{w}(\rho_{\kappa}) - F_{w}(\rho_L) 
- \frac{F_{w}(\hat{\rho})-F_{w}(\rho_L)}{\hat{\rho}-\rho_L} \left(\rho_{\kappa} - \rho_{L}\right)
\right) > 0.
\end{align*}
Furthermore, we also have
\begin{align*}
\kappa \leq \hat{v}:\quad&
\frac{p^{-1}(w)}{c^{-}} \, \Upsilon_\kappa(\bU_L,\hat{\bU},\sigma) = 
\frac{F_{w}(\hat{\rho}) - F_{w}(\rho_L)}{\hat{\rho}-\rho_L} \left( \rho_L - \hat{\rho} \right) 
+ F_{w}(\hat{\rho}) - F_{w}(\rho_L) = 0,
\\
\kappa \geq v_L:\quad& 
\frac{p^{-1}(w)}{c^{-}} \, \Upsilon_\kappa(\bU_L,\hat{\bU},\sigma) =
\frac{F_{w}(\hat{\rho}) - F_{w}(\rho_L)}{\hat{\rho}-\rho_L} \left( \hat{\rho} - \rho_L \right) 
+ F_{w}(\rho_L) - F_{w}(\hat{\rho}) = 0.
\end{align*}

We now assume $\sigma=0$ and show that
\[\Upsilon_\kappa(\bU_-,\bU_+,0) =
S(\kappa) \left|c^{+}-c^{-}\right| - \left( \mathfrak{Q}_\kappa(v_+,c^+) - \mathfrak{Q}_\kappa(v_-,c^-) \right) \geq 0.\]
If $v_\pm = 0$, then for any $\kappa \in [0,\ov]$ we have $\Upsilon_\kappa(\bU_-,\bU_+,0) = S(\kappa) \left( \left|c^{+}-c^{-}\right| + c^{+}-c^{-} \right) \geq 0$ because $S(0)=0$.
Assume that $v_\pm \neq 0$.
In this case, by \eqref{e:Slipknot} and \eqref{e:pezza} we have $w_- = w_+$; denote by $w$ this common value.
We distinguish the following cases:
\begin{itemize}
\item 
If $\kappa \in (v_-,v_+)$ then $c^+ < c^-$ and by \eqref{e:Katatonia} we have
\begin{align*}
p^{-1}(w)\,\Upsilon_\kappa(\bU_-,\bU_+,0) &=
\left(c^{-}-c^{+}\right) F_{w}(\rho_{\kappa}) + c^+ \left(F_{w}(\rho_+)-F_{w}(\rho_{\kappa})\right) - c^- \left(F_{w}(\rho_{\kappa})-F_{w}(\rho_-)\right)
\\&=
2c^+ \left(F_{w}(\rho_+)-F_{w}(\rho_{\kappa})\right)
> 0.
\end{align*}
\item 
If $\kappa \in (v_+,v_-)$ then $c^- < c^+$ and by \eqref{e:Katatonia} we have
\begin{align*}
p^{-1}(w)\,\Upsilon_\kappa(\bU_-,\bU_+,0) &=
\left(c^{+}-c^{-}\right) F_{w}(\rho_{\kappa}) + c^+ \left(F_{w}(\rho_{\kappa})-F_{w}(\rho_+)\right) - c^- \left(F_{w}(\rho_-)-F_{w}(\rho_{\kappa})\right)
\\&=
2c^+ \left(F_{w}(\rho_{\kappa})-F_{w}(\rho_+)\right)
> 0.
\end{align*}
\item 
If $\kappa \leq \min\{v_-,v_+\}$ then by \eqref{e:Katatonia} we have
\begin{align*}
p^{-1}(w)\,\Upsilon_\kappa(\bU_-,\bU_+,0) &=
\left|c^{+}-c^{-}\right| F_{w}(\rho_{\kappa}) + c^+ \left(F_{w}(\rho_+)-F_{w}(\rho_{\kappa})\right) - c^- \left(F_{w}(\rho_-)-F_{w}(\rho_{\kappa})\right)
\\&=
\left(\left|c^{+}-c^{-}\right|-c^{+}+c^{-}\right) F_{w}(\rho_{\kappa})
\geq 0.
\end{align*}
\item 
If $\kappa \geq \max\{v_-,v_+\}$ then by \eqref{e:Katatonia} we have
\begin{align*}
p^{-1}(w)\,\Upsilon_\kappa(\bU_-,\bU_+,0) &=
\left|c^{+}-c^{-}\right| F_{w}(\rho_{\kappa}) + c^+ \left(F_{w}(\rho_{\kappa})-F_{w}(\rho_+)\right) - c^- \left(F_{w}(\rho_{\kappa})-F_{w}(\rho_-)\right)
\\&=
\left(\left|c^{+}-c^{-}\right|+c^{+}-c^{-}\right) F_{w}(\rho_{\kappa})
\geq 0.
\end{align*}
\end{itemize}
This concludes the proof of the former claim.

To prove the latter claim, consider a $1$-shock $(\bU_{-},\bU_{+})$ with $\rho_->\rho_+$ and $\sigma = \sigma(\bU_-,\bU_{+},c^{-})<0$.
In this case, take $w=w_{-}=w_{+}$ and observe that for $\kappa \in (v_- , v_{+})$
\begin{align*}
\frac{p^{-1}(w)}{c^{-}} \, \Upsilon_\kappa(\bU_-,\bU_+,\sigma) &= 
\frac{F_{w}(\rho_{+}) - F_{w}(\rho_-)}{\rho_{+}-\rho_-} \left( 2\rho_{\kappa} - \rho_{+} - \rho_{-} \right)
+ F_{w}(\rho_{+}) - 2F_{w}(\rho_{\kappa}) + F_{w}(\rho_{-})
\\&=
-2 \left(
F_{w}(\rho_{\kappa}) - F_{w}(\rho_-) 
- \frac{F_{w}(\rho_{+}) - F_{w}(\rho_-)}{\rho_{+}-\rho_-} \left(\rho_{\kappa} - \rho_{-}\right)
\right) < 0.
\end{align*}
This concludes the proof.
\end{proof}

\section{Proof of Theorem~\ref{t:uniqueness}}
\label{s:proofuni}

This section is devoted to the Lipschitz estimate with respect to the initial data of the $v$-component of entropy solutions satisfying \eqref{e:Lipesti}.
Recall that $L \doteq \vref \, V'(0) = \max_{v \in [0,\ov]} \left| \vref \, V'(v) - V(v) \right| > 0$.
For simplicity, we denote $v(\bU^{*})$ and $v(\bU^*_{o})$ by $v^{*}$ and $v^{*}_{o}$, respectively, where $*\,\in\{a,b\}$.
Take $\oc \doteq \max_{x\in\R} c(x)$.

\begin{lemma}
\label{l:Moodring}
For every $v^{a},v^{b}\in[0,\ov]$ we have
\[
\left|S(v^{b})-S(v^{a})\right|
\leq
L \left|e^{-v^{b}/\vref}-e^{-v^{a}/\vref}\right|.
\]
\end{lemma}

\begin{proof}
Set $r^{*} \doteq e^{-v^{*}/\vref}$ for $*\,\in\{a,b\}$.
Then
\[\left|S(v^{b})-S(v^{a})\right| = \left|V\left(-\vref\,\ln(r^{b})\right) r^{b} - V\left(-\vref\,\ln(r^{a})\right) r^{a}\right| \leq L \left| r^{b} - r^{a} \right|\]
because, by \eqref{e:V}\textsubscript{$2$}, for any $r \in [e^{-\ov/\vref},1]$ we have
\[\left|\frac{\d}{\d r} \left( V\left(-\vref\,\ln(r)\right) r \right)\right| = \vref \, V'\left(-\vref\,\ln(r)\right) - V\left(-\vref\,\ln(r)\right) \leq L.
\]
This concludes the proof.
\end{proof}

\begin{proposition}
The entropy solutions $\bU_{a}$ and $\bU_{b}$ corresponding to two initial data $\bU_{o}^{a}$ and $\bU_{o}^{b}$ satisfy the following estimate for any $T>0$ and $R>0$:
\[\int_{-R}^{R} \left|e^{-v^{b}(T,x)/\vref} - e^{-v^{a}(T,x)/\vref}\right| \d x \leq \int_{-R-T \, \oc \, L}^{R+T \, \oc \, L} \left|e^{-v^{b}(0,x)/\vref} - e^{-v^{a}(0,x)/\vref}\right| \d x.\]
\end{proposition}

\begin{proof}
Take $\kappa^{a}, \kappa^{b} \in [0,\ov]$ and a smooth function $\varphi = \varphi(t,x,s,y) \geq 0$ with compact support in the set $\D^{2}$.
By Definition~\ref{d:entro2}, the entropy solutions $\bU_{a}(t,x)$ and $\bU_{b}(s,y)$ satisfy condition \eqref{e:entro2}:
\begin{align*}
\iint_{\D} \left( \mathfrak{E}_{\kappa^{a}}(v^{a}) \, \partial_t\varphi + \mathfrak{Q}_{\kappa^{a}}(v^{a},c) \, \partial_x\varphi \right) \d x\, \d t
+ \int_{\R} \mathfrak{E}_{\kappa^{a}}\left(v_{o}^{a}\right) \, \varphi(0,x,s,y) \, \d x
\\
+ S(\kappa^{a}) \sum_{i=0}^{I_{c}} \left|c(\xi_{i}^{+})-c(\xi_{i}^{-})\right| \int_0^{+\infty} \varphi(t,\xi_{i},s,y) \, \d t
&\geq 0,
\\
\iint_{\D} \left( \mathfrak{E}_{\kappa^{b}}(v^{b}) \, \partial_s\varphi + \mathfrak{Q}_{\kappa^{b}}(v^{b},c) \, \partial_y\varphi \right) \d y\, \d s
+ \int_{\R} \mathfrak{E}_{\kappa^{b}}\left(v_{o}^{b}\right) \, \varphi(t,x,0,y) \, \d y
\\
+ S(\kappa^{b}) \sum_{i=0}^{I_{c}} \left|c(\xi_{i}^{+})-c(\xi_{i}^{-})\right| \int_0^{+\infty} \varphi(t,x,s,\xi_{i}) \, \d s
&\geq 0.
\end{align*}

Set $\kappa^{a} = v^{b}(s,y)$ in the former inequality and integrate with respect to $s$ and $y$.
Similarly, set $\kappa^{b} = v^{a}(t,x)$ in the latter inequality and integrate with respect to $t$ and $x$.
Then add the two results and get
\begin{align*}
\iiiint_{\D^{2}} \left[ \left|e^{-v^{b}/\vref} - e^{-v^{a}/\vref}\right| \left(\partial_t\varphi+\partial_s\varphi\right) 
- \left| S(v^{b}) - S(v^{a}) \right| \left( c(x) \, \partial_x\varphi + c(y) \, \partial_y\varphi \right) \right] \d x\, \d t \, \d y \, \d s
\\
+ \iint_{\D} \int_{\R} \left|e^{-v^{b}/\vref} - e^{-v^{a}_{o}/\vref}\right| \, \varphi(0,x,s,y) \, \d x \, \d y \, \d s
+ \int_{\R} \iint_{\D} \left|e^{-v^{a}/\vref} - e^{-v^{b}_{o}/\vref}\right| \, \varphi(t,x,0,y) \, \d t \, \d x \, \d y
\\
+ \sum_{i=0}^{I_{c}} \left|c(\xi_{i}^{+})-c(\xi_{i}^{-})\right| \left[ \iint_{\D} \int_0^{+\infty} S(v^{b}) \, \varphi(t,\xi_{i},s,y) \, \d t \, \d s \, \d y
+ \iint_{\D} \int_0^{+\infty} S(v^{a}) \, \varphi(t,x,s,\xi_{i}) \, \d s \, \d x \, \d t \right]
&\geq 0.
\end{align*}

For $\eta>0$ sufficiently small, we introduce the Lipschitz function
\[\Theta_{\eta}(x) \doteq
\left\{\begin{array}{@{}l@{\quad}l@{}}
-\frac{x-\xi_{i}}{\eta}-2&\hbox{if } x\in\bigcup_{i=1}^{I_{c}} (\xi_{i}-3\eta,\xi_{i}-2\eta)
\\[5pt]
0&\hbox{if } x\in\bigcup_{i=1}^{I_{c}}[\xi_{i}-2\eta,\xi_{i}+2\eta]
\\[5pt]
\frac{x-\xi_{i}}{\eta}-2&\hbox{if } x\in\bigcup_{i=1}^{I_{c}} (\xi_{i}+2\eta,\xi_{i}+3\eta)
\\[5pt]
1&\hbox{otherwise}.
\end{array}\right.\]
Observe that $\Theta_{\eta}$ vanishes in a neighbourhood of $x=\xi_{i}$ for any $i\in\{1,\ldots,I_{c}\}$ and $\|\Theta_{\eta} - 1\|_{\L1(\R)} \to 0$ as $\eta\to0$.
Let $\delta_{\varepsilon}$, with $\varepsilon\in(0,\eta)$, be a smooth approximation of the Dirac mass at the origin.
More precisely, given a $\C\infty$-function $\delta \colon \R \to [0,1]$ such that $\int_{\R} \delta(x) \, \d x = 1$ and $\delta(x) = 0$ for all $|x|\geq1$, we take $\delta_{\varepsilon}(z) = \delta(z/\varepsilon)/\varepsilon$.
Introduce the Lipschitz test function
\[\varphi(t,x,s,y) = \Phi\left(\frac{t+s}{2},\frac{x+y}{2}\right) \delta_{\varepsilon}\left(\frac{t-s}{2}\right) \delta_{\varepsilon}\left(\frac{x-y}{2}\right) \Theta_{\eta}(x) \, \Theta_{\eta}(y),\]
where $\Phi = \Phi(T,X) \geq 0$ is a smooth function with compact support in the half plane $T>0$.
For completeness, we note that, by $\Lloc1$-continuity, the entropy inequality extends to non-negative Lipschitz test functions with compact support.
For any $(t,x,s,y) \in \D^{2}$ we have
\begin{align*}
\partial_t\varphi + \partial_s\varphi
={}&
\partial_{T}\Phi\left(\frac{t+s}{2},\frac{x+y}{2}\right) \delta_{\varepsilon}\left(\frac{t-s}{2}\right) \delta_{\varepsilon}\left(\frac{x-y}{2}\right) \Theta_{\eta}(x) \, \Theta_{\eta}(y),
\\
c(x) \, \partial_x\varphi + c(y) \, \partial_y\varphi
={}&
c(x) \left( \partial_x\varphi + \partial_y\varphi \right) + \left( c(y) - c(x) \right) \partial_y\varphi
\\={}&
c(x) \Bigl[
\begin{aligned}[t]
&\partial_{X}\Phi\left(\frac{t+s}{2},\frac{x+y}{2}\right) \Theta_{\eta}(x) \, \Theta_{\eta}(y)
\\
+{}&\Phi\left(\frac{t+s}{2},\frac{x+y}{2}\right) \left( \Theta_{\eta}'(x) \, \Theta_{\eta}(y) + \Theta_{\eta}(x) \, \Theta_{\eta}'(y) \right) \Bigr] \, \delta_{\varepsilon}\left(\frac{t-s}{2}\right) \delta_{\varepsilon}\left(\frac{x-y}{2} \right) .
\end{aligned}
\end{align*}
In the last equality above, we used the fact that $\left( c(y) - c(x) \right) \partial_y\varphi$ vanishes.
Indeed, if $\left|\frac{x-y}{2}\right| \geq \varepsilon$ then both $\delta_{\varepsilon}\left(\frac{x-y}{2}\right)$ and $\delta_{\varepsilon}'\left(\frac{x-y}{2}\right)$ vanish, and consequently $\partial_y\varphi=0$.
Assume then $|x-y| < 2\varepsilon$.
In this case, if $|y-\xi_{i}| < 2\eta$ then both $\Theta_{\eta}(y)$ and $\Theta_{\eta}'(y)$ vanish, and hence $\partial_y\varphi=0$.
Assume therefore that $|x-y| < 2\varepsilon$ and $|y-\xi_{i}| \geq 2\eta$ for any $i\in\{1,\ldots,I_{c}\}$. In this case, there exists $i\in\{0,\ldots,I_{c}\}$ such that both $x$ and $y$ belong to $\Xi_{i}$ because $\varepsilon<\eta$, and therefore $c(x)=c(y)$.
This concludes the proof of the claim.

By exploiting the above expressions we deduce that
\begin{align*}
\iiiint_{\D^{2}} \Bigl\{
\left[ \left|e^{-\frac{v^{b}}{\vref}} - e^{-\frac{v^{a}}{\vref}}\right| \partial_{T}\Phi\left(\frac{t+s}{2},\frac{x+y}{2}\right)
- \left| S(v^{b}) - S(v^{a}) \right| c(x) \, \partial_{X}\Phi\left(\frac{t+s}{2},\frac{x+y}{2}\right)  \right] \Theta_{\eta}(x) \, \Theta_{\eta}(y)
\\
- \left| S(v^{b}) - S(v^{a}) \right| c(x) \, \Phi\left(\frac{t+s}{2},\frac{x+y}{2}\right) \left( \Theta_{\eta}'(x) \, \Theta_{\eta}(y) + \Theta_{\eta}(x) \, \Theta_{\eta}'(y) \right)
\Bigr\} \, \delta_{\varepsilon}\left(\frac{t-s}{2}\right) \delta_{\varepsilon}\left(\frac{x-y}{2} \right) \d x\, \d t \, \d y \, \d s 
\\
+ \iint_{\D} \int_{\R} \left|e^{-v^{b}/\vref} - e^{-v^{a}_{o}/\vref}\right| \, \Phi\left(\frac{s}{2},\frac{x+y}{2}\right) \delta_{\varepsilon}\left(-\frac{s}{2}\right) \delta_{\varepsilon}\left(\frac{x-y}{2}\right) \Theta_{\eta}(x) \, \Theta_{\eta}(y) \, \d x \, \d y \, \d s
\\
+ \int_{\R} \iint_{\D} \left|e^{-v^{a}/\vref} - e^{-v^{b}_{o}/\vref}\right| \, \Phi\left(\frac{t}{2},\frac{x+y}{2}\right) \delta_{\varepsilon}\left(\frac{t}{2}\right) \delta_{\varepsilon}\left(\frac{x-y}{2}\right) \Theta_{\eta}(x) \, \Theta_{\eta}(y) \, \d t \, \d x \, \d y
&\geq 0.
\end{align*}
We let $\varepsilon>0$ go to zero and obtain
\begin{align}
\label{e:Ghost}
&\iint_{\Dr}
\left[ \left|e^{-v^{b}(t,x)/\vref} - e^{-v^{a}(t,x)/\vref}\right| \partial_{t}\Phi(t,x)
- \left| S\left(v^{b}(t,x)\right) - S\left(v^{a}(t,x)\right) \right| c(x) \, \partial_{x}\Phi(t,x) \right] \Theta_{\eta}(x)^{2} \, \d x\, \d t
\\
\geq{}& 2 \iint_{\Dr} \left| S\left(v^{b}(t,x)\right) - S\left(v^{a}(t,x)\right) \right| c(x) \, \Phi(t,x) \, \Theta_{\eta}(x) \, \Theta_{\eta}'(x) \, \d x\, \d t
\nonumber
\end{align}
because $\Phi$ has compact support in the half plane $T>0$.
The second double integral above writes
\begin{align*}
&\frac{2}{\eta} \int_{0}^{+\infty} \sum_{i=1}^{I_{c}} \Biggl[
\begin{aligned}[t]
&\int_{\xi_{i}-3\eta}^{\xi_{i}-2\eta} \left| S\left(v^{b}(t,x)\right) - S\left(v^{a}(t,x)\right) \right| c(x) \, \Phi(t,x) \, \left(\frac{x-\xi_{i}}{\eta}+2\right) \d x
\\+{}
&\int_{\xi_{i}+2\eta}^{\xi_{i}+3\eta} \left| S\left(v^{b}(t,x)\right) - S\left(v^{a}(t,x)\right) \right| c(x) \, \Phi(t,x) \left(\frac{x-\xi_{i}}{\eta}-2\right) \d x
\Biggr] \d t
\end{aligned}
\\
={}&2 \int_{0}^{+\infty} \sum_{i=1}^{I_{c}} \Biggl[
\begin{aligned}[t]
&\int_{2}^{3} \left| S\left(v^{b}(t,\xi_{i}-\eta\,z)\right) - S\left(v^{a}(t,\xi_{i}-\eta\,z)\right) \right| c(\xi_{i}-\eta\,z) \, \Phi(t,\xi_{i}-\eta\,z) \, (2-z) \, \d z
\\+{}
&\int_{2}^{3} \left| S\left(v^{b}(t,\xi_{i}+\eta\,z)\right) - S\left(v^{a}(t,\xi_{i}+\eta\,z)\right) \right| c(\xi_{i}+\eta\,z) \, \Phi(t,\xi_{i}+\eta\,z) \, (z-2) \, \d z
\Biggr] \d t
\end{aligned}
\end{align*}
where in the above last equality we take $x=\xi_{i}-\eta\,z$ and $x=\xi_{i}+\eta\,z$.
We now let $\eta>0$ go to zero and get
\begin{align*}
&2 \int_{0}^{+\infty} \sum_{i=1}^{I_{c}} \Biggl[
\begin{aligned}[t]
&\int_{2}^{3} \left| S\left(v^{b}(t,\xi_{i}^{-})\right) - S\left(v^{a}(t,\xi_{i}^{-})\right) \right| c(\xi_{i}^{-}) \, (2-z) \, \d z
\\+{}
&\int_{2}^{3} \left| S\left(v^{b}(t,\xi_{i}^{+})\right) - S\left(v^{a}(t,\xi_{i}^{+})\right) \right| c(\xi_{i}^{+}) \, (z-2) \, \d z
\Biggr] \, \Phi(t,\xi_{i}) \, \d t
\end{aligned}
\\={}&
\int_{0}^{+\infty} \sum_{i=1}^{I_{c}} \biggl[
\begin{aligned}[t]
&-\left| \ts\left(v^{b}(t,\xi_{i}^{-}),c(\xi_{i}^{-})\right) - \ts\left(v^{a}(t,\xi_{i}^{-}),c(\xi_{i}^{-})\right) \right|
\\
&+ \left| \ts\left(v^{b}(t,\xi_{i}^{+}),c(\xi_{i}^{+})\right) - \ts\left(v^{a}(t,\xi_{i}^{+}),c(\xi_{i}^{+})\right) \right|
\biggr] \, \Phi(t,\xi_{i}) \, \d t = 0
\end{aligned}
\end{align*}
because the quantity inside the square brackets vanishes.
Indeed, by Remark~\ref{r:traces} \ref{i:Puscifer} we have $\ts\left(v^{*}(t,\xi_{i}^{-}),c(\xi_{i}^{-})\right) = \ts\left(v^{*}(t,\xi_{i}^{+}),c(\xi_{i}^{+})\right)$ for $*\,\in\{a,b\}$.

The first double integral \eqref{e:Ghost} converges, by dominated convergence, to
\[\iint_{\Dr}
\left[ \left|e^{-v^{b}(t,x)/\vref} - e^{-v^{a}(t,x)/\vref}\right| \partial_{t}\Phi(t,x)
- \left| S\left(v^{b}(t,x)\right) - S\left(v^{a}(t,x)\right) \right| c(x) \, \partial_{x}\Phi(t,x) \right] \d x\, \d t.\]
Hence we obtain
\[
\iint_{\Dr} \left[
\left|e^{-v^{b}/\vref} - e^{-v^{a}/\vref}\right| \partial_{t}\Phi - \left| S(v^{b}) - S(v^{a}) \right| c(x) \, \partial_{x}\Phi \right] \d x\, \d t
\geq 0.
\]
For $\varepsilon>0$ sufficiently small, we take
\[\Phi_{\varepsilon}(t,x) \doteq \left( \int_{x-R-(T-t) \, \oc \, L+\varepsilon}^{x+R+(T-t) \, \oc \, L} \delta_{\varepsilon}(y) \, \d y \right) \left( \int_{t-T+\varepsilon}^{t} \delta_{\varepsilon}(s) \, \d s \right)\]
in the above estimate and let $\varepsilon$ go to zero.
We then obtain
\[\int_{-R-T \, \oc \, L}^{R+T \, \oc \, L} \left|e^{-v^{b}(0,x)/\vref} - e^{-v^{a}(0,x)/\vref}\right| \d x 
- \int_{-R}^{R} \left|e^{-v^{b}(T,x)/\vref} - e^{-v^{a}(T,x)/\vref}\right| \d x \geq \int_{0}^{T} I(t) \, \d t,\]
where
\begin{align*}
I(t) \doteq{}&
\oc \, L \left|e^{-v^{b}(t,-R-(T-t)\,\oc \, L)/\vref} - e^{-v^{a}(t,-R-(T-t)\,\oc \, L)/\vref}\right|
\\&
+ \oc \, L \left|e^{-v^{b}(t,R+(T-t)\,\oc \, L)/\vref} - e^{-v^{a}(t,R+(T-t)\,\oc \, L)/\vref}\right|
\\&
+ \left| S\left(v^{b}(t,-R-(T-t)\,\oc \, L)\right) - S\left(v^{a}(t,-R-(T-t)\,\oc \, L)\right) \right| c\left(-R-(T-t)\,\oc \, L\right) 
\\&
- \left| S\left(v^{b}(t,R+(T-t)\,\oc \, L)\right) - S\left(v^{a}(t,R+(T-t)\,\oc \, L)\right) \right| c\left(R+(T-t)\,\oc \, L\right)
\end{align*}
is non-negative by Lemma~\ref{l:Moodring}.
This completes the proof.
\end{proof}

At last, to conclude the proof of Theorem~\ref{t:uniqueness} and get the estimate \eqref{e:Lipesti} it is sufficient to observe that
\[\frac{1}{\vref} \, e^{-\ov/\vref}  \left|v^{b}-v^{a}\right| \leq \left|e^{-v^{b}/\vref}-e^{-v^{a}/\vref}\right| \leq \frac{1}{\vref} \left|v^{b}-v^{a}\right|.\]

\section{Proof of Theorem~\ref{t:LornaShore2}}
\label{s:proof}

In this section, we prove Theorem~\ref{t:LornaShore2}.
To this end, we apply the wave-front tracking algorithm to construct a sequence $\{\bU^{n}\}_n$ that converges to an entropy solution of the Cauchy problem \eqref{e:2x2}, \eqref{e:ini2}.
The estimates needed to establish compactness are derived in the $\tbW$-coordinate, see \eqref{e:tbW}. 
This motivates the introduction of the following notation
\begin{equation}
\label{e:tTV}
\tTV(\bU,c) \doteq \TV(\tbW(\bU,c)) = \TV(\ts(\bU,c)) + \TV(w(\bU)).
\end{equation}

\begin{remark}
The introduction of $\tTV$ has the following motivation.
The total variation of the approximate solution, when expressed in the conservative variables $\bU$ or in the Riemann invariant coordinates $\bW$, may well increase after an interaction, see Remark~\ref{r:Malmsteen} below.
On the contrary, the total variation of the approximate solution expressed in the $\tbW$-coordinates does not increase after any interaction, see Theorem~\ref{t:Deftones}.
\end{remark}

\subsection{Approximation of \texorpdfstring{$\bU_o$}{}}

We now approximate the initial datum $\bU_o$ with a $\PC$-function $\bU_o^n$.
Here $\PC$ is the space of right-continuous piecewise constant functions having a finite number of discontinuities (and taking a finite number of values).
Recall that $c \colon \R \to [\uc,+\infty)$ takes the form \eqref{e:def_c}
\[c(x) = \sum_{i = 0}^{I_{c}} c_{i+1/2} \, \mathbbm{1}_{\Xi_{i}}(x),\]
for some constants $c_{i+1/2} \geq \uc > 0$ such that $c_{i-1/2} \neq c_{i+1/2}$, and 
\begin{align*}
&\Xi_{0} \doteq (-\infty,\xi_{1}),&
&\Xi_{i} \doteq [\xi_i,\xi_{i+1}),\ i \in \{1,\ldots,I_{c}-1\},&
&\Xi_{I_{c}} \doteq [\xi_{I_{c}},+\infty),
\end{align*}
with $\xi_i < \xi_{i+1}$.

By hypotheses, the constants $\uc, \oV>0$ satisfy \eqref{e:cV}, namely $\oV \leq \uc \, S(\ov)$.
Fix $n \in \N$ sufficiently large, and take $\varepsilon_{1}^{n} \doteq \oV/2^{n}$ and $\varepsilon_{2}^{n} \doteq -\uw/2^{n}$.
On each strip $[0,+\infty) \times \Xi_{i}$ we introduce the grid
\begin{equation}
\label{e:NineInchNails}
\calU_{i+1/2}^{n} \doteq \left\{ \bU \in \calU_{c_{i+1/2}} : \tbW(\bU,c_{i+1/2}) \in (\varepsilon_{1}^{n}\Z \cap [0,\oV]) \times (\varepsilon_{2}^{n}\Z \cap [\uw,0]) \right\}
\end{equation}
with $\calU_{c}$ defined in \eqref{e:dominva}, and $\tbW$ in \eqref{e:tbW}.
Approximate $\bU_o$ with a piecewise constant function $\bU_o^{n} \in \PC(\R;\calU)$ such that
\begin{equation}
\label{e:TOOL}
x\in\Xi_{i} \Longrightarrow \bU_o^{n}(x) \in \calU_{i+1/2}^{n},
\end{equation}
and
\begin{align}
\label{e:uoncond1}
&\|\tbW(\bU^{n}_o,c)\|_{\L\infty(\R)} \leq \|\tbW(\bU_o,c)\|_{\L\infty(\R)},&
&\|\tbW(\bU_o^{n},c) - \tbW(\bU_o,c)\|_{\L\infty(\R)} \leq \frac{1}{n},
\\
\label{e:uoncond2}
&\tTV(\bU_o^{n},c) \leq \tTV(\bU_o,c),&
&\int_{-n}^n \|\bU_o^{n}(x)-\bU_o(x)\| \, \d x \leq \frac{1}{n},
\\
\label{e:uoncond3}
&\lim_{n\to+\infty} \lim_{x\to-\infty} \bU_o^{n}(x) = \lim_{x\to-\infty} \bU_o(x),&
&\lim_{n\to+\infty} \lim_{x\to+\infty} \bU_o^{n}(x) = \lim_{x\to+\infty} \bU_o(x),
\end{align}
where $\tTV(\bU,c)$ is defined in \eqref{e:tTV}.

\subsection{Approximate Riemann solvers}
\label{s:Karnivool}

Consider the Riemann solver $\RS_{c^-,c^+}$ corresponding to $c^-,c^+ > 0$ and described in Definition~\ref{d:RS}.
We aim to introduce approximate Riemann solvers
\begin{align*}
&\RS_{i-1/2}^{n} \approx \RS_{c_{i-1/2},c_{i-1/2}},&
&\RS_i^{n} \approx \RS_{c_{i-1/2},c_{i+1/2}},
\end{align*}
such that
\begin{equation}
\label{e:Polyphia}
\begin{aligned}
(\bU_{L},\bU_{R}) \in \calU_{i-1/2}^{n} \times \calU_{i-1/2}^{n}
&\Longrightarrow
\RS_{i-1/2}^{n}[\bU_{L},\bU_{R}](\nu) \in \calU_{i-1/2}^{n} \hbox{ for all }\nu\in\R,
\\
(\bU_{L},\bU_{R}) \in \calU_{i-1/2}^{n} \times \calU_{i+1/2}^{n}
&\Longrightarrow
\RS_i^{n}[\bU_{L},\bU_{R}](\nu) \in \calU_{i\pm1/2}^{n} \hbox{ for all }\pm\nu>0.
\end{aligned}
\end{equation}
We first establish the following result.

\begin{lemma}
The following properties hold.
\begin{itemize}
\item
If $(\bU_{L}, \bU_{R}) \in \calU_{i-1/2}^{n} \times \calU_{i-1/2}^{n}$, then $\hat{\bU}= \check{\bU} \in \calU_{i-1/2}^{n}$, where $\hat{\bU}, \check{\bU} \in \calU$ are implicitly defined by \eqref{e:Avatar} with $c^{\pm} = c_{i-1/2}$.
\item
If $(\bU_{L}, \bU_{R}) \in \calU_{i-1/2}^{n} \times \calU_{i+1/2}^{n}$, then $(\hat{\bU}, \check{\bU}) \in \calU_{i-1/2}^{n} \times \calU_{i+1/2}^{n}$, where $\hat{\bU}, \check{\bU} \in \calU$ are implicitly defined by \eqref{e:Avatar} with $c^{\pm} = c_{i\pm1/2}$.
\end{itemize}
\end{lemma}
\begin{proof}
We prove the second property, the first follows by analogous arguments.
By \eqref{e:Leprous} we have
\begin{align*}
\ts(\hat{\bU},c_{i-1/2}) = \ts(\check{\bU},c_{i+1/2}) = \ts(\bU_{R},c_{i+1/2}) &\in \varepsilon_{1}^{n}\Z \cap [0,\oV],&
w(\hat{\bU}) = w(\check{\bU}) = w(\bU_L) &\in \varepsilon_{2}^{n}\Z \cap [\uw,0].
\end{align*}
This concludes the proof.
\end{proof}
The above lemma leads us to discretise the $1$-rarefaction waves.
More precisely, if $(\bU_{L},\bU_{R}) \in \calU_{i-1/2}^{n} \times \calU_{i+1/2}^{n}$ is such that $\RS_{c_{i-1/2},c_{i+1/2}}[\bU_{L},\bU_{R}]$ has a $1$-rarefaction, namely, if $(\hat{\bU},\check{\bU}) \in \calU_{i-1/2}^{n} \times \calU_{i+1/2}^{n}$ satisfies
\[
\ts(\hat{\bU},c_{i-1/2}) = \ts(\bU_L,c_{i-1/2}) + J\,\varepsilon_{1}^{n},\quad J \in \N,
\]
then we let 
\begin{equation*}
\RS_i^{n}[\bU_L,\bU_R](\nu) =
\left\{\begin{array}{@{}l@{\quad\hbox{if }}l}
\bU_L & \nu < \sigma_{i-1/2}^{n}(\bU_L,\bU_{1/2}),
\\
\bU_{j+1/2} & 
\begin{aligned}[t]
\sigma_{i-1/2}^{n}(\bU_{j-1/2},\bU_{j+1/2}) < \nu \leq \sigma_{i-1/2}^{n}(\bU_{j+1/2},\bU_{j+3/2}),\ j\in\{0,\ldots,J-2\},
\end{aligned}
\\
\hat{\bU} & \sigma_{i-1/2}^{n}(\bU_{J-3/2},\hat{\bU}) \leq \nu < 0,
\\
\check{\bU} & 0 \leq \nu < c_{i+1/2} \, V(v(\bU_{R})),
\\
\bU_{R} & \nu \geq c_{i+1/2} \, V(v(\bU_{R})),
\end{array}\right.
\end{equation*}
where $\bU_{j-1/2} \in \calU_{i-1/2}^{n}$, $j\in\{0,\ldots,J\}$, are implicitly defined by
\begin{align}
\label{e:Plini}
&\ts(\bU_{j-1/2},c_{i-1/2}) = \ts(\bU_{L},c_{i-1/2}) + j \, \varepsilon_{1}^{n},&
&w(\bU_{j-1/2}) = w(\bU_{L}),
\end{align}
so that $\bU_{-1/2} = \bU_{L}$ and $\bU_{J-1/2} = \hat{\bU}$, and
\begin{equation}
\label{e:Haken}
\sigma_{i-1/2}^{n}(\bU_{a},\bU_{b}) \doteq c_{i-1/2} \, \frac{F_{w(\bU_{b})}(\rho_{b}) - F_{w(\bU_{a})}(\rho_{a})}{\rho_{b}-\rho_{a}}.
\end{equation}
Observe that $\lambda_{1}(\bU_{j-1/2},c_{i-1/2}) < \sigma_{i-1/2}^{n}(\bU_{j-1/2},\bU_{j+1/2}) < \lambda_{1}(\bU_{j+1/2},c_{i-1/2})$, moreover \eqref{e:Plini}, \eqref{e:TheMarsVolta}\textsubscript{$1$} imply that
\[\rho_{j-1/2} = F_{w(\bU_L)}^{-1}\left( \frac{\ts(\bU_L,c_{i-1/2})+j\,\varepsilon_{1}^{n}}{c_{i-1/2}} \, p^{-1}\left(w(\bU_L)\right) \right).\]

The construction of $\RS_{i-1/2}^{n}$ is analogous and is therefore omitted.

\subsection{Approximate solution}

For sufficiently small times, the approximate solution $\bU^{n}$ is obtained by gluing together the solutions computed using approximate Riemann solvers introduced in Section~\ref{s:Karnivool}.
More precisely, we first apply $\RS_i^{n}$ at the points $(t,x) = (0,\xi_i)$ of discontinuity of $c$, and $\RS_{i+1/2}^{n}$ at every discontinuity of $\bU_o^{n}$ in the interval $\mathring\Xi_i$.
Then, whenever two wave-fronts interact inside one of the strips $(0,+\infty) \times \mathring\Xi_i$, we prolong the approximate solution $\bU^{n}$ by applying $\RS_{i+1/2}^{n}$.
Furthermore, each time a wave-front reaches $x=\xi_{i}$ we extend $\bU^{n}$ by applying $\RS_i^{n}$.
In what follows, we refer to the times at which the approximate Riemann solvers are applied as \emph{interaction times}.

For $t>0$ sufficiently small, $\bU^{n}$ is then piecewise constant with jumps along a finite number of polygonal lines.
Conventionally, we assume that $\bU^{n}$ is right-continuous in space and left-continuous in time.
Then we can write
\begin{equation*}
\bU^{n}(t,x) = \sum_{j=0}^{I_{x}^{n}(t)} \bU_{j+1/2}^{n} \, \mathbbm{1}_{\calI_{j}^{n}(t)}(x),
\end{equation*}
where 
\begin{align*}
\calI_{0}^{n}(t) &\doteq (-\infty,x_{1}^{n}(t)),&
\calI_{j}^{n}(t) &\doteq [x_j^{n}(t),x_{j+1}^{n}(t)),\ j\in\{1,\ldots,I_{x}^{n}(t)-1\},&
\calI_{I_{x}^{n}(t)}^{n}(t) &\doteq [x_{I_{x}^{n}(t)}^{n}(t),+\infty),
\end{align*}
with $x_j^{n}(t) < x_{j+1}^{n}(t)$, $\bU_{j-1/2}^{n} \neq \bU_{j+1/2}^{n}$, and $t \mapsto I_{x}^{n}(t)$ constant between two interaction times.
Observe that $x_j^{n}(0^+)$ coincides with the position of a discontinuity of either $\bU_o^{n}$ or $c$.
Moreover, by construction we have
\[
x\in\Xi_{i} \Longrightarrow \bU^{n}(t,x) \in \calU_{i+1/2}^{n}.
\]

For clarity, we emphasise that $\bU^{n}(t,\cdot\,)$ may be continuous at $x=\xi_i$. Indeed, by \eqref{e:UI=UJ}\textsubscript{$2$}, $(\bU_{L},\bU_{R})\in\calU_{i-1/2}^{n}\times\calU_{i+1/2}^{n}$ is such that $\RS_i^{n}[\bU_L,\bU_R]$ is continuous at $\nu=0$ if and only if the following conditions are satisfied:
\begin{align}
\label{e:Vulfpeck}
&v(\bU_R) = 0,&
&w(\bU_L) = w(\bU_R).
\end{align}

In the following proposition, we provide a general expression for the propagation speed that holds for any wave-front.

\begin{proposition}
\label{p:StrappingYoungLad}
Away from the interaction times we have
\begin{equation*}
\dot{x}_j^{n}(t) = \sigma^{n}\left(\bU_{j-1/2}^{n},\bU_{j+1/2}^{n}, x_j^{n}(t)\right),
\end{equation*}
with
\begin{equation}
\label{e:Vola}
\sigma^{n}(\bU_-,\bU_+, x) \doteq
\frac{c(x^{+}) \, F_{w(\bU_{+})}(\rho_{+}) - c(x^{-}) \, F_{w(\bU_{-})}(\rho_{-})}{\rho_{+}-\rho_{-}}.
\end{equation}
Moreover, any discontinuity of $\bU^{n}$ satisfies the Rankine-Hugoniot condition \eqref{e:RH}.
\end{proposition}
\begin{proof}
By Proposition~\ref{p:InFlames2a}, any discontinuity that does not correspond to a discretised $1$-rarefaction satisfies the Rankine-Hugoniot condition \eqref{e:RH}, which implies \eqref{e:Vola}.
On the other hand, if the discontinuity corresponds to a discretised $1$-rarefaction, its propagation speed is given by \eqref{e:Haken} and is strictly negative.
This implies that there exists an index $i$ such that $x_j^{n}(t) \in \mathring{\Xi}_i$, and therefore $c(x_j^{n}(t)^{\pm}) = c_{i+1/2}$.
The conclusion then follows immediately, see \eqref{e:Haken}.
This concludes the proof.
\end{proof}

In the next theorem we list the main properties of the approximate solution $\bU^{n}$.
The goal of the following theorem is to uniformly bound $\tTV(\bU^{n}(t),c)$, the number of $c$-waves in $\bU^n(t)$, denoted by $I_{c}^{n}(t)$, the number of $1$-waves, $I_{1}^{n}(t)$, and $2$-contact discontinuities, $I_{2}^{n}(t)$, in $\bU^n(t)$ for any $t>0$ to ensure global existence of the approximate solution $\bU^n$.

\begin{theorem}
\label{t:Deftones}
Fix $n \in \N$.
For all $t>0$ we have that:
\begin{enumerate}[label={\bf(\arabic*)}]
\item\label{i:Deftones}
$\tTV(\bU^{n}(t),c) \leq \tTV(\bU_o^{n},c)$;
\item
the number of $c$-waves $I_{c}^{n}(t)$ is bounded from above by $I_{c}$;
\item
the number of $2$-waves $I_{2}^{n}(t)$ coincides with $I_{2}^{n}(0^+)$;
\item
the number of $1$-waves $I_{1}^{n}(t)$ is bounded from above by $I_{1}^{n}(0^+)$.
\end{enumerate}
At last, $\bU^{n}$ is globally well defined.
\end{theorem}
\begin{proof}
By slightly changing the speed of some waves, we can assume that only two waves interact at any single time.
Denote by $t>0$ a time at which an interaction occurs and by $\bU_L,\bU_M,\bU_R \in \calU$ the states involved in the interaction.
We want to estimate
\[\Delta\tTV \doteq \tTV(\bU^{n}(t^+),c) - \tTV(\bU^{n}(t^-),c).\]
We distinguish the following interaction types:
\begin{enumerate}[label={\bf(\roman*)}]
\item 
\label{i:A}
a $2$-contact discontinuity reaches $x=\xi_{i+1}$;
\item 
\label{i:B}
a $1$-wave interacts with a stationary discontinuity at $x=\xi_{i+1}$;\item 
\label{i:C}
a $1$-wave reaches $x=\xi_{i+1}$ and $\bU^n(t)$ is continuous at $x=\xi_{i+1}$;
\item 
\label{i:D}
a $2$-contact discontinuity interacts with a $1$-wave in $\mathring\Xi_{i}$;
\item 
\label{i:E}
two $1$-waves interact in $\mathring\Xi_{i}$.
\end{enumerate}
Below we make use of the properties established in claim~\ref{i:Puscifer} of Remark~\ref{r:traces}.

\begin{itemize}[leftmargin=*]
\begin{figure}[!htb]\centering
\begin{tikzpicture}[every node/.style={anchor=south west,inner sep=3pt},x=5mm, y=5mm, black]

\draw[-{Latex[scale=1.6]}] (0,0) -- (10,0) node[below] {$\ts$};
\draw[-{Latex[scale=1.6]}] (0,-4) -- (0,1) node[left] {\strut$w$};
\coordinate (A) at (5,-1);
\coordinate (B) at (5,-3);
\draw[fill=black] (A) circle (2pt) node[above] {\strut$(\bU_L,c_{i+1/2}), (\bU_{\star},c_{i+3/2})$};
\draw[fill=black] (B) circle (2pt) node[below] {\strut$(\bU_M,c_{i+1/2}), (\bU_R,c_{i+3/2})$};
\draw[dashed, thick, midarrow] (A) -- (B);

\begin{scope}[shift={(15,0)}]
\coordinate (C2m) at (0,-3);
\coordinate (C2p) at (10,-.5);
\coordinate (A) at (5,-4);
\coordinate (B) at (5,-1.5);
\coordinate (C) at (5,1);
\coordinate (UL) at (0,0);
\coordinate (UM) at (4,-3.5);
\coordinate (UR) at (10,-3);
\coordinate (US) at (10,.75);
\draw[fill=black] (B) circle (2pt);
\node[right, inner sep=0pt] at (UL) {\strut$(\bU_L,c_{i+1/2})$};
\node[left, inner sep=0pt] at (UM) {\strut$(\bU_M,c_{i+1/2})$};
\node[left] at (UR) {\strut$(\bU_R,c_{i+3/2})$};
\node[left] at (US) {\strut$(\bU_{\star},c_{i+3/2})$};
\draw[thick, midarrow] (A) node[below] {\strut$\xi_{i+1}$} -- node[sloped, below] {\strut$c$-wave} (B);
\draw[thick, midarrow] (B) -- node[sloped, above] {\strut$c$-wave} (C);
\draw[thick, midarrow] (C2m) -- node[sloped, above] {\strut$2$-con.\,disco.} (B);
\draw[thick, midarrow] (B) -- node[sloped, above] {\strut$2$-con.\,disco.} (C2p);
\end{scope}

\end{tikzpicture}
\caption{The interaction \ref{i:A} considered in the proof of Theorem~\ref{t:Deftones}, represented in the $\tbW$-coordinates introduced in~\eqref{e:tbW} (left) and in the $(x,t)$-plane (right).}
\label{f:A}
\end{figure}

\item 
Consider the interaction \ref{i:A}, see \figurename~\ref{f:A}.
In this case, before the interaction, the $2$-contact discontinuity $(\bU_{L},\bU_{M})$ moves in $\mathring\Xi_{i}$ with strictly positive propagation speed, hence $v(\bU_L) = v(\bU_M) > 0$.
This implies that $\bU_{M} \neq \bU_{R}$, see \eqref{e:Vulfpeck}, and that $(\bU_{M},\bU_{R})$ cannot be a stationary $2$-contact discontinuity.
Hence, $(\bU_{M},\bU_{R})$ is a $c$-wave .
By claim~\ref{i:Puscifer} of Remark~\ref{r:traces} we have
\begin{align*}
&\ts(\bU_L,c_{i+1/2}) = \ts(\bU_M,c_{i+1/2}) = \ts(\bU_R,c_{i+3/2}),&
&w(\bU_M) = w(\bU_R).
\end{align*}
As a consequence $\RS_{i+1}^{n}[\bU_L,\bU_R]$ involves a stationary wave $(\bU_{L},\bU_{\star})$ and a $2$-contact discontinuity $(\bU_{\star},\bU_{R})$, hence the number of waves does not change after the interaction.
Moreover we have
\begin{align*}
&\ts(\bU_L,c_{i+1/2}) = \ts(\bU_{\star},c_{i+3/2}) = \ts(\bU_R,c_{i+3/2}),&
&w(\bU_L) = w(\bU_{\star}),
\end{align*}
and therefore
\[\Delta\tTV = |w(\bU_\star) - w(\bU_R)| - |w(\bU_L) - w(\bU_M)| = 0.\]

\begin{figure}[!htb]\centering
\begin{tikzpicture}[every node/.style={anchor=south west,inner sep=3pt},x=5mm, y=5mm, black]

\draw[-{Latex[scale=1.6]}] (0,0) -- (10,0) node[below] {$\ts$};
\draw[-{Latex[scale=1.6]}] (0,-5) -- (0,1) node[left] {\strut$w$};
\coordinate (A) at (0,-1.5);
\coordinate (B) at (0,-4);
\coordinate (C) at (2,-1.5);
\coordinate (Uh) at (2,-1);
\coordinate (Uc) at (2,-2);
\coordinate (D) at (2,-4);
\draw[fill=black] (A) circle (2pt) node[left] {\strut$(\bU_L,c_{i+1/2})$};
\draw[fill=black] (C) circle (2pt);
\node[right] at (Uh) {\strut$(\hat{\bU},c_{i+1/2})$};
\node[right] at (Uc) {\strut$(\check{\bU},c_{i+3/2})$};
\draw[fill=black] (B) circle (2pt) node[left] {\strut$(\bU_M,c_{i+1/2})$};
\draw[fill=black] (D) circle (2pt) node[right] {\strut$(\bU_R,c_{i+3/2})$};
\draw[dashed, thick, midarrow] (A) -- (B);
\draw[dashed, thick, midarrow] (A) -- (C);
\draw[dashed, thick, midarrow] (B) -- (D);
\draw[dashed, thick, midarrow] (C) -- (D);
\draw[dotted, thick] (C) -- (2,0) node[above] {\strut$\varepsilon_{1}^{n}$};

\begin{scope}[shift={(15,0)}]
\coordinate (C2m) at (9.25,-4.5);
\coordinate (C2p) at (0,-.5);
\coordinate (A) at (5,-5);
\coordinate (B) at (5,-1.5);
\coordinate (C) at (5,1);
\coordinate (CC) at (10,.25);
\coordinate (UL) at (.5,-3.5);
\coordinate (UM) at (8.5,-4.5);
\coordinate (UH) at (1,.5);
\coordinate (UC) at (5.25,.8);
\coordinate (UR) at (10,-1.5);
\draw[fill=black] (B) circle (2pt);
\node[right, inner sep=0pt] at (UL) {\strut$(\bU_L,c_{i+1/2})$};
\node[left, inner sep=0pt] at (UM) {\strut$(\bU_M,c_{i+3/2})$};
\node[left, inner sep=0pt] at (UR) {\strut$(\bU_R,c_{i+3/2})$};
\node[right, inner sep=0pt] at (UH) {\strut$(\hat{\bU},c_{i+1/2})$};
\node[right, inner sep=0pt] at (UC) {\strut$(\check{\bU},c_{i+3/2})$};
\draw[thick, midarrow] (A) node[below] {\strut$\xi_{i+1}$} -- node[sloped, above] {\strut$2$-con.\,disco} (B);
\draw[thick, midarrow] (B) -- node[sloped, above] {\strut$c$-wave} (C);
\draw[thick, midarrow] (C2m) -- node[sloped, above] {\strut$1$-rare.} (B);
\draw[thick, midarrow] (B) -- node[sloped, above] {\strut$1$-rare.} (C2p);
\draw[thick, midarrow] (B) -- node[sloped, above] {\strut$2$-con.\,disco.} (CC);
\end{scope}

\end{tikzpicture}
\caption{The interaction corresponding to item~\ref{i:B} considered in the proof of Theorem~\ref{t:Deftones}, represented in the $\tbW$-coordinates (left) and in the $(x,t)$-plane (right), in the case the stationary discontinuity $(\bU_L,\bU_M)$ is a $2$-contact discontinuity.}
\label{f:Ba}
\end{figure}

\item 
Consider the interaction \ref{i:B} and assume that the stationary discontinuity at $x=\xi_{i+1}$ is a $2$-contact discontinuity, see \figurename~\ref{f:Ba}.
In this case $\ts(\bU_L,c_{i+1/2}) = \ts(\bU_M,c_{i+3/2}) = 0$ and $w(\bU_L) \neq w(\bU_M)$.
This implies that the $1$-wave $(\bU_M,\bU_R)$ is a $1$-rarefaction with $\ts(\bU_R,c_{i+3/2}) = \varepsilon_{1}^{n}$ and $w(\bU_M) = w(\bU_R)$.
As a consequence $\RS_{i+1}^{n}[\bU_L,\bU_R]$ involves a $1$-rarefaction $(\bU_{L},\hat{\bU})$, a stationary $c$-wave $(\hat{\bU},\check{\bU})$ and a $2$-contact discontinuity $(\check{\bU},\bU_{R})$, where $\hat{\bU}$ and $\check{\bU}$ are given in \eqref{e:Avatar}.
Therefore, the number of waves increases by one after the interaction and
\[
\Delta\tTV =
\begin{aligned}[t]
&|\ts(\bU_L,c_{i+1/2}) - \ts(\hat{\bU},c_{i+1/2})| + |w(\check{\bU}) - w(\bU_R)| 
\\-{}&
|w(\bU_L) - w(\bU_M)| - |\ts(\bU_M,c_{i+3/2}) - \ts(\bU_R,c_{i+3/2})| = 0.
\end{aligned}
\]

\begin{figure}[!htb]\centering
\begin{tikzpicture}[every node/.style={anchor=south west,inner sep=3pt},x=6mm, y=5mm, black]

\draw[-{Latex[scale=1.6]}] (0,0) -- (8,0) node[below] {$\ts$};
\draw[-{Latex[scale=1.6]}] (0,-4) -- (0,1) node[left] {\strut$w$};
\coordinate (A) at (6,-2);
\coordinate (B) at (2,-2);
\draw[fill=black] (A) circle (2pt) node[above] {\strut$(\bU_L,c_{i+1/2})$} node[below] {\strut$(\bU_M,c_{i+3/2})$};
\draw[fill=black] (B) circle (2pt) node[above] {\strut$(\bU_{\star},c_{i+1/2})$} node[below] {\strut$(\bU_R,c_{i+3/2})$};
\draw[dashed, thick, midarrow] (A) -- (B);

\begin{scope}[shift={(10,0)}]
\draw[-{Latex[scale=1.6]}] (0,0) -- (8,0) node[below] {$\ts$};
\draw[-{Latex[scale=1.6]}] (0,-4) -- (0,1) node[left] {\strut$w$};
\coordinate (A) at (3.5,-2);
\coordinate (B) at (5,-2);
\draw[dotted, thick] (A) -- (3.5,0) node[above] {\strut$j\,\varepsilon_{1}^{n}$};
\draw[dotted, thick] (B) -- (5,0) node[above] {\strut$\quad (j+1)\,\varepsilon_{1}^{n}$};

\draw[fill=black] (A) circle (2pt) node[left] {\strut$\displaystyle\genfrac{}{}{0pt}{}{(\bU_L,c_{i+1/2})}{\displaystyle (\bU_M,c_{i+3/2})}$};
\draw[fill=black] (B) circle (2pt) node[right] {\strut$\displaystyle \genfrac{}{}{0pt}{}{(\bU_{\star},c_{i+1/2})}{(\bU_R,c_{i+3/2})}$};
\draw[dashed, thick, midarrow] (A) -- (B);
\end{scope}

\begin{scope}[shift={(20,0)}]
\coordinate (C2m) at (10,-3);
\coordinate (C2p) at (0,-.5);
\coordinate (A) at (5,-4);
\coordinate (B) at (5,-1.5);
\coordinate (C) at (5,1);
\coordinate (UL) at (0,-2.5);
\coordinate (UM) at (10,-3.5);
\coordinate (UR) at (10,0);
\coordinate (US) at (0,.5);
\draw[fill=black] (B) circle (2pt);
\node[right, inner sep=0pt] at (UL) {\strut$(\bU_L,c_{i+1/2})$};
\node[left, inner sep=0pt] at (UM) {\strut$(\bU_M,c_{i+3/2})$};
\node[left, inner sep=0pt] at (UR) {\strut$(\bU_R,c_{i+3/2})$};
\node[right, inner sep=0pt] at (US) {\strut$(\bU_{\star},c_{i+1/2})$};
\draw[thick, midarrow] (A) node[below] {\strut$\xi_{i+1}$} -- node[sloped, above] {\strut$c$-wave} (B);
\draw[thick, midarrow] (B) -- node[sloped, below] {\strut$c$-wave} (C);
\draw[thick, midarrow] (C2m) -- node[sloped, above] {\strut$1$-wave} (B);
\draw[thick, midarrow] (B) -- node[sloped, above] {\strut$1$-wave} (C2p);
\end{scope}

\end{tikzpicture}
\caption{The interaction corresponding to item~\ref{i:B} considered in the proof of Theorem~\ref{t:Deftones}, represented in the $\tbW$-coordinates (left and center) and in the $(x,t)$-plane (right), in the case the stationary discontinuity $(\bU_L,\bU_M)$ is a $c$-wave with $\ts(\bU_R,c_{i+3/2}) > 0$.}
\label{f:Bb}
\end{figure}

\item 
Consider the interaction \ref{i:B} and assume that the stationary discontinuity at $x=\xi_{i+1}$ is a $c$-wave with $\ts(\bU_R,c_{i+3/2}) > 0$, see \figurename~\ref{f:Bb}.
In this case $\ts(\bU_L,c_{i+1/2}) = \ts(\bU_M,c_{i+3/2}) > 0$ and $w(\bU_L) = w(\bU_M)$.
Then, the $1$-wave $(\bU_M,\bU_R)$ is either a $1$-shock with $\ts(\bU_R,c_{i+3/2}) < \ts(\bU_M,c_{i+3/2})$ and $w(\bU_M) = w(\bU_R)$, or a $1$-rarefaction with $\ts(\bU_R,c_{i+3/2}) = \ts(\bU_M,c_{i+3/2}) + \varepsilon_{1}^{n}$ and $w(\bU_M) = w(\bU_R)$.
In the former case $\RS_{i+1}^{n}[\bU_L,\bU_R]$ involves a $1$-shock $(\bU_{L},\bU_{\star})$ and a $c$-wave $(\bU_{\star},\bU_{R})$, whereas in the latter case a $1$-rarefaction $(\bU_{L},\bU_{\star})$ and a $c$-wave $(\bU_{\star},\bU_{R})$.
In both the cases, $\bU_{\star}$ satisfies
\begin{align*}
&\ts(\bU_{\star},c_{i+1/2}) = \ts(\bU_R,c_{i+3/2}),&
&w(\bU_{\star}) = w(\bU_R).
\end{align*}
Therefore, the number of waves does not change after the interaction and
\[
\Delta\tTV =
|\ts(\bU_L,c_{i+1/2}) - \ts(\bU_{\star},c_{i+1/2})| - |\ts(\bU_M,c_{i+3/2}) - \ts(\bU_R,c_{i+3/2})| = 0.
\]

\begin{figure}[!htb]\centering
\begin{tikzpicture}[every node/.style={anchor=south west,inner sep=3pt},x=5mm, y=5mm, black]

\draw[-{Latex[scale=1.6]}] (0,0) -- (10,0) node[below] {$\ts$};
\draw[-{Latex[scale=1.6]}] (0,-4) -- (0,1) node[left] {\strut$w$};
\coordinate (A) at (6,-2);
\coordinate (B) at (0,-2);
\draw[fill=black] (A) circle (2pt) node[above] {\strut$(\bU_L,c_{i+1/2})$} node[below] {\strut$(\bU_M,c_{i+3/2})$};
\draw[fill=black] (B) circle (2pt) node[left] {\strut$(\bU_R,c_{i+3/2})$};
\draw[dashed, thick, midarrow] (A) -- (B);

\begin{scope}[shift={(15,0)}]
\coordinate (C2m) at (10,-3);
\coordinate (C2p) at (0,-.5);
\coordinate (A) at (5,-4);
\coordinate (B) at (5,-1.5);
\coordinate (C) at (5,1);
\coordinate (UL) at (0,-2.5);
\coordinate (UM) at (10,-3.5);
\coordinate (UR) at (10,0);
\coordinate (US) at (0,.5);
\draw[fill=black] (B) circle (2pt);
\node[right, inner sep=0pt] at (UL) {\strut$(\bU_L,c_{i+1/2})$};
\node[left, inner sep=0pt] at (UM) {\strut$(\bU_M,c_{i+3/2})$};
\node[left, inner sep=0pt] at (UR) {\strut$(\bU_R,c_{i+3/2})$};
\node[right, inner sep=0pt] at (US) {\strut$(\bU_{\star},c_{i+1/2})$};
\draw[thick, midarrow] (A) node[below] {\strut$\xi_{i+1}$} -- node[sloped, above] {\strut$c$-wave} (B);
\draw[thick, dashed] (B) -- (C);
\draw[thick, midarrow] (C2m) -- node[sloped, above] {\strut$1$-shock} (B);
\draw[thick, midarrow] (B) -- node[sloped, above] {\strut$1$-shock} (C2p);
\end{scope}

\end{tikzpicture}
\caption{The interaction corresponding to item~\ref{i:B} considered in the proof of Theorem~\ref{t:Deftones}, represented in the $\tbW$-coordinates (left and center) and in the $(x,t)$-plane (right), in the case the stationary discontinuity $(\bU_L,\bU_M)$ is a $c$-wave with $\ts(\bU_R,c_{i+3/2}) = 0$.}
\label{f:Bbb}
\end{figure}

\item 
Consider the interaction \ref{i:B} and assume that the stationary discontinuity at $x=\xi_{i+1}$ is a $c$-wave with $\ts(\bU_R,c_{i+3/2}) = 0$, see \figurename~\ref{f:Bbb}.
In this case $\ts(\bU_L,c_{i+1/2}) = \ts(\bU_M,c_{i+3/2}) > 0$ and $w(\bU_L) = w(\bU_M)$.
Then, the $1$-wave $(\bU_M,\bU_R)$ is a $1$-shock with $0 = \ts(\bU_R,c_{i+3/2}) < \ts(\bU_M,c_{i+3/2})$ and $w(\bU_M) = w(\bU_R)$.
Hence $\RS_{i+1}^{n}[\bU_L,\bU_R]$ involves a $1$-shock $(\bU_{L},\bU_{R})$ and is continuous at $\nu=0$.
Therefore, the number of waves decreases by one after the interaction and
\[
\Delta\tTV =
|\ts(\bU_L,c_{i+1/2}) - \ts(\bU_{R},c_{i+1/2})| - |\ts(\bU_M,c_{i+3/2}) - \ts(\bU_R,c_{i+3/2})| = 0.
\]

\begin{figure}[!htb]\centering
\begin{tikzpicture}[every node/.style={anchor=south west,inner sep=3pt},x=5mm, y=5mm, black]

\draw[-{Latex[scale=1.6]}] (0,0) -- (10,0) node[below] {$\ts$};
\draw[-{Latex[scale=1.6]}] (0,-4) -- (0,1) node[left] {\strut$w$};
\coordinate (A) at (0,-2);
\coordinate (B) at (2,-2);
\draw[dotted, thick] (B) -- (2,0) node[above] {\strut$\varepsilon_{1}^{n}$};

\draw[fill=black] (A) circle (2pt) node[left] {\strut$\displaystyle\genfrac{}{}{0pt}{}{(\bU_L,c_{i+1/2})}{\displaystyle (\bU_M,c_{i+3/2})}$};
\draw[fill=black] (B) circle (2pt) node[right] {\strut$\displaystyle \genfrac{}{}{0pt}{}{(\bU_{\star},c_{i+1/2})}{(\bU_R,c_{i+3/2})}$};
\draw[dashed, thick, midarrow] (A) -- (B);

\begin{scope}[shift={(13,0)}]
\coordinate (C2m) at (10,-3);
\coordinate (C2p) at (0,-.5);
\coordinate (A) at (5,-4);
\coordinate (B) at (5,-1.5);
\coordinate (C) at (5,1);
\coordinate (UL) at (0,-2.5);
\coordinate (UM) at (10,-3.5);
\coordinate (UR) at (10,0);
\coordinate (US) at (0,.5);
\draw[fill=black] (B) circle (2pt);
\node[right, inner sep=0pt] at (UL) {\strut$(\bU_L,c_{i+1/2})$};
\node[left, inner sep=0pt] at (UM) {\strut$(\bU_M,c_{i+3/2})$};
\node[left, inner sep=0pt] at (UR) {\strut$(\bU_R,c_{i+3/2})$};
\node[right, inner sep=0pt] at (US) {\strut$(\bU_{\star},c_{i+1/2})$};
\draw[thick, dashed] (A) node[below] {\strut$\xi_{i+1}$} -- (B);
\draw[thick, midarrow] (B) -- node[sloped, above] {\strut$c$-wave} (C);
\draw[thick, midarrow] (C2m) -- node[sloped, above] {\strut$1$-rare.} (B);
\draw[thick, midarrow] (B) -- node[sloped, above] {\strut$1$-rare.} (C2p);
\end{scope}

\end{tikzpicture}
\caption{The interaction corresponding to item~\ref{i:C} considered in the proof of Theorem~\ref{t:Deftones}, represented in the $\tbW$-coordinates (left) and in the $(x,t)$-plane (right).}
\label{f:C}
\end{figure}

\item 
Consider the interaction \ref{i:C}, see \figurename~\ref{f:C}.
In this case $\ts(\bU_L,c_{i+1/2}) = \ts(\bU_M,c_{i+3/2}) = 0$ and $w(\bU_L) = w(\bU_M)$.
Then, the $1$-wave $(\bU_M,\bU_R)$ is a $1$-rarefaction with $\ts(\bU_R,c_{i+3/2}) = \varepsilon_{1}^{n}$ and $w(\bU_M) = w(\bU_R)$.
Hence $\RS_{i+1}^{n}[\bU_L,\bU_R]$ involves a $1$-rarefaction $(\bU_{L},\bU_{\star})$ and a $c$-wave $(\bU_{\star},\bU_{R})$, where $\bU_{\star}$ satisfies
\begin{align*}
&\ts(\bU_{\star},c_{i+1/2}) = \ts(\bU_R,c_{i+3/2}),&
&w(\bU_{\star}) = w(\bU_L) = w(\bU_R).
\end{align*}
Therefore, the number of waves increases by one after the interaction and
\[
\Delta\tTV =
|\ts(\bU_L,c_{i+1/2}) - \ts(\bU_{\star},c_{i+1/2})| - |\ts(\bU_M,c_{i+3/2}) - \ts(\bU_R,c_{i+3/2})| = 0.
\]

\item 
Consider the interaction \ref{i:D}.
In this case we have
\begin{align*}
&\ts(\bU_L,c_{i+1/2}) = \ts(\bU_M,c_{i+1/2}) \neq \ts(\bU_R,c_{i+1/2}),&
&w(\bU_L) \neq w(\bU_M) = w(\bU_R).
\end{align*}
Hence $\RS_{i+1/2}^{n}[\bU_L,\bU_R]$ involves a $1$-wave $(\bU_{L},\bU_{\star})$ and a $2$-contact discontinuity $(\bU_{\star},\bU_{R})$, where $\bU_{\star}$ satisfies
\begin{align*}
&\ts(\bU_{\star},c_{i+1/2}) = \ts(\bU_R,c_{i+1/2}),&
&w(\bU_{\star}) = w(\bU_L).
\end{align*}
Therefore, the number of waves does not change after the interaction and
\[
\Delta\tTV =
\begin{aligned}[t]
&|\ts(\bU_L,c_{i+1/2}) - \ts(\bU_{\star},c_{i+1/2})| + |w(\bU_\star) - w(\bU_R)| 
\\-{}&
|w(\bU_L) - w(\bU_M)| - |\ts(\bU_{M},c_{i+1/2}) - \ts(\bU_R,c_{i+1/2})| = 0.
\end{aligned}
\]

\item 
Consider the interaction \ref{i:E}.
We distinguish the following cases.
\begin{itemize}
\item 
If two $1$-shocks interact, then
\begin{align*}
&\ts(\bU_L,c_{i+1/2}) > \ts(\bU_M,c_{i+1/2}) > \ts(\bU_R,c_{i+1/2}),&
&w(\bU_L) = w(\bU_M) = w(\bU_R).
\end{align*}
As a consequence $\RS_{i+1/2}^{n}[\bU_L,\bU_R]$ has just a $1$-shock.
Therefore, the number of waves decreases by one after the interaction and
\[\Delta\tTV = 
\begin{aligned}[t]
&|\ts(\bU_L,c_{i+1/2}) - \ts(\bU_R,c_{i+1/2})| 
\\-{}&
|\ts(\bU_L,c_{i+1/2}) - \ts(\bU_M,c_{i+1/2})| - |\ts(\bU_M,c_{i+1/2}) - \ts(\bU_R,c_{i+1/2})| = 0.
\end{aligned}
\]
\item 
If $(\bU_L,\bU_M)$ is a $1$-shock and $(\bU_M,\bU_R)$ is a $1$-rarefaction, then
\begin{align*}
&\ts(\bU_L,c_{i+1/2}) > \ts(\bU_R,c_{i+1/2}) > \ts(\bU_M,c_{i+1/2}),&
&w(\bU_L) = w(\bU_M) = w(\bU_R).
\end{align*}
As a consequence $\RS_{i+1/2}^{n}[\bU_L,\bU_R]$ has just a $1$-shock.
Therefore, the number of waves decreases by one after the interaction and
\begin{align*}
\Delta\tTV &= 
\begin{aligned}[t]
&|\ts(\bU_L,c_{i+1/2}) - \ts(\bU_R,c_{i+1/2})|
\\-{}&
|\ts(\bU_L,c_{i+1/2}) - \ts(\bU_M,c_{i+1/2})| - |\ts(\bU_M,c_{i+1/2}) - \ts(\bU_R,c_{i+1/2})|
\end{aligned}
\\&= -2 \left( \ts(\bU_R,c_{i+1/2})-\ts(\bU_M,c_{i+1/2}) \right) = -2 \varepsilon_{1}^{n} < 0.
\end{align*}
\item 
If $(\bU_L,\bU_M)$ is a $1$-rarefaction and $(\bU_M,\bU_R)$ is a $1$-shock, then
\begin{align*}
&\ts(\bU_M,c_{i+1/2}) > \ts(\bU_L,c_{i+1/2}) > \ts(\bU_R,c_{i+1/2}),&
&w(\bU_L) = w(\bU_M) = w(\bU_R).
\end{align*}
As a consequence $\RS_{i+1/2}^{n}[\bU_L,\bU_R]$ has just a $1$-shock.
Therefore, the number of waves decreases by one after the interaction and
\begin{align*}
\Delta\tTV &= 
\begin{aligned}[t]
&|\ts(\bU_L,c_{i+1/2}) - \ts(\bU_R,c_{i+1/2})|
\\-{}&
|\ts(\bU_L,c_{i+1/2}) - \ts(\bU_M,c_{i+1/2})| - |\ts(\bU_M,c_{i+1/2}) - \ts(\bU_R,c_{i+1/2})|
\end{aligned}
\\&= -2 \left( \ts(\bU_M,c_{i+1/2})-\ts(\bU_L,c_{i+1/2}) \right) = -2 \varepsilon_{1}^{n} < 0.
\end{align*}
\end{itemize}
\end{itemize}
Observe that the above case by case study shows that interaction times cannot accumulate in finite time.
Indeed, any $1$-wave, $2$-contact discontinuity, or $c$-wave can cross any other wave at most once.
Furthermore, we proved that $\Delta\tTV  \leq 0$.
By \eqref{e:uoncond2}\textsubscript{$1$} we then have $\tTV(\bU^{n}(t),c) \leq \tTV(\bU_o^{n},c) \leq \tTV(\bU_o,c)$.
As a consequence, the number of waves is bounded by $\tTV(\bU_o,c) \left(\frac{1}{\varepsilon_{1}^{n}}+\frac{1}{\varepsilon_{2}^{n}}\right) + I_{c}$.
This concludes the proof.
\end{proof}

\begin{remark}
\label{r:Malmsteen}
It is worth noting that the total variation of $\bU^{n}(t)$ may well increase after an interaction.
Indeed, if for instance at time $t>0$ a $2$-contact discontinuity $(\bU_L,\bU_M)$ interacts with a $1$-shock $(\bU_M,\bU_R)$ at $x \in \mathring\Xi_{i}$ and
\begin{align*}
&v(\bU_L) = v(\bU_M) > v(\bU_R),&
&w(\bU_L) > w(\bU_M) = w(\bU_R),
\end{align*}
then, after the interaction, $\bU^n$ locally has at $x$ a $1$-shock $(\bU_L,\bU_\star)$ and a $2$-contact discontinuity $(\bU_\star,\bU_R)$, with $\bU_\star$ such that
\begin{align*}
&v(\bU_L) > v(\bU_\star) = v(\bU_R),&
&w(\bU_L) = w(\bU_\star) > w(\bU_R).
\end{align*}
As a consequence, by \eqref{e:TheMarsVolta2}\textsubscript{$1$} and \eqref{e:pm1} we have
\begin{align*}
\TV\left(\bU(t^+)\right) - \TV\left(\bU(t^-)\right)
&=
(2\rho_\star-\rho_L-\rho_R) - (\rho_L+\rho_R-2\rho_M)
\\&=
\frac{2}{\Rref} \left( p^{-1}\left(w(\bU_L)\right) - p^{-1}\left(w(\bU_R)\right) \right) \left( p^{-1}\left(-v(\bU_R)\right) - p^{-1}\left(-v(\bU_L)\right) \right) > 0.
\end{align*}
We observe that the total variation of $\bW(\bU)$ does not increase after such an interaction.
In fact, this property holds for any interaction occurring in $\mathring{\Xi}_i$; however, it may fail for interactions taking place at the discontinuities of $c$.
Indeed, if for instance at time $t>0$ a $1$-rarefaction $(\bU_L,\bU_R)$ reaches $x=\xi_i$ and $\bU^n(t^-,\xi_i^-) = \bU^n(t^-,\xi_i^+) = \bU_L$, then as already observed in the proof of Theorem~\ref{t:Deftones}, see \figurename~\ref{f:C}, we have
\begin{align*}
&v(\bU_L) = 0 < v(\bU_R),&
&w(\bU_L) = w(\bU_R).
\end{align*}
After such interaction, $\bU^n$ locally has at $x=\xi_i$ a $1$-rarefaction $(\bU_L,\bU_\star)$ and a $c$-wave $(\bU_\star,\bU_R)$, with $\bU_\star$ such that
\begin{align*}
&\ts(\bU_{\star},c_{i-1/2}) = \ts(\bU_R,c_{i+1/2}),&
&w(\bU_{\star}) = w(\bU_L) = w(\bU_R).
\end{align*}
If $c_{i-1/2} < c_{i+1/2}$, then by \eqref{e:SnarkyPuppy} it follows
\[
c_{i-1/2} \, S(v(\bU_\star)) = c_{i+1/2} \, S(v(\bU_R))
\Longrightarrow
v(\bU_\star) = S^{-1}\left( \frac{c_{i+1/2}}{c_{i-1/2}} \, S(v(\bU_R)) \right) > v(\bU_R),\]
because $S$ defined in \eqref{e:S} is strictly increasing.
As a consequence, we have
\begin{align*}
\TV\left(\bW\left(\bU(t^+)\right)\right) - \TV\left(\bW\left(\bU(t^-)\right)\right)
&=
\left( 2v(\bU_\star) - v(\bU_L) - v(\bU_R) \right) - \left( v(\bU_R) - v(\bU_L) \right)
\\&=
2\left( v(\bU_\star) - v(\bU_R) \right) > 0.
\end{align*}
\end{remark}

\begin{corollary}
\label{c:Unweak}
$\bU^{n}$ is a weak solution of the Cauchy problem for the system \eqref{e:2x2} with initial condition $\bU(0,x) = \bU_o^{n}(x)$ in the sense of Definition~\ref{d:ws2}.
\end{corollary}
\begin{proof}
It is sufficient to observe that any discontinuity of the approximate solution $\bU^{n}$ satisfies the Rankine-Hugoniot condition \eqref{e:RH}, see Proposition~\ref{p:StrappingYoungLad}, and that it is globally defined by Theorem~\ref{t:Deftones}.
This concludes the proof.
\end{proof}

\subsection{Convergence}

We now prove that $\{\tbW(\bU^{n},c)\}_n$ converges in $\Lloc1$ (along a subsequence) to $\tbW(\bU,c)$ where $\bU$ is, in fact, an entropy solution in the sense of Definition~\ref{d:entro2}.

For any $t>0$, by \eqref{e:tTV}, claim~\ref{i:Deftones} of Theorem~\ref{t:Deftones}, and \eqref{e:uoncond2}\textsubscript{$1$} we have that
\begin{equation}
\label{e:Karnivool}
\TV(\tbW(\bU^{n}(t),c)) = \tTV(\bU^{n}(t),c) \leq \tTV(\bU_o^{n},c) \leq \tTV(\bU_o,c).
\end{equation}
Furthermore, by construction, see \eqref{e:Leprous} and \figurename~\ref{f:Puscifer}, and by \eqref{e:uoncond1}\textsubscript{$1$} it follows that
\[
\|\tbW(\bU^{n}(t),c)\|_{\L\infty(\R)} = \|\tbW(\bU^{n}_o,c)\|_{\L\infty(\R)} \leq \|\tbW(\bU_o,c)\|_{\L\infty(\R)}.
\]

Observe that by \eqref{e:BlacklitCanopy}, the absolute value of the propagation speed of any wave is bounded by 
\[\|c\|_{\L\infty(\R)} \, \max\left\{ \vref \, V'(0) , V(\ov) \right\}.\]
As a consequence, by \eqref{e:uoncond3} and \eqref{e:Karnivool} for any $0 \leq s < t$ there holds
\begin{equation}
\label{e:BlindMelon}
\|\tbW(\bU^{n}(t),c)-\tbW(\bU^{n}(s),c)\|_{\L1(\R)} \leq \widetilde{L} \, (t-s),
\end{equation}
with
\begin{equation}
\label{e:tL}
\widetilde{L} \doteq \tTV(\bU_o,c) \, \|c\|_{\L\infty(\R)} \, \max\left\{ \vref \, V'(0) , V(\ov) \right\}.
\end{equation}

By \eqref{e:Karnivool}, Helly's theorem, and the Cantor diagonal process, one finds a subsequence, still denoted $\{\tbW(\bU^{n},c)\}_n$, such that $\{\tbW(\bU^{n}(t,\cdot\,),c)\}_n$ is convergent in $\Lloc1(\R;\calU)$, for any rational $t \in [0,+\infty)$.
Then, \eqref{e:BlindMelon} implies that $\{\tbW(\bU^{n}(t,\cdot\,),c)\}_n$ is Cauchy in $\Lloc1(\R;\calU)$ for all $t \in [0,+\infty)$, and hence $\{\tbW(\bU^{n},c)\}_n$ converges in $\Lloc1$ to some function $\tbW^{\infty} = (\ts^{\infty},w^{\infty})$ of locally bounded variation on $\D$.
Moreover, $\tbW^{\infty}$ satisfies the following estimates for any $0 \leq s < t$ and $a<b$
\begin{gather}
\nonumber
\left\|\tbW^{\infty}(s)\right\|_{\L\infty(\R)} \leq \left\|\tbW(\bU_o)\right\|_{\L\infty(\R)},
\qquad
\TV\left(\tbW^{\infty}(s)\right) \leq \TV\left(\tbW(\bU_{o},c)\right),
\\
\label{e:Manowar}
\int_a^b \left\|\tbW^{\infty}(t,x)-\tbW^{\infty}(s,x)\right\| \, \d x \leq \widetilde{L} \, (t-s).
\end{gather}

From \eqref{e:TOOL}, \eqref{e:Polyphia}, \eqref{e:NineInchNails} and \eqref{e:dominva}, it follows that
\[t\geq0,\ x \in \Xi_i \Longrightarrow \bU^n(t,x) \in \calU_{i+1/2}^{n} \subset \calU_{c_{i+1/2}} \Longrightarrow \ts\left(\bU^n(t,x),c_{i+1/2}(x)\right) \in [0,\oV].\]
As a consequence, also $\ts^{\infty}$ takes values in $[0,\oV]$ and therefore by \eqref{e:cV}, \eqref{e:S}, and \eqref{e:pm1} we have that
\[
0 \leq \frac{\ts^{\infty}}{c} \, p^{-1}(w^{\infty}) \leq
\frac{\oV}{\uc} \, p^{-1}(w^{\infty}) \leq 
V(\ov) \, \frac{\Rref}{p^{-1}(\ov)} \, p^{-1}(w^{\infty}) 
= p^{-1}(w^{\infty}-\ov)\, V(\ov).
\]
By~\eqref{e:F} and Lemma~\ref{l:F}, there exists a unique state
$\bU^{\infty}$,
determined by~\eqref{e:TheMarsVolta}, such that
$\tbW^{\infty} = \tbW(\bU^{\infty},c)$.
More explicitly, the components of $\bU^{\infty} = (\rho^{\infty}, q^{\infty})$ are given by
\begin{align*}
\rho^{\infty}
&= F_{w^{\infty}}^{-1}\left(
\frac{\ts^{\infty}}{c}\, p^{-1}(w^{\infty})
\right),
&
q^{\infty}
&= w^{\infty}\,
F_{w^{\infty}}^{-1}\left(
\frac{\ts^{\infty}}{c}\, p^{-1}(w^{\infty})
\right).
\end{align*}

\begin{lemma}
\label{l:Gojira}
$\{\bU^{n}\}_{n}$ converges to $\bU^{\infty}$ in $\Lloc1$.
\end{lemma}
\begin{proof}
For simplicity of notation, we omit the dependence on $(t,x)$ and let
\begin{align*}
&\ts^{n} = \ts(\bU^{n},c),&
&\ts^{\infty} = \ts(\bU^{\infty},c),&
&w^{n} = w(\bU^{n}),&
&w^{\infty} = w(\bU^{\infty}).
\end{align*}
By \eqref{e:SnarkyPuppy} and \eqref{e:TheMarsVolta2}\textsubscript{$1$} it follows that $\rho = p^{-1}\left(w-S^{-1}(\ts/c)\right)$.
Observe that
\begin{align*}
&\Lip(p^{-1}) = \frac{\Rref}{\vref},&
&\Lip(S^{-1}) = \frac{\vref \, p^{-1}(\ov)}{\left(\vref \, V'(\ov) - V(\ov)\right) \Rref},
\end{align*}
and therefore
\begin{align*}
|\rho^{n}-\rho^{\infty}| &=
\left| 
p^{-1}\left(w^{n} - S^{-1}\left(\ts^{n}/c\right)\right)
-
p^{-1}\left(w^{\infty} - S^{-1}\left(\ts^{\infty}/c\right)\right)
\right|
\\&\leq
\frac{\Rref}{\vref} \left|w^{n}-w^{\infty}\right| + \frac{p^{-1}(\ov)}{\uc \left(\vref \, V'(\ov) - V(\ov)\right)} \left|\ts^{n} - \ts^{\infty}\right|.
\end{align*}
Hence, have that $\{\rho^{n}\}_{n}$ converges to $\rho^{\infty}$ in $\Lloc1$.
At last, we have
\[
|q^{n}-q^{\infty}| =
|w^{n} \, \rho^{n} - w^{\infty} \, \rho^{\infty}| \leq
-w^{n} \, |\rho^{n} - \rho^{\infty}| + |w^{n} - w^{\infty}| \, \rho^{\infty}
\leq
-\uw \, |\rho^{n} - \rho^{\infty}| + |w^{n} - w^{\infty}| \, \Rref,
\]
and therefore $\{q^{n}\}_{n}$ converges to $q^{\infty}$ in $\Lloc1$.
This concludes the proof.
\end{proof}

\begin{lemma}
\label{l:LipUt}
There exists $L_{\bU}>0$ such that for any $0\leq s<t$ and $a<b$ we have
\begin{equation}
\label{e:LipUt}
\int_a^b \|\bU^{\infty}(t,x)-\bU^{\infty}(s,x)\| \, \d x \leq L_{\bU} \, (t-s).
\end{equation}
\end{lemma}
\begin{proof}
By applying an analogous procedure as in the proof of Lemma~\ref{l:Gojira}, we get
\begin{align*}
\int_a^b \left|\rho^{\infty}(t,x)-\rho^{\infty}(s,x)\right| \d x
={}&
\int_a^b \left|
p^{-1}\left(w^{\infty}(t,x) - S^{-1}\left(\frac{\ts^{\infty}(t,x)}{c}\right)\right)
- p^{-1}\left(w^{\infty}(s,x) - S^{-1}\left(\frac{\ts^{\infty}(s,x)}{c}\right)\right)
\right| \d x
\\\leq{}&
\int_a^b \left(
\frac{\Rref}{\vref} \left|w^{\infty}(t,x) - w^{\infty}(s,x)\right| + \frac{p^{-1}(\ov)}{\uc \left(\vref \, V'(\ov) - V(\ov)\right)} \left|\ts^{\infty}(t,x) - \ts^{\infty}(s,x)\right|
\right) \d x.
\end{align*}
By \eqref{e:Manowar} we have then
\[\int_a^b \left|\rho^{\infty}(t,x)-\rho^{\infty}(s,x)\right| \d x \leq
\left( \frac{\Rref}{\vref} + \frac{p^{-1}(\ov)}{\uc \left(\vref \, V'(\ov) - V(\ov)\right)} \right) \widetilde{L} \, (t-s),
\]
where $\widetilde{L}$ is given in \eqref{e:tL}.
Similarly, by \eqref{e:Manowar} we get
\begin{align*}
&\int_a^b \left|q^{\infty}(t,x)-q^{\infty}(s,x)\right| \d x = \int_a^b \left|w^{\infty}(t,x) \, \rho^{\infty}(t,x) - w^{\infty}(s,x) \, \rho^{\infty}(s,x)\right| \d x
\\\leq{}&
\|\rho^{\infty}\|_{\L\infty(\D)} \int_a^b \left| w^{\infty}(t,x) -w^{\infty}(s,x) \right| \d x 
+ \|w^{\infty}\|_{\L\infty(\D)} \, \int_a^b \left| \rho^{\infty}(t,x) -\rho^{\infty}(s,x) \right| \d x
\\\leq{}&
\left( \Rref - \uw \left( \frac{\Rref}{\vref} + \frac{p^{-1}(\ov)}{\uc \left(\vref \, V'(\ov) - V(\ov)\right)} \right) \right) \widetilde{L} \, (t-s) .
\end{align*}
Therefore, it is sufficient to take
\[L_{\bU}  \doteq \left(
\Rref + (1 - \uw) \left( \frac{\Rref}{\vref} + \frac{p^{-1}(\ov)}{\uc \left(\vref \, V'(\ov) - V(\ov)\right)} \right)
\right) \widetilde{L}.\]
This concludes the proof.
\end{proof}

\subsection{Consistency of the scheme}

In this section we prove that $\bU^{\infty}$ is an entropy solution by showing first that it is a weak solution.
By \eqref{e:uoncond2}\textsubscript{$2$}, \eqref{e:LipUt} and Lemma~\ref{l:Gojira}, for a.e.\ $x\in\R$ we have the existence of the strong trace $\bU^{\infty}(0^{+},x)$ and that $\bU^{\infty}(t,x) \to \bU_{o}(x)$ as $t\downarrow0$ in $\Lloc1(\R)$.
Hence, the initial condition \eqref{e:ini2} holds almost everywhere.
For this reason we can consider the integral conditions \eqref{e:ws2} and \eqref{e:entro2} for test functions $\varphi$ with compact support in $\Dr$.

\begin{proposition}
$\bU^{\infty}$ is a weak solution of the Cauchy problem \eqref{e:2x2}, \eqref{e:ini2} in the sense of Definition~\ref{d:ws2}.
\end{proposition}

\begin{proof}
We prove that $\bU$ satisfies \eqref{e:ws2}
\begin{equation*}
\iint_{\Dr} \left( \partial_t\varphi + c \, V\left(v(\bU^{\infty})\right) \partial_x\varphi \right) \bU^{\infty} \,\d x\, \d t = 0,
\end{equation*}
for any test function $\varphi \in \Cc\infty(\Dr;\R)$.
Since $\bU^{n}$ is uniformly bounded with $\rho^{n} \geq \ur > 0$, and the flux of system \eqref{e:2x2} is uniformly continuous with respect to $\bU$, uniformly with respect to $x$, on $\calU$, by Lemma~\ref{l:Gojira} it is sufficient to prove that
\begin{align*}
&\lim_{n\to+\infty} \iint_{\Dr} \left( \partial_t\varphi + c \, V\left(v(\bU^{n})\right) \partial_x\varphi \right) \bU^{n} \,\d x\, \d t = 0.
\end{align*}
The double integral above is equal to zero because $\bU^{n}$ is a weak solution, see Corollary~\ref{c:Unweak}.
This concludes the proof.
\end{proof}

Establishing that $\bU^{\infty}$ is an entropy solution is a more involved task.
The proof is based on the following result.

\begin{lemma}
\label{l:StevieRayVaughan}
For any test function $\varphi \in \Cc\infty(\Dr;[0,+\infty))$ and $\kappa \in [0,\ov]$ we have
\begin{equation}
\label{e:StevieRayVaughan}
\liminf_{n\to+\infty} \Biggl[
\begin{aligned}[t]
\iint_{\Dr} \left( \mathfrak{E}_\kappa\left(v(\bU^{n})\right) \partial_t\varphi + \mathfrak{Q}_\kappa\left(v(\bU^{n}),c\right) \partial_x\varphi \right) \d x\, \d t
\\
+ S(\kappa) \sum_{i=0}^{I_{c}} \left|c(\xi_{i}^{+})-c(\xi_{i}^{-})\right| \int_0^{+\infty} \varphi(t,\xi_{i}) \, \d t 
&\Biggr]
\geq 0.
\end{aligned}
\end{equation}

\end{lemma}
\begin{proof}
Choose $T>0$ so that $\varphi$ vanishes outside the strip $(0,T)\times\R$.
Introduce
\[A_{n}(\kappa) \doteq \int_{0}^{T} \sum_{j\in R^{n}(t)} \Upsilon_\kappa\left(\bU_{j-1/2}^{n},\bU_{j+1/2}^{n},\sigma_{i+1/2}^{n}(\bU_{j-1/2}^{n},\bU_{j+1/2}^{n})\right) \varphi\left(t,x_{j}^{n}(t)\right) \d t,
\]
where $R^{n}(t)$ is the set of indexes corresponding to the discretised rarefactions of $\bU^{n}(t)$,  $\Upsilon_\kappa$ is given in \eqref{e:Uk} and $\sigma_{i+1/2}^{n}$ in \eqref{e:Haken}.
By using the Gauss-Green formula and Proposition~\ref{p:InFlames2} we have that
\[
\iint_{\Dr} \left( \mathfrak{E}_\kappa\left(v(\bU^{n})\right) \partial_t\varphi + \mathfrak{Q}_\kappa\left(v(\bU^{n}),c\right) \partial_x\varphi \right) \d x\, \d t
+ S(\kappa) \sum_{i=0}^{I_{c}} \left|c(\xi_{i}^{+})-c(\xi_{i}^{-})\right| \int_0^{+\infty} \varphi(t,\xi_{i}) \, \d t
\geq A_{n}(\kappa).
\]

Fix $j\in R^{n}(t)$ and consider the corresponding discretised rarefaction wave. 
Since its propagation speed is strictly negative, we assume that it belongs to $\mathring\Xi_i$, as this holds for a.e.~$t>0$.
Let $w = w(\bU_{j\pm1/2}^{n})$ and $\rho_{\kappa} = p^{-1}(w-\kappa)$.
If $\kappa \in ( v(\bU_{j-1/2}^{n}) , v(\bU_{j+1/2}^{n}) )$ then
\begin{align*}
0<{}&-\frac{p^{-1}(w)}{c_{i+1/2}} \, \Upsilon_\kappa\left(\bU_{j-1/2}^{n},\bU_{j+1/2}^{n},\sigma_{i+1/2}^{n}(\bU_{j-1/2}^{n},\bU_{j+1/2}^{n})\right) 
\\={}& 
-\frac{F_{w}(\rho_{j+1/2}^{n}) - F_{w}(\rho_{j-1/2}^{n})}{\rho_{j+1/2}^{n}-\rho_{j-1/2}^{n}} \left( \rho_{j+1/2}^{n} + \rho_{j-1/2}^{n} - 2\rho_{\kappa} \right)
- F_{w}(\rho_{j+1/2}^{n}) + 2F_{w}(\rho_{\kappa}) - F_{w}(\rho_{j-1/2}^{n})
\\={}&
2 \left[
F_{w}(\rho_{\kappa}) - F_{w}(\rho_{j-1/2}^{n}) 
- \frac{F_{w}(\rho_{j+1/2}^{n}) - F_{w}(\rho_{j-1/2}^{n})}{\rho_{j+1/2}^{n}-\rho_{j-1/2}^{n}} \left(\rho_{\kappa} - \rho_{j-1/2}^{n}\right)
\right].
\end{align*}
By the concavity of $F_{w}$, the quantity between the square brackets is strictly positive.
Furthermore, the regularity of $F_w$ and the linear interpolation error estimate yield
\begin{align*}
0 &< F_{w}(\rho_{\kappa}) - F_{w}(\rho_{j-1/2}^{n}) 
- \frac{F_{w}(\rho_{j+1/2}^{n}) - F_{w}(\rho_{j-1/2}^{n})}{\rho_{j+1/2}^{n}-\rho_{j-1/2}^{n}} \left(\rho_{\kappa} - \rho_{j-1/2}^{n}\right) 
\\&\leq
\frac{1}{8} \, \left(\rho_{j+1/2}^{n}-\rho_{j-1/2}^{n}\right)^{2} \, \max_{\rho\in[\rho_{j+1/2}^{n}-\rho_{j-1/2}^{n}]} \left|F_{w}''(\rho)\right|
\\&\leq
\frac{\vref}{8 \, \ur} \, \left(\rho_{j+1/2}^{n}-\rho_{j-1/2}^{n}\right)^{2} \, \max_{v\in[0,\ov]} \left(V'(v)-\vref\,V''(v)\right).
\end{align*}
By \eqref{e:SnarkyPuppy}, the monotonicity of $S$ and \eqref{e:LipFm1} we have
\begin{align*}
0 &< \rho_{j-1/2}^{n} - \rho_{j+1/2}^{n} = 
F_{w}^{-1}\left(\frac{\ts(\bU_{j-1/2}^{n},c_{i+1/2})\,p^{-1}(w)}{c_{i+1/2}}\right) - F_{w}^{-1}\left(\frac{\ts(\bU_{j+1/2}^{n},c_{i+1/2})\,p^{-1}(w)}{c_{i+1/2}}\right)
\\&\leq
\frac{\Lip(F_{w}^{-1})}{c_{i+1/2}} \, \left( \ts(\bU_{j+1/2}^{n},c_{i+1/2}) - \ts(\bU_{j-1/2}^{n},c_{i+1/2}) \right) \,p^{-1}(w)
\\&\leq
\frac{p^{-1}(w)}{\left(\vref \, V'(\ov) - V(\ov)\right) c_{i+1/2}} \, \left( \ts(\bU_{j+1/2}^{n},c_{i+1/2}) - \ts(\bU_{j-1/2}^{n},c_{i+1/2}) \right).
\end{align*}
The above estimates ensure that
\begin{align*}
0 &\leq -\Upsilon_\kappa\left(\bU_{j-1/2}^{n},\bU_{j+1/2}^{n},\sigma_{i+1/2}^{n}(\bU_{j-1/2}^{n},\bU_{j+1/2}^{n})\right) 
\\&\leq 
\frac{p^{-1}(w) \, \vref}{4 \ur \left(\vref \, V'(\ov) - V(\ov)\right)^{2} c_{i+1/2}} \, \left( \ts(\bU_{j+1/2}^{n},c_{i+1/2}) - \ts(\bU_{j-1/2}^{n},c_{i+1/2}) \right)^{2} 
\, \max_{v\in[0,\ov]} \left(V'(v)-\vref\,V''(v)\right)
\\&\leq 
C_{\rm ref} \, \varepsilon_{1}^{n} \left( \ts(\bU_{j+1/2}^{n},c_{i+1/2}) - \ts(\bU_{j-1/2}^{n},c_{i+1/2}) \right),
\end{align*}
where
\[C_{\rm ref} \doteq \frac{\Rref \, \vref}{4 \ur \left(\vref \, V'(\ov) - V(\ov)\right)^{2} \uc} \, \max_{v\in[0,\ov]} \left(V'(v)-\vref\,V''(v)\right).\]
The above estimate holds also for $\kappa \notin ( v(\bU_{j-1/2}^{n}) , v(\bU_{j+1/2}^{n}) )$.
Indeed, if $\kappa < v(\bU_{j-1/2}^{n})$, then
\begin{align*}
&\frac{p^{-1}(w)}{c_{i+1/2}} \, \Upsilon_\kappa\left(\bU_{j-1/2}^{n},\bU_{j+1/2}^{n},\sigma_{i+1/2}^{n}(\bU_{j-1/2}^{n},\bU_{j+1/2}^{n})\right) =
\\={}& 
\frac{F_{w}(\rho_{j+1/2}^{n}) - F_{w}(\rho_{j-1/2}^{n})}{\rho_{j+1/2}^{n}-\rho_{j-1/2}^{n}} \left( \rho_{j-1/2}^{n} - \rho_{j+1/2}^{n} \right)
+ F_{w}(\rho_{j+1/2}^{n}) - F_{w}(\rho_{j-1/2}^{n})=0.
\end{align*}
Analogously, we have $\Upsilon_\kappa\left(\bU_{j-1/2}^{n},\bU_{j+1/2}^{n},\sigma_{i+1/2}^{n}(\bU_{j-1/2}^{n},\bU_{j+1/2}^{n})\right) = 0$ also in the case $\kappa \geq v(\bU_{j+1/2}^{n})$.

By Theorem~\ref{t:Deftones} and \eqref{e:uoncond2}\textsubscript{$1$}, we have $\TV(\ts(\bU^{n}(t),c))$ is bounded by $\tTV(\bU_o,c)$ uniformly in $n$ and $t$, and therefore
\begin{align*}
0 &\leq -A_{n}(\kappa) \leq 
-C_{\rm ref} \, \varepsilon_{1}^{n} \int_{0}^{T} \sum_{j\in R^{n}(t)} \left( \ts(\bU_{j+1/2}^{n},c_{i+1/2}) - \ts(\bU_{j-1/2}^{n},c_{i+1/2}) \right)
 \varphi\left(t,x_{j}^{n}(t)\right) \d t
 \\&\leq
-C_{\rm ref} \, \|\varphi\|_{\L\infty(\R)} \, T \, \tTV(\bU_o,c) \, \varepsilon_{1}^{n} \xrightarrow[n\to 0]{} 0,
\end{align*}
and this concludes the proof.
\end{proof}

\begin{proposition}
\label{p:Uentro}
$\bU^{\infty}$ is an entropy solution of \eqref{e:2x2}, \eqref{e:ini2} in the sense of Definition~\ref{d:entro2}.
\end{proposition}
\begin{proof}
The proof of the entropy inequality \eqref{e:entro2} relies on Lemma~\ref{l:StevieRayVaughan}.
The a.e.\ convergence of $\{\bU^{n}\}_{n}$ to $\bU^{\infty}$, ensures that the first line in \eqref{e:StevieRayVaughan} yields the first line in \eqref{e:entro2} for any $\varphi \in \Cc\infty(\Dr;[0,+\infty))$.
The second line of \eqref{e:entro2} does not involve $\bU^{n}$ and coincides with the first line in \eqref{e:entro2}.
This concludes the proof.
\end{proof}

\subsection{Uniqueness}

As already proved in Proposition~\ref{p:Uentro}, $\bU^{\infty}$ is an entropy solution.
In particular, it is a weak solution of the Cauchy problem
\begin{equation}
\left\{\begin{array}{@{}l@{}}
\partial_t \rho + \partial_x \left( b \, \rho \right)=0
\\
\rho(0,x)=\rho_{o}(x),
\end{array}\right.
\label{e:transport}
\end{equation}
where $b \doteq c \, V(v(\bU^{\infty}))$ is in $\L\infty((0,T)\times\R;\R)$ and $\rho_{o} \in \Lloc1(\R)$.
Recall that $\rho^{\infty}$ takes values in $[\ur,\Rref]$, hence $b$ is nearly incompressible in the sense of \cite[Definition~1.1]{Gusev}.
As a consequence, the Cauchy problem \eqref{e:transport} admits a unique weak solution by \cite[Theorem~1.2]{Gusev}.

Assume that also $\widetilde{\bU}^{\infty}$ is an entropy solution of the Cauchy problem \eqref{e:2x2}, \eqref{e:ini2} with the same initial datum $\bU_{o}=(\rho_{o},q_{o})$.
By Theorem~\ref{t:uniqueness} we have that $v(\bU^{\infty})$ and $v(\widetilde{\bU}^{\infty})$ coincide.
As a result, both $\rho^{\infty}$ and $\widetilde{\rho}^{\infty}$ are weak solutions of \eqref{e:transport} with the same $b$ and initial datum $\rho_{o}$.
This implies that $\rho^{\infty} \equiv \widetilde{\rho}^{\infty}$ and, consequently $w(\bU^{\infty}) \equiv w(\widetilde{\bU}^{\infty})$ (because $w=v+p(\rho)$) and therefore $q^{\infty} \equiv \widetilde{q}^{\infty}$ (because $q=w\,\rho$).

\subsection*{Acknowledgment}
F.A.C., S.F.\ and M.D.R.\ are members of Gruppo Nazionale per l’Analisi Matematica, la Probabilit\`a e le loro Applicazioni (GNAMPA) of the Istituto Nazionale di Alta Matematica (INdAM). 
F.A.C. and S.F. are partially supported by the Ministry of University and Research (MUR), Italy, under the grant PRIN 2020 - Project N.~20204NT8W4, \lq\lq Non-linear evolution PDEs, fluid dynamics and transport equations: theoretical foundations and applications\rq\rq.
F.A.C.~recognizes financial support from the Project \lq\lq Leggi di conservazione con termini nonlocali e applicazioni al traffico veicolare\rq\rq, Progetti di Ateneo 2025, University of L'Aquila.
S.F. is partially supported by the Italian \lq\lq National Centre for HPC, Big Data and Quantum Computing\rq\rq\ - Spoke 5 \lq\lq Environment and Natural Disasters\rq\rq, by the InterMaths Network, \url{www.intermaths.eu}, by the INdAM-GNAMPA project 2026 code CUP E53C25002010001 \lq\lq Modelli di reazione-diffusione-trasporto: dall'analisi alle applicazioni\rq\rq and the INdAM-GNAMPA project 2025 code CUP E5324001950001 \lq\lq Teoria e applicazioni dei modelli evolutivi: trasporto ottimo, metodi variazionali e approssimazioni particellari deterministiche\rq\rq. 
M.D.R.~recognizes financial support from the PRIN 2022 Project \lq\lq Modeling, Control and Games through Partial Differential Equations\rq\rq, CUP~D53D23005620006, funded by the European Union-Next Generation EU, and from \lq\lq INdAM - GNAMPA Project\rq\rq, CUP~E5324001950001.

\bibliographystyle{abbrv}
\bibliography{biblio}

\end{document}